%% file: LNPnotes.tex
\documentclass{svmult}

\input{LNPpreamble}

\begin{document}

\title{Introduction to Categories and Categorical Logic}
\author{Samson Abramsky and Nikos Tzevelekos}
\institute{Oxford University Computing Laboratory\\Wolfson Building, Parks Road, Oxford OX1 3QD, U.K.}
\maketitle

\preface
The aim of these notes is to provide a succinct, accessible introduction to some of the basic ideas of category theory and categorical logic. The notes are based on a lecture course given at Oxford over the past few years. They contain numerous exercises, and hopefully will prove useful for self-study by those seeking a first introduction to the subject, with fairly minimal prerequisites. The coverage is by no means comprehensive, but should provide a good basis for further study; a  guide to further reading is included.

The main prerequisite is a basic familiarity with the elements of discrete mathematics: sets, relations and functions. An Appendix contains a summary of what we will need, and it may be useful to review this first.
In addition, some prior exposure to abstract algebra\HY vector spaces and linear maps, or groups and group homomorphisms\HY would be helpful.

\setcounter{tocdepth}{2}
\tableofcontents\listoftables

\section{Introduction}
Why study categories\,---\,what are they good for? We can offer a range of answers for readers coming from different backgrounds:
\begin{itemize}
\item For \textbf{mathematicians}: category theory organises your previous mathematical experience in a new and powerful way, revealing new connections and structure,
and allows you to ``think bigger thoughts''.
\item For  \textbf{computer scientists}: category theory gives a precise handle on important notions such as compositionality, abstraction, represen{\-}tation-independence, genericity and more.
Otherwise put, it provides the fundamental mathematical structures underpinning many key programming concepts.
\item For  \textbf{logicians}: category theory gives a syntax-independent view of the fundamental structures of logic, and opens up new kinds of models and interpretations.
\item For  \textbf{philosophers}: category theory opens up a fresh approach to structuralist foundations of mathematics and science; and an alternative to the traditional focus
on set theory.
\item For \textbf{physicists}: category theory offers new ways of formulating physical theories in a structural form. There have \textit{inter alia}~been some striking recent applications to quantum information and
computation.
\end{itemize}

\subsection{From Elements To Arrows}
Category theory can be seen as a ``generalised theory of functions'',  where the focus is shifted from the pointwise, set-theoretic view of functions,
to an abstract view of functions as \emph{arrows}.

Let us briefly recall the arrow notation for functions between sets.\footnote{A review of basic ideas about sets, functions and relations, and some of the notation we will be using, is provided in Appendix A.} A
function $f$ with \emph{domain} $X$ and \emph{codomain} $Y$ is denoted by: $f : X \rarr Y$.
\[ \mbox{Diagrammatic notation:} \;\;  X \stackrel{f}{\longrightarrow} Y\,.\]
The fundamental operation on functions is \emph{composition}: if $f : X \rarr Y$ and $g : Y \rarr Z$, then we can define $g \circ f: X \rarr Z$ by $g \circ f(x) := g(f(x))$.\footnote{We shall use the notation ``:=" for ``is defined to be" throughout these notes.} Note that, in order for the composition to be defined, the codomain of $f$ must be the same as the domain of $g$.
\[  \mbox{Diagrammatic notation:}\;\; X \stackrel{f}{\longrightarrow} Y \stackrel{g}{\longrightarrow} Z\,. \]
Moreover, for each set $X$ there is an \emph{identity function} on $X$, which is denoted by:
\[ \id[X] : X \longrightarrow X\qquad \id[X](x) := x \,. \]
These operations are  governed by the \emph{associativity law} and the \emph{unit laws}. For $f : X \rarr Y$, $g: Y \rarr Z$, $h : Z \rarr W$:
\[ (h \circ g) \circ f = h \circ (g \circ f)\,, \qquad f\circ\id[X] = f = \id[Y]\circ f\,. \]
Notice that these equations are formulated purely in terms of the algebraic operations on functions, without any reference to the elements of the sets $X$, $Y$, $Z$, $W$.
We will refer to any concept pertaining to functions which can be defined purely in terms of composition and identities  as \emph{arrow-theoretic}.
We will now take a first step towards learning to ``think with arrows'' by seeing how we can replace some familiar definitions couched in terms of elements by arrow-theoretic equivalents; this will lead us towards the notion of category.

We say that a function $f : X \longrightarrow Y$ is:
\begin{center}
\begin{tabular}{lll}
\emph{injective} & if &
$\forall x,x'\in X.\;f(x) = f(x') \;\; \Longrightarrow \;\; x = x'$\,, \\
\emph{surjective} & if &
$\forall y \in Y. \, \exists x \in X. \, f(x) = y$\,, \\
\\
\emph{monic} & if &
$\forall g,h.\;f \circ g = f \circ h \;\; \Longrightarrow \;\; g = h$\,, \\
\emph{epic} & if & $\forall g,h.\;g \circ f = h \circ f \;\; \Longrightarrow \;\; g = h$\,.
\end{tabular}
\end{center}
Note that injectivity and surjectivity are formulated in terms of elements, while epic and monic are arrow-theoretic.
\begin{myproposition} Let $f:X\rarr Y$. Then,
\begin{enumerate}\renewcommand{\theenumi}{{\rm\arabic{enumi}}}
\item $f$ is injective iff $f$ is monic.
\item $f$ is surjective iff $f$ is epic.
\end{enumerate}
\end{myproposition}
\proof We show 1. Suppose $f : X \rarr Y$ is injective, and that $f \circ g = f \circ h$, where $g, h : Z \rarr X$. Then, for all $z \in Z$:
\[ f(g(z)) = f \circ g(z) = f \circ h(z) = f(h(z))\,. \]
Since $f$ is injective, this implies $g(z) = h(z)$. Hence we have shown that
\[ \forall z \in Z. \; g(z) = h(z)\,, \]
and so we can conclude that $g = h$. So $f$ injective implies $f$ monic.
\\
For the converse, fix a one-element set $\Term = \{ \bullet \}$. Note that elements $x \in X$ are in 1--1 correspondence with functions $\bar{x} :
\Term \rarr X$, where $\bar{x}(\bullet) := x$. Moreover, if $f(x) = y$ then $\bar{y} = f \circ \bar{x}$\,. Writing injectivity in these terms, it
amounts to the following.
\[ \forall x, x' \in X. \; f \circ \bar{x} = f \circ \bar{x}' \;\; \Longrightarrow \;\; \bar{x} = \bar{x}'  \]
Thus we see that being injective is a \emph{special case} of being monic. \qed
\begin{myexercise}
Show that $f:X\rarr Y$ is surjective iff it is epic.
\end{myexercise}

%

\subsection{Categories Defined}
\begin{mydefinition}
A \boldemph{category} $\CC$ consists of:
\begin{itemize}
  \item A collection $\Ob(\CC)$ of \boldemph{objects}. Objects are denoted by $A$, $B$, $C$, etc.
  \item A collection $\Arr(\CC)$ of \boldemph{arrows} (or \boldemph{morphisms}). Arrows are denoted by $f$, $g$, $h$, etc.
  \item Mappings $\dm, \cod : \Arr(\CC) \rarr \Ob(\CC)$, which assign to each arrow $f$ its \boldemph{domain} $\dm(f)$ and its \boldemph{codomain} $\cod(f)$. An arrow $f$ with domain $A$ and codomain $B$  is written
  $f:A\rarr B$.  For each pair of objects $A$, $B$, we define the  set
  \[ \CC(A,B) := \{ f \in \Arr(\CC) \mid f : A \rarr B \} \, . \]
  We refer to $\CC(A, B)$ as a \boldemph{hom-set}. Note that distinct hom-sets are \emph{disjoint}.
    \item For any triple of objects $A$, $B$, $C$, a \boldemph{composition} map
    \[ c_{A,B,C} : \CC (A, B) \times \CC (B, C) \longrightarrow \CC (A, C)\,. \]
  $c_{A, B, C}(f, g)$ is written \ $g\circ\!f$ \ (or sometimes $f ; g$). Diagrammatically:
  \[ A \stackrel{f}{\longrightarrow} B \stackrel{g}{\longrightarrow} C \]
  \item For each object $A$, an \boldemph{identity} arrow $\id[A]:A\rightarrow A$.
\end{itemize}
The above must satisfy the following axioms.
\[ h \circ (g \circ f) = (h \circ g) \circ f \,,\qquad f \circ \id[A] = f = \id[B] \circ f\,. \]
whenever the domains and codomains of the arrows match appropriately so that the compositions are well-defined. \deq
\end{mydefinition}

\subsection{Diagrams in Categories}
\boldemph{Diagrammatic reasoning} is an important tool in category theory. The basic cases are commuting triangles and squares. To say that the following triangle commutes
\[
\xymatrix@=13mm{A\ar[r]^f\ar[dr]_h & B\ar[d]^g\\ & C}
\]
is exactly equivalent to asserting the equation $g \circ f = h$. Similarly, to say that the following square commutes
\[
\xymatrix@=13mm{A\ar[r]^f\ar[d]_h & B\ar[d]^g\\ C\ar[r]_k & D}
\]
means exactly that $g \circ f = k \circ h$. For example, the equations
\[ h \circ (g \circ f) = (h \circ g) \circ f \,,\qquad f \circ \id[A] = f = \id[B] \circ f\,, \]
can be expressed by saying that the following diagrams commute.
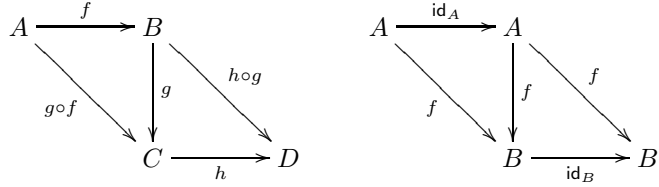
\[ \xymatrix@=13mm{A \ar[r]^f\ar[dr]_{g\circ f} & B\ar[d]^g\ar[dr]^{h\circ g} \\ & C\ar[r]_h & D} \qquad
   \xymatrix@=13mm{A \ar[r]^{\id[A]}\ar[dr]_{f} & A\ar[d]^f\ar[dr]^{f} \\ & B\ar[r]_{\id[B]} & B}
\]
As these examples illustrate, most of the diagrams we shall use will be ``pasted together'' from triangles and squares: the commutation of the diagram as a whole will then reduce to the commutation of the constituent triangles and squares.

We turn to the general case. The formal definition is slightly cumbersome; we give it anyway for reference.

\begin{mydefinition}
\label{diagdef}
We define a \boldemph{graph} to be a collection of \boldemph{vertices} and \boldemph{directed edges}, where each edge $e : v \rarr w$ has a specified source vertex $v$ and target vertex $w$.
Thus graphs are like categories without composition and identities.\footnote{This would be a ``multigraph'' in normal parlance, since multiple edges between a given pair of vertices are allowed.}
A \boldemph{diagram in a category $\CC$} is a graph  whose vertices are labelled with objects of $\CC$ and whose edges are labelled with arrows of $\CC$, such that, if $e : v \rarr w$  is labelled with $f : A \rarr B$, then we must have $v$ labelled by $A$ and $w$ labelled by $B$.  We say that such a  diagram  \emph{commutes} if  any two paths in it with common source and target, and at least one of which has length greater than 1, are equal. That is, given paths
\[ A \labarrow{f_{1}} C_{1} \labarrow{f_{2}} \cdots C_{n-1} \labarrow{f_{n}} B  \qquad \mbox{and} \qquad
A \labarrow{g_{1}} D_{1} \labarrow{g_{2}} \cdots D_{m-1} \labarrow{g_{m}} B,
\]
if $ \max(n,m) > 1$ then
\[ f_{n} \circ \cdots \circ f_{1} = g_{m} \circ \cdots \circ g_{1} \, . \]
 \deq[-1]
\end{mydefinition}

\noindent To illustrate this definition, to say that the following diagram commutes
\[
\xymatrix@=13mm{E\ar[r]^e & A\ar@<2pt>[r]^f\ar@<-2pt>[r]_g & B}
\]
amounts to the assertion that $f \circ e = g \circ e$\,; it does \emph{not} imply that $f = g$.

\subsection{Examples}
Before we proceed to our first examples of categories, we shall present some background material on partial orders, monoids and topologies, which will provide running examples throughout these notes.

\paragraph{Partial orders}
A \emph{partial order} is a structure $(P, \leq )$ where $P$ is a set and $\leq$ is a binary relation on $P$ satisfying:
\begin{itemize}
\item[$\bullet$] $x \leq x$ \hfill(Reflexivity)
\item[$\bullet$] $x \leq y \; \wedge y \leq x \;\; \Rightarrow \;\; x = y$  \hfill(Antisymmetry)
\item[$\bullet$] $x \leq y \; \wedge \; y \leq z \;\; \Rightarrow \;\; x \leq z$  \hfill(Transitivity)
\end{itemize}
For example, $(\mathbb{R}, \leq )$ and $(\mathcal{P}(X), \subseteq )$ are partial orders, and so are strings with the sub-string relation.

If $P$, $Q$ are partial orders, a map $h : P \rarr Q$ is a \emph{partial order homomorphism} (or \emph{monotone function}) if:
\[ \forall x,y\in P.\; x \leq y \; \Longrightarrow \; h(x) \leq h(y)\,. \]
Note that homomorphisms are {closed under composition}, and that {identity maps} are homomorphisms.

\paragraph{Monoids} A \emph{monoid} is a structure $(M, \cdot , 1)$ where $M$ is a set,
\[ \uscore\cdot\uscore : M \times M \longrightarrow M \]
is a binary operation, and $1 \in M$, satisfying the following axioms.
\[ (x \cdot y) \cdot z = x \cdot (y \cdot z)\,, \qquad 1 \cdot x = x = x\cdot 1\,. \]
For example, $(\mathbb{N}, + , 0)$ is a monoid, and so are strings with string-concatenation. Moreover, groups are special kinds of monoids.

If $M$, $N$ are monoids, a map $h:M\rarr N$ is a \emph{monoid homomorphism} if
\[ \forall m_1,m_2\in M.\;h(m_1\cdot m_2) = h(m_1)\cdot h(m_2)\,, \qquad h(1) = 1\,.  \]

\begin{myexercise}
\label{monhom}
Suppose that $G$ and $H$ are groups (and hence monoids), and that $h : G \rarr H$ is a monoid homomorphism. Prove that $h$ is a group homomorphism.
\end{myexercise}

\paragraph{Topological spaces}
A topological space is a pair $(X, T_{X})$ where $X$ is a set, and $T_{X}$ is a family of subsets of $X$ such that
\begin{itemize}
  \item $\varnothing,X \in T_{X}$,
  \item if $U, V \in T_{X}$ then $U \cap V \in T_{X}$,
  \item if $\{ U_{i} \}_{i \in I}$ is any family in $T_{X}$, then $\bigcup_{i \in I} U_{i} \in T_{X}$\,.
\end{itemize}
A \emph{continuous map} $f : (X, T_{X}) \rarr (Y, T_{Y})$ is a function $f : X \rarr Y$ such that, for all $U \in T_{Y}$, $f^{-1}(U) \in T_{X}$.
\newlines{2}
Let us now see some first examples of  categories.
\begin{itemize}
\item {Any kind of mathematical structure, together with structure preserving functions, forms a category. E.g.}
\begin{itemize}
\item {\textbf{Set} (sets and functions)}
\item {\textbf{Mon} (monoids and monoid homomorphisms)}
\item {\textbf{Grp} (groups and group homomorphisms)}
\item {$\text{\textbf{Vect}}_{k}$ (vector spaces over a field $k$, and linear maps)}
\item {\textbf{Pos} (partially ordered sets and monotone functions)}
\item {\textbf{Top} (topological spaces and continuous functions)}
\end{itemize}
\item \Rel: objects are sets, arrows $R : X \rightarrow Y$ are \emph{relations} $R \subseteq X \times Y$.
Relational composition:
\[ R;S (x, z) \;\; \Longleftrightarrow \;\; \exists y. \, R(x, y) \; \wedge \; S(y, z) \]
\item Let $k$ be a field (for example, the real or complex numbers). Consider the following category $\text{\textbf{Mat}}_{k}$. The objects are natural numbers. A morphism $M : \mathbf{n} \rarr \mathbf{m}$ is an $\mathbf{n} \times \mathbf{m}$ matrix with entries in $k$. Composition is matrix multiplication, and the identity on $\mathbf{n}$ is the $\mathbf{n} \times \mathbf{n}$ diagonal matrix.
\item[$\diamond$] Monoids are one-object categories. Arrows correspond to the elements of the monoid, with the monoid operation being arrow-composition
    and the monoid unit being the identity arrow.
\item[$\diamond$] {A category in which for each pair of objects $A$, $B$ there is at most one morphism from $A$ to $B$ is the same thing as a
  \boldemph{preorder}, \ie a reflexive and transitive relation.}
\end{itemize}
Note that our first class of examples illustrate the idea of categories as \emph{mathematical contexts}; settings in which various mathematical theories can be developed. Thus for example, \tbf{Top} is the context for general topology, \tbf{Grp} is the context for group theory, etc.

On the other hand, the last two examples illustrate that many important mathematical structures themselves appear as categories of particular kinds. The fact that two such different kinds of structures as monoids and posets should appear as extremal versions of categories is also rather striking.

This ability to capture mathematics both ``in the large'' and ``in the small'' is a first indication of the flexibility and power of categories.

\begin{myexercise}
Check that \Mon, $\Vect_k$, \Pos\ and \Top\ are indeed categories.
\end{myexercise}
\begin{myexercise}
Check carefully that  monoids correspond exactly to one-object categories. Make sure you understand the difference between such a category and \Mon. (For example: how many objects does \Mon\  have?)
\end{myexercise}
\begin{myexercise}
Check carefully that  preorders correspond exactly to categories in which each homset has at most one element. Make sure you understand the difference between such a category and \Pos. (For example: how big can homsets in \Pos\  be?)
\end{myexercise}
\subsection{First Notions}
Many important mathematical notions can be expressed at the general level of categories.
\begin{mydefinition}
Let $\CC$ be a category. A morphism $f : X \rarr Y$ in $\CC$ is:\\[1.5mm]
\begin{tabular}{lll}
$\bullet$ \ \boldemph{monic} (or a  \boldemph{monomorphism}) & if &  $f \circ g = f \circ h \; \Longrightarrow \; g = h$\,, \\[1.5mm]
$\bullet$ \ \boldemph{epic}  (or an \boldemph{epimorphism})  & if &  $g \circ f = h \circ f \; \Longrightarrow \; g = h$\,.
\end{tabular}\\[1.5mm]
An \boldemph{isomorphism} in $\CC$ is an arrow $i:A\rarr B$ such that there exists an arrow $j : B \rarr A$ --- the \boldemph{inverse} of $i$ --- satisfying
\[ j \circ i = \id[A]\,, \qquad i \circ j = \id[B]\,. \]
\deq[-1]
\end{mydefinition}
We denote isomorphisms by $i:A\xrightarrow{\cong}B$, and write $i^{-1}$ for the inverse of $i$. We say that $A$ and $B$ \boldemph{are
isomorphic}, $A\cong B$, if there exists some $i:A\xrightarrow{\cong}B$.

\begin{myexercise}
Show that the inverse, if it exists, is unique.
\end{myexercise}
\begin{myexercise}
Show that $\cong$ is an equivalence relation on the objects of a category.
\end{myexercise}
As we saw previously, in $\Set$ monics are injections and epics are surjections. On the other hand, isomorphisms in \Set\ correspond exactly to bijections, in \Grp\
to group isomorphisms, in \Top\ to homeomorphisms, in \Pos\ to order isomorphisms, etc.
\begin{myexercise}
Verify these claims.
\end{myexercise}
Thus we have at one stroke captured the key notion of isomorphism in a form which applies
to \emph{all} mathematical contexts. This is a first taste of the level of generality which category theory naturally affords.

We have already identified monoids as one-object categories. We can now identify groups as \emph{exactly those one-object categories in which every arrow is an isomorphism}. This also leads to a natural generalisation, of considerable importance in current mathematics: a \boldemph{groupoid} is a category in which every morphism is an isomorphism.

\paragraph{Opposite Categories and Duality}
The directionality of arrows within a category $\CC$ can be reversed without breaking the conditions of being a category; this yields the notion of \boldemph{opposite category}.

\begin{mydefinition}
Given a category $\CC$, \boldemph{the opposite category} $\CCop$ is given by taking the same objects as $\CC$, and
\[ \CCop(A, B) := \mathcal{C}(B, A)\,. \]
Composition and identities are inherited from $\CC$. \deq
\end{mydefinition}
Note that if we have
\[ A \stackrel{f}{\longrightarrow} B \stackrel{g}{\longrightarrow} C
\]
in $\CCop$, this means
\[ A \stackrel{f}{\longleftarrow} B \stackrel{g}{\longleftarrow} C \]
in $\CC$, so composition  $g \circ f$ in $\CCop$ is defined as $f \circ g$ in $\CC$!

Consideration of opposite categories leads to a \boldemph{principle of duality}: a statement $S$ is true about $\CC$ if and only if its dual (\ie the one
obtained from $S$ by reversing all the arrows) is true about $\op{\CC}$. For example,
\[ \text{A morphism $f$ is {monic} in $\mathcal{C}^{\mathsf{op}}$ if and only if it is {epic} in $\mathcal{C}$\,.} \]
Indeed, $f$ is monic in $\CCop$ iff for all $g, h : C \rarr B$ in $\CCop$,
\[ f \circ g = f \circ h \; \Longrightarrow \; g = h\,, \]
iff for all $g, h : B \rarr C$ in $\CC$,
\[ g \circ f = h \circ f \; \Longrightarrow \; g = h\,, \]
iff $f$ is epic in $\CC$. We say that monic and epic are \emph{dual notions}.

\begin{myexercise}
If $P$ is a preorder, for example $(\mathbb{R}, \leq )$, describe $\op{P}$ explicitly.
\end{myexercise}

\paragraph{Subcategories} Another way to obtain new categories from old ones is by restricting their objects or arrows.

\begin{mydefinition}
Let $\CC$ be a category. Suppose that we are given collections
\[ \Ob(\DD)\subseteq \Ob(\CC)\,,\quad\forall A,B\in \Ob(\DD).\; \DD(A,B)\subseteq\CC(A,B)\,. \]
We say that $\DD$ is a \boldemph{subcategory} of $\CC$ if
\[ A \in \Ob(\DD) \; \Rightarrow \; \id[A] \in \DD(A, A), \quad f \in \DD(A, B), g  \in \DD(B, C) \; \Rightarrow \; g \circ f \in \DD(A, C) \, , \]
and hence $\DD$ itself is a category.
In particular, $\DD$ is:
\begin{itemize}
  \item A \boldemph{full} subcategory of $\CC$ if  for any $A,B\in \Ob(\DD)$, $\DD(A,B)=\CC(A,B)$.
  \item A \boldemph{lluf} subcategory of $\CC$ if $\Ob(\DD)=\Ob(\CC)$. \deq
\end{itemize}
\end{mydefinition}
For example, \Grp\ is a full subcategory of \Mon\ (by Exercise~\ref{monhom}), and \Set\ is a lluf subcategory of \Rel.

\paragraph{Simple cats} We close this section with some very basic examples of categories.
\begin{itemize}
  \item $\1$ is the category with one object and one arrow, that is,
  \begin{gather*}
    \1:= \ \xymatrix@M=1pt{\bullet\ar@(dr,dl)[]}\\[-1mm]
  \end{gather*}
  where the arrow is necessarily $\id[\bullet]$\,. Note that, although we say that $\1$ is \emph{the} one-object/one-arrow category, there is by no
  means a \emph{unique} such category. This is explained by the intuitively evident fact that any two such categories are isomorphic. (We will define what it means for categories to be isomorphic later.)
  \item In two-object categories, there is the one with two arrows, $\2:= \ \xymatrix@=2mm{\bullet & \bullet}$\,, and also:
  \[ \2_{\scriptscriptstyle\rarr}:= \ \xymatrix@=8mm@M=2pt{\bullet\ar[r] & \bullet}\quad,
  \quad \2_{\scriptscriptstyle\rightrightarrows}:= \ \xymatrix@=10mm@M=2pt{\bullet\ar@/^/[r]\ar@/_/[r] & \bullet}\quad,
  \quad \2_{\scriptscriptstyle\rightleftarrows}:= \ \xymatrix@=10mm@M=2pt{\bullet\ar@/^/[r] & \bullet\ar@/^/[l]}\quad\dots \]
  Note that we have omitted identity arrows for economy.   Categories with only identity arrows, like $\1$ and $\2$, are called \emph{discrete categories}.
\end{itemize}
\begin{myexercise}
How many categories $\CC$ with $\Ob(\CC)=\{\bullet\}$ are there? (Hint: what do such categories correspond to?)
\end{myexercise}

\subsection{Exercises}
\begin{enumerate}\renewcommand{\theenumi}{\textbf{\arabic{enumi}}}
    \item Consider the following properties of an arrow $f$ in a category $\CC$.
    \begin{itemize}
        \item $f$ is \emph{split monic} if for some $g$, $g \circ f$ is an  identity arrow.
        \item $f$ is \emph{split epic} if for some $g$, $f \circ g$ is an  identity arrow.
    \end{itemize}
    \begin{enumerate}
        \item Prove that if $f$ and $g$ are arrows such that $g \circ f$ is monic, then $f$ is monic.
        \item Prove that, if $f$ is split epic then it is epic.
        \item Prove that, if $f$ and $g \circ f$ are iso then $g$ is iso.
        \item Prove that, if $f$ is monic and split epic then it is iso.
        \item In the category $\mathbf{Mon}$ of monoids and monoid homomorphisms, consider the inclusion map
        \[ i : (\mathbb{N},{+}, 0) \longrightarrow (\mathbb{Z}, {+}, 0) \]
        of natural numbers into the integers. Show that this arrow is both monic and epic. Is it an iso?
    \end{enumerate}
      The \textbf{Axiom of Choice} in Set Theory states that, if $\{ X_{i} \}_{i \in I}$ is a family of non-empty sets, we can form a set
      $X=\{ x_{i}\mid i\in I \}$ where $x_{i} \in X_{i}$ for all $i \in I$.
    \begin{enumerate}\setcounter{enumii}{5}
      \item Show that in $\mathbf{Set}$ an arrow which is epic is split epic. Explain why this needs the Axiom of Choice.
      \item Is it always the case that an arrow which is epic is split epic? Either prove that it is, or give a counter-example.
    \end{enumerate}
    \item Give a description of partial orders as categories of a special kind.
\end{enumerate}

\section{Some Basic Constructions}\label{s:BasConstr}

We shall now look at a number of basic constructions which appear throughout mathematics, and which acquire their proper general form in the language of categories.

\subsection{Initial and Terminal Objects}

A first such example is that of initial and terminal objects. While apparently trivial, they are actually both important and useful, as we shall see in the sequel.

\begin{mydefinition}
An object $I$ in a category $\CC$ is \boldemph{initial} if, for every object $A$, \emph{there exists a unique arrow} from $I$ to $A$, which we write
$\iota_A : I \rarr A$. \\
A \boldemph{terminal} object in $\CC$ is an object $T$ such that, for every object $A$, \emph{there exists a unique arrow} from $A$ to $T$, which we write $\tau_A: A \rarr T$.\deq
\end{mydefinition}
Note that initial and terminal objects are dual notions: $T$ is terminal in $\CC$ iff it is initial in $\CCop$. We sometimes write $\ena$ for the
terminal object and $\zero$ for the initial one.
Note also the assertions of \emph{unique existence} in the definitions. This is one of the \textit{leitmotifs} of category theory; we shall encounter it again in a conceptually deeper form in section~\ref{s:Univ}.

Let us examine initial and terminal objects in our standard example categories.
\begin{itemize}
\item In \textbf{Set}, the empty set is an initial object while any one-element set $\{\bullet\}$ is terminal.
\item In \textbf{Pos}, the poset $(\varnothing,\varnothing)$ is an initial object while $(\{\bullet\}, \{ (\bullet,\bullet)\})$ is terminal.
\item In \textbf{Top}, the space $(\varnothing,\{\varnothing\})$ is an initial object while $(\{\bullet\},\{\varnothing,\{\bullet\}\})$ is terminal.
\item In $\text{\textbf{Vect}}_k$, the one-element space $\{0\}$ is both initial and terminal.
\item In a poset, seen as a category, an initial object is a least element, while a terminal object is a greatest element.
\end{itemize}
\begin{myexercise}
Verify these claims.
In each case, identify the canonical arrows.
\end{myexercise}
\begin{myexercise}
Identify the initial and terminal objects in \tbf{Rel}.
\end{myexercise}
\begin{myexercise}
Suppose that a monoid, viewed as a category, has either an initial or a terminal object. What must the monoid be?
\end{myexercise}
We shall now establish a fundamental fact: initial and terminal objects are \emph{unique up to (unique) isomorphism}. As we shall see, this is characteristic of all such ``universal'' definitions.
For example, the apparent arbitrariness in the fact that any singleton set is a terminal object in \tbf{Set} is answered by the fact that what counts is the property of being terminal; and this suffices to ensure that any two concrete objects having this property must be isomorphic to each other.

The proof of the proposition, while elementary, is a first example of distinctively categorical reasoning.
\begin{myproposition}
\label{inuniqueprop}
  If $I$ and $I'$ are initial objects in the category $\CC$ then there exists a unique isomorphism $I\iso I'$.
\end{myproposition}
\proof%
Since $I$ is initial and $I'$ is an object of $\CC$, there is a unique arrow $\iota_{I'}:I\rarr I'$. We claim that $\iota_{I'}$ is  an isomorphism.
\\
Since $I'$ is initial and $I$ is an object in $\CC$, there is an arrow $\iota'_I:I'\rarr I$. Thus we obtain $\iota_{I'};\iota'_I:I\rarr I$, while we
also have the identity morphism $\id[I]:I\rarr I$. But $I$ is initial and therefore there exists a \emph{unique} arrow from $I$ to $I$, which means that
$\iota_{I'};\iota'_I=\id[I]$. Similarly, $\iota'_{I};\iota_{I'}=\id[I']$, so $\iota_{I'}$ is indeed an isomorphism. \qed[2]
Hence, initial objects are ``unique up to (unique) isomorphism'', and we can (and do) speak of \emph{the} initial object (if any such exists). Similarly for terminal objects.

\begin{myexercise}
Let $\CC$ be a category with an initial object $\zero$. For any object $A$, show the following.
\begin{itemize}
  \item If $A\cong\zero$ then $A$ is an initial object.
  \item If there exists a monomorphism $f:A\rarr\zero$ then $f$ is an iso, and hence $A$ is initial.
\end{itemize}
\end{myexercise}

\subsection{Products and Coproducts}

\subsubsection{Products}
We now consider one of the most common constructions in mathematics: the formation of ``direct products''. Once again, rather than giving a case-by-case construction of direct products in each mathematical context we encounter, we can express once and for all a general notion of product, meaningful in any category --- and such that, if a product exists, it is characterised uniquely up to unique isomorphism, just as for initial and terminal objects. Given a particular mathematical context, \ie a category, we can then verify whether on not the product exists in that category. The concrete construction appropriate to the context will enter only into the proof of \emph{existence}; all of the useful \emph{properties} of the product follow from the general definition. Moreover, the categorical notion of product has a \emph{normative} force; we can test whether a concrete construction works as intended by verifying that it satisfies the general definition.

In set theory, the cartesian product is defined in terms of the ordered pair:
\[ X \times Y := \{ (x, y) \mid x \in X \; \wedge \; y \in Y \} . \]
It turns out that ordered pairs can be \emph{defined} in set theory, e.g.~as
\[ ( x, y ) := \{ \{ x, y \}, y\} . \]
Note that in no sense is such a definition canonical. The essential \emph{properties} of ordered pairs are:
\begin{enumerate}
\item We can retrieve the first and second components $x$, $y$ of the ordered pair $(x, y )$, allowing \emph{projection functions} to be defined:
\[ \pi_{1} : (x, y) \mapsto x, \qquad \pi_{2} : (x, y) \mapsto y \, . \]
\item The information about first and second components completely determines the ordered pair:
\[ (x_{1}, x_{2} ) =  ( y_{1}, y_{2} ) \;\; \Longleftrightarrow \;\; x_{1} = y_{1} \; \wedge \; x_{2} = y_{2} . \]
\end{enumerate}

\noindent The categorical definition expresses these properties in arrow-theoretic terms, meaningful in any category.

\begin{mydefinition}
Let $A$, $B$ be objects in a category $\CC$. An $A$,$B$--pairing is a triple $(P, p_{1}, p_{2})$ where $P$ is an object, $p_1 : P \rarr A$ and $p_2 : P
\rarr B$. A morphism of $A$,$B$--pairings
\[ f : (P, p_{1}, p_{2}) \longrightarrow (Q, q_{1}, q_{2}) \]
is a morphism $f : P \rarr Q$ in $\CC$ such that $q_1 \circ f  = p_1$ and $q_2 \circ f = p_2$\,, \ie the following diagram commutes.
\[
\xymatrix@=13mm{
& P\ar[dl]_{p_{1}}\ar[dr]^{p_{2}}\ar[d]^f  \\
A & Q\ar[l]^{q_{1}}\ar[r]_{q_{2}} & B}\]
The $A$,$B$--pairings form a category $\mathbf{Pair}(A, B)$.
We say that $(A \times B, \pi_{1}, \pi_{2})$ is a \boldemph{product} of $A$ and $B$ if it
is \emph{terminal} in $\mathbf{Pair}(A, B)$.\deq
\end{mydefinition}
\begin{myexercise}
Verify that $\mathbf{Pair}(A, B)$ is a category.
\end{myexercise}
Note that products are specified by triples $A\overset{\pi_1}{\llarr}A\times B\overset{\pi_2}{\lrarr}B$, where $\pi_i$'s are called \emph{projections}. For
economy (and if projections are obvious) we may say that $A\times B$ is the product of $A$ and $B$. We say that $\CC$ \boldemph{has (binary) products} if
each pair of objects $A,B$ has a product in $\CC$. A direct consequence of the definition, by Proposition~\ref{inuniqueprop},  is that if products exist, they are unique up to (unique) isomorphism.

Unpacking the uniqueness condition from $\mathbf{Pair}(A, B)$ back to $\CC$ we obtain a more concise definition of products which we  use in
practice.

\begin{mydefinition}[Equivalent definition of product]\label{d:prod}
Let $A,B$ be objects in a category $\CC$. A product of $A$ and $B$ is an object $A\times B$ together with a pair of arrows
$A\overset{\pi_1}{\llarr}A\times B\overset{\pi_2}{\lrarr}B$ such that for every triple $A\overset{f}{\llarr}C\overset{g}{\lrarr}B$ there exists a
\emph{unique} morphism
\[ \langle f, g \rangle : C \longrightarrow A \times B \]
such that the following diagram commutes.
\[
\begin{aligned}
\xymatrix@=13mm{A & A\times B\ar[l]_{\pi_{1}}\ar[r]^{\pi_{2}} & B \\  & C\ar[ul]^f\ar[ur]_g\ar@{-->}[u]_{\ang{f,g}}}
\end{aligned}\qd[3]
\left(\begin{aligned} \pi_{1} \circ \langle f, g \rangle = f \\ \pi_{2} \circ \langle f, g \rangle = g
\end{aligned}\right)\]\deq[-1]
\end{mydefinition}
We call  $\ang{f,g}$ the \emph{pairing} of  $f$ and $g$.

Note that the above diagram features a \emph{dashed arrow}. Our intention with such diagrams is always to express the following idea: if the undashed part of the diagram commutes, then \emph{there exists a unique arrow} (the dashed one) such that the whole diagram commutes. In any case, we shall always spell out the intended statement explicitly.

We look at how this definition works in our standard example categories.
\begin{itemize}
\item In \textbf{Set}, products are the usual cartesian products.
\item In \textbf{Pos}, products are cartesian products with the pointwise order.
\item In \textbf{Top}, products are cartesian products with the product topology.
\item In $\text{\textbf{Vect}}_k$,  products are direct sums.
\item In a poset, seen as a category, products are \emph{greatest lower bounds}.
\end{itemize}
\begin{myexercise}
Verify these claims.
\end{myexercise}
The following proposition shows that the uniqueness of the pairing arrow can be specified purely equationally, by the equation:
\[ \forall h:C\rarr A\times B.\ h=\ang{\pi_1\circ h,\pi_2\circ h} \]

\begin{myproposition}\label{p:prod}
For any triple $A\overset{\pi_1}{\llarr}A\times B\overset{\pi_2}{\lrarr}B$ the following statements are equivalent.
\begin{description}
\item[(I)] For any triple $A\overset{f}{\llarr}C\overset{g}{\lrarr}B$ there exists a {unique} morphism $\ang{f,g}:C\rarr A\times B$
                        such that $\pi_1\circ\ang{f,g}=f$ and $\pi_2\circ\ang{f,g}=g$.
\item[(II)]For any triple $A\overset{f}{\llarr}C\overset{g}{\lrarr}B$ there exists a morphism $\ang{f,g}:C\rarr A\times B$
                        such that $\pi_1\circ\ang{f,g}=f$ and $\pi_2\circ\ang{f,g}=g$, and
                        moreover, for any $h:C\rarr A\times B$, $h=\ang{\pi_1\circ h,\pi_2\circ h}$.
\end{description}
\end{myproposition}
\proof
For (I)$\Rightarrow$(II), take any $h:C\rarr A\times B$\,; we need to show $h=\ang{\pi_1\circ h,\pi_2\circ h}$. We have
\[ \xymatrix{A & C\ar[l]_-{\pi_1\circ h}\ar[r]^-{\pi_2\circ h} & B} \]
and hence, by (I), there exists unique $k:C\rarr A\times B$ such that
\[ \pi_1\circ k= \pi_1\circ h \qd\land\qd \pi_2\circ k= \pi_2\circ h \tag{$*$}\]
Note now that ($*$) holds both for $k:=h$ and $k:=\ang{\pi_1\circ h,\pi_2\circ h}$, the latter because of (I). Hence, $h=\ang{\pi_1\circ h,\pi_2\circ
h}$.
\\
For (II)$\Rightarrow$(I), take any triple $A\overset{f}{\llarr}C\overset{g}{\lrarr}B$. By (II), we have that there exists an arrow $\ang{f,g}:C\rarr
A\times B$ such that $\pi_1\circ\ang{f,g}=f$ and $\pi_2\circ\ang{f,g}=g$. We need to show it is the unique such. Let $k:C\rarr A\times B$ s.t.
\[ \pi_1\circ k= f \qd\land\qd \pi_2\circ k= g \]
Then, by (II),
\[ k=\ang{\pi_1\circ k,\pi_2\circ k}=\ang{f,g} \]
as required. \qed[2]
In the following proposition we give some useful properties of products. First, let us introduce some notation for arrows: given $f_{1} : A_{1} \rarr
B_{1}$, $f_{2} : A_{2} \rarr B_{2}$, define
\[ f_{1} \times f_{2}  := \langle f_{1} \circ \pi_{1},  f_{2} \circ \pi_{2} \rangle : A_{1} \times A_{2} \lrar B_{1} \times B_{2}. \]

\begin{myproposition}\label{p:times} For any $f:A\rarr B$, $g:A\rarr C$, $h:A'\rarr A$, and any $p:B\rarr B'$, $q:C\rarr C'$,
\begin{itemize}
    \item $\langle f, g \rangle \circ h = \langle f \circ h, g \circ h \rangle$,
    \item $(p \times q) \circ \langle f, g \rangle = \langle p \circ f, q \circ g \rangle.$
\end{itemize}
\end{myproposition}
\proof For the first claim we have:
\[ \langle f, g \rangle \circ h = \langle \pi_{1} \circ (\langle f, g \rangle \circ h), \pi_{2} \circ (\langle f, g \rangle \circ h) \rangle =
    \langle f \circ h, g \circ h \rangle . \]
And for the second:
\[ \begin{array}{rcl}
(p \times q) \circ \langle f, g \rangle & = & \langle p \circ \pi_{1},  q \circ \pi_{2} \rangle \circ \langle f, g \rangle \\
& = &   \langle p \circ \pi_{1} \circ \langle f, g \rangle,  q \circ \pi_{2}  \circ \langle f, g \rangle \rangle \\
& = &   \langle p \circ f,  q \circ g \rangle .
\end{array}
\]
\qed[-1]

\paragraph{General Products} The notion of products can be generalised to arbitrary arities as follows.
A product for a family of objects $\{A_i \}_{i \in I}$ in a category $\CC$ is an object $P$ and morphisms
\[ p_i : P \longrightarrow A_i \quad (i \in I) \]
such that, for all objects $B$ and arrows
\[ f_i : B \longrightarrow A_i \quad (i \in I) \]
there is a \emph{unique} arrow
\[ g : B \longrightarrow P \]
such that, for all $i \in I$, the following diagram commutes.
\[ \xymatrix{B\ar[rr]^g\ar[dr]_{f_i} && P\ar[dl]^{p_i}\\ & A_i }\]
As before, if such a product exists, it is unique up to (unique) isomorphism. We write $P=\prod_{i\in I}A_i$ for the product object, and
$g=\ang{f_i\mid i\in I}$ for the unique morphism in the definition.

\begin{myexercise}
    What is the product of the empty family?
\end{myexercise}
\begin{myexercise}
Show that if a category has binary and nullary products then it has all finite products.
\end{myexercise}

\subsubsection{Coproducts}
We now investigate the dual notion to products: namely coproducts.
Formally, coproducts in $\CC$ are just products in  $\CCop$, interpreted back
in $\CC$ . We spell out the definition.
\begin{mydefinition}
Let $A,B$ be objects in a category $\CC$. A \emph{coproduct} of $A$ and $B$ is an object $A+B$ together with a pair of arrows
$A\overset{\incl[1]}{\lrarr}A+B\overset{\incl[2]}{\llarr}B$ such that for every triple $A\overset{f}{\lrarr}C\overset{g}{\llarr}B$ there exists a
\emph{unique} morphism
\[ [f, g] : A+B \longrightarrow C \]
such that the following diagram commutes.
\[
\begin{aligned}
\xymatrix@=13mm{A\ar[r]^-{\incl[1]}\ar[dr]_f & A+B\ar@{-->}[d]_{[f,g]} & B\ar[l]_-{\incl[2]}\ar[dl]^g \\ & C}
\end{aligned}\qd[3]
\left(\begin{aligned} {[f,g]}\circ\incl[1]  = f \\ [f,g]\circ\incl[2] = g
\end{aligned}\right)\]
\deq[-1]
\end{mydefinition}
We call the $\incl[i]$'s \emph{injections} and $[f,g]$ the \emph{copairing} of $f$ and $g$. As with pairings, uniqueness of copairings can be specified by an equation:
\[ \forall h : A+B \rarr C. \; h = [h\circ\incl[1],h\circ\incl[2]] \]

\paragraph{Coproducts in $\mathbf{Set}$}
This is given by \emph{disjoint union} of sets, which can be defined concretely e.g. by
\[ X + Y := \{ 1 \} \times X \, \cup \, \{ 2 \} \times Y . \]
We can define \emph{injections}
\[ \xymatrix{X\ar[r]^-{\incl[1]} & X+Y & Y\ar[l]_-{\incl[2]}} \]
\[ \mathsf{in}_1 (x) := (1, x)\,, \qquad \mathsf{in}_2 (y) := (2, y)\,. \]
Also, given functions $f: X \longrightarrow Z$ and $g : Y\longrightarrow Z$, we can define
\[ [f, g] : X + Y \longrightarrow Z \]
\[ [f,g](1,x) := f(x)\,, \qquad [f, g](2,y) := g(y)\,. \]
\begin{myexercise}
Check that this construction does yield coproducts in \tbf{Set}.
\end{myexercise}
Note that this example suggests that coproducts allow for \emph{definition by cases}.

Let us examine coproducts for some of our other standard examples.
\begin{itemize}
\item In \textbf{Pos}, disjoint unions (with the inherited orders) are coproducts.
\item In \textbf{Top}, topological disjoint unions are coproducts.
\item In $\text{\textbf{Vect}}_k$, direct sums are coproducts.
\item In a poset, \emph{least upper bounds} are coproducts.
\end{itemize}
\begin{myexercise}
Verify these claims.
\end{myexercise}
\begin{myexercise}
Dually to products, express coproducts as initial objects of a category $\mathbf{Copair}(A, B)$ of $A$,$B$--copairings.
\end{myexercise}

\subsection{Pullbacks and Equalisers}
We shall consider two further constructions of interest: \emph{pullbacks} and \emph{equalisers}.
\subsubsection{Pullbacks}
\begin{mydefinition}
Consider a pair of morphisms $A \overset{f}{\longrightarrow} C \overset{g}{\longleftarrow} B$.
The \boldemph{pull-back} of  $f$ along $g$ is a pair $A \overset{p}{\longleftarrow} D \overset{q}{\longrightarrow} B$ such that $f\circ p=g\circ q$
%
and, for any pair $A \overset{p'}{\llarr} D' \overset{q'}{\lrarr} B$ such that $f\circ p'=g\circ q'$, there exists a unique $h:D'\rarr D$ such that
$p'=p\circ h$ and $q'=q\circ h$. Diagrammatically,
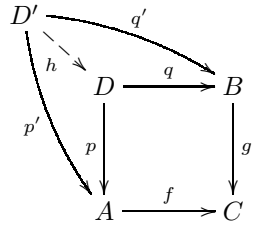
\[ \xymatrix@=5mm{D'\ar@{-->}[dr]_h\ar@/^/[drrr]^{q'}\ar@/_/[dddr]_{p'}\\& D\ar[rr]^q\ar[dd]_p && B\ar[dd]^g \\\\ & A\ar[rr]^f && C} \]\deq[-1]
\end{mydefinition}

\begin{myexample}
\begin{itemize}\item[]
\item In \textbf{Set} the pullback of $A \overset{f}{\longrightarrow} C \overset{g}{\longleftarrow} B$ is defined as a \emph{subset of the cartesian product}:
\[ A \times_C B \; = \; \{ (a, b) \in A \times B \mid f(a) = g(b) \}. \]
For example, consider a category $\CC$ with
\[ \Arr(\CC) \overset{\dm}{\longrightarrow} \Ob(\CC)  \overset{\cod}{\longleftarrow} \Arr(\CC) \, . \]
Then the pullback of $\dm$ along $\cod$ is the set of \emph{composable morphisms}, \ie pairs of morphisms $(f, g)$ in $\CC$ such that $f\circ g$ is well-defined.

\item In \textbf{Set} again, subsets (\ie inclusion maps) pull back to subsets:
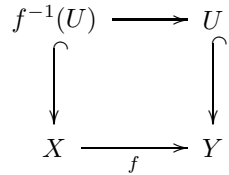
\[
\xymatrix@=10mm@M=2mm{f^{-1}(U)\ar[r]\ar@{^{(}->}[d] & U\ar@{^{(}->}[d] \\ X\ar[r]_f & Y}
\]
\end{itemize}
\end{myexample}

\begin{myexercise}
Let $\CC$ be a category with a terminal object $\ena$. Show that, for any $A,B\in Ob(\CC)$, the pullback of $A\xrightarrow{\tau_A}\ena\xleftarrow{\tau_B}B$ is the product of $A$ and $B$, if it exists.
\end{myexercise}
Just as for products, pullbacks can equivalently be described as terminal objects in suitable categories. Given a pair of morphisms $A \overset{f}{\longrightarrow} C
\overset{g}{\longleftarrow} B$, we define an \emph{$(f, g)$--cone} to be a triple $(D,p,q)$ such that the following diagram commutes.
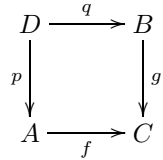
\[ \xymatrix@=10mm{ D\ar[r]^q\ar[d]_p & B\ar[d]^g \\ A\ar[r]_f & C} \]
A morphism of $(f, g)$--cones $h : (D_1 , p_1 , q_1 ) \rarr (D_2 , p_2 , q_2 )$ is a morphism $h : D_1 \rarr D_2$ such that the following diagram
commutes.
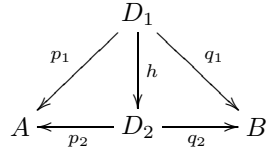
\[ \xymatrix@=10mm{& D_1\ar[dl]_{p_1}\ar[d]^h\ar[dr]^{q_1} \\ A & D_2\ar[l]^{p_2}\ar[r]_{q_2} & B} \]
We can thus form a category \textbf{Cone}$(f, g)$. A pull-back of  $f$ along $g$, if it exists,  is exactly a terminal object of \textbf{Cone}$(f, g)$. Once again, this shows the uniqueness of pullbacks up to unique isomorphism.

\subsubsection{Equalisers}
\begin{mydefinition}
Consider a pair of parallel arrows $\xymatrix@1{A\ \ar@<-.7mm>[r]_g\ar@<.7mm>[r]^f &\ B}$. An \boldemph{equaliser} of $(f,g)$ is an arrow $e:E\rarr A$
such that $f\circ e = g\circ e$ and, for any arrow $h:D\rarr A$ such that $f\circ h = g\circ h$, there is a unique $\hat{h}:D\rarr E$ so that
$h=e\circ\hat{h}$.
Diagrammatically,
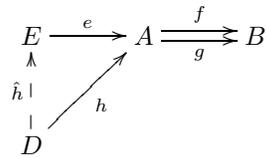
\[ \xymatrix@=10mm{E\ar[r]^e & A\ar@<-.7mm>[r]_g\ar@<.7mm>[r]^f & B\\ D\ar[ur]_h\ar@{-->}[u]^{\hat{h}}} \]
\deq[-1]
\end{mydefinition}
As for products, uniqueness of the arrow from $D$ to $E$ can be expressed equationally:
\[ \forall k:D\rarr E.\; \widehat{e\circ k}=k\,. \]
\begin{myexercise}
Why is $\widehat{e\circ k}$ well-defined for any $k:D\rarr E$?
Prove that the above equation is equivalent to the uniqueness requirement.
\end{myexercise}
\begin{myexample}
In $\Set$, the equaliser of $f,g$ is given by the inclusion
\[ \{ x\in A\ |\ f(x)=g(x)\} \hookrightarrow A \, . \]
This allows \emph{equationally defined subsets} to be defined as equalisers.
For example, consider the pair of maps $\xymatrix@1{\RR^{2}\ \ar@<-.7mm>[r]_g\ar@<.7mm>[r]^f &\ \RR}$, where
\[ f : (x, y) \mapsto x^{2} + y^{2}, \qquad g : (x, y) \mapsto 1 \, . \]
Then, the equaliser is the unit circle as a subset of $\RR^{2}$.
\end{myexample}

\subsection{Limits and Colimits}
The notions we have introduced so far are all special cases of a general notion of \emph{limits} in categories, and the dual notion of \emph{colimits}.
\begin{center}\renewcommand{\arraystretch}{1.4}
\begin{tabular}{|@{\;}l@{\quad}|@{\;}l@{\quad}|}\hline
\textbf{Limits} & \textbf{Colimits} \\ \hline
Terminal Objects & Initial Objects \\
Products & Coproducts \\
Pullbacks & Pushouts \\
Equalisers & Coequalisers\\\hline
\end{tabular}\captionof{table}{Examples of Limits and Colimits.}
\end{center}
An important aspect of studying any kind of mathematical structure is to see what limits and colimits the category of such structures has. We shall
return to these ideas shortly.

\subsection{Exercises}
\begin{enumerate}\renewcommand{\theenumi}{\textbf{\arabic{enumi}}}
\item      Give an example of a category where some pair of objects lacks a product or coproduct.
\item (\emph{Pullback lemma}) Consider the following commutative diagram.
\[\xymatrix@=10mm{A\ar[r]^f\ar[d]^u & B\ar[r]^g\ar[d]^v & C\ar[d]^w \\ D\ar[r]_h & E\ar[r]_i & F} \]
Given that the right hand square $BCEF$ and the outer square $ACDF$ are pullbacks, prove that the left hand square $ABDE$ is a pullback.

\item Consider $A \overset{f}{\longrightarrow} C \overset{g}{\longleftarrow} B$ with pullback $A\overset{p}{\longleftarrow}D\overset{q}{\longrightarrow}B$.
For each $A\overset{p'}{\longleftarrow}D'\overset{q'}{\longrightarrow}B'$ with $f\circ p'=g\circ q'$, let $\phi(p',q'):D'\rarr D$ be the arrow dictated
by the pullback condition. Express uniqueness of $\phi(p',q')$ equationally.

\end{enumerate}

\section{Functors}

Part of the ``categorical philosophy'' is:
\begin{center}
\fbox{\em Don't just look at the objects; take the morphisms into account too.}
\end{center}
We can also apply this to categories!

\subsection{Basics}
A ``morphism of categories'' is a \emph{functor}.

\begin{mydefinition}
A \boldemph{functor} $F : \CC \rarr \DD$ is given by:
\begin{itemize}
\item An object-map, assigning an object $FA$ of $\DD$ to every object $A$ of $\CC$.
\item An arrow-map, assigning an arrow $Ff : FA \rarr FB$ of $\DD$ to every arrow $f:A\rarr B$ of $\CC$,
    in such a way that composition and identities are preserved:
\[ F(g \circ f) = Fg \circ Ff\,, \qquad F \id[A] = \id[FA] . \]\deq[-1]
\end{itemize}
\end{mydefinition}
Note that we use the same symbol to denote the object- and arrow-maps; in practice, this never causes confusion. Since functors preserve domains and codomains of arrows, for each pair of objects $A$, $B$ of $\CC$, there is a well-defined map
\[ F_{A, B} : \CC(A, B) \rarr \DD(FA, FB) \, . \]
The conditions expressing preservation of composition and identities are called  \emph{functoriality}.

%
\begin{myexample}
Let $(P,\leq)$, $(Q,\leq)$ be preorders (seen as categories). A functor $F:(P,\leq)\lrarr(Q,\leq)$ is specified by an object-map, say $F:P\rarr Q$, and
an appropriate arrow-map. The arrow-map corresponds to the condition
\[ \forall p_1,p_2\in P.\, p_1\leq p_2 \implies F(p_1)\leq F(p_2)\,, \]
\ie to monotonicity of $F$. Moreover, the functoriality conditions are trivial since in the codomain $(Q,\leq)$ all hom-sets are singletons.\\
Hence, a functor between preorders is just a monotone map.
\end{myexample}
\begin{myexample}
Let $(M,\cdot,1)$, $(N,\cdot,1)$ be monoids. A functor $F:(M,\cdot,1)\lrarr(N,\cdot,1)$ is specified by a trivial object map (monoids are categories
with a single object) and an arrow-map, say $F:M\rarr N$. The functoriality conditions correspond to
\[ \forall m_1,m_2\in M.\, F(m_1\cdot m_2)=F(m_1)\cdot F(m_2)\,,\qquad F(1)=1\,,\]
\ie to $F$ being a monoid homomorphism.\\
Hence, a functor between monoids is just a monoid homomorphism.
\end{myexample}
Other examples are the following.
\begin{itemize}
\item Inclusion of a sub-category, $\CC\hookrightarrow\DD$, is a functor (by taking the identity map for object- and arrow-map).
\item The \emph{covariant} powerset functor $\PP:\mathbf{Set} \rarr \mathbf{Set}$:
\[ X \mapsto \PP(X)\,, \qquad (f : X \rarr Y) \mapsto \PP(f) := S \mapsto \{ f(x) \mid x \in S \} . \]
%
\item $U : \Mon \rarr \Set$ is the `forgetful' or `underlying' functor which sends a monoid to its set
  of elements, `forgetting' the algebraic structure, and sends a homomorphism to the corresponding function between sets. There are
  similar forgetful functors for other categories of structured sets. Why are these trivial-looking functors useful? --- We shall see!

\item Group theory examples. The assignment of the commutator sub-group of a group extends to a functor from $\textbf{Group}$ to $\textbf{Group}$;
    and the assignment of the quotient by this normal subgroup extends to a functor from $\textbf{Group}$ to $\textbf{AbGroup}$.
    The assignment of the centraliser of a group does not!
\item More sophisticated examples: e.g.~\emph{homology}. The basic idea of algebraic topology is that there are functorial assignments of algebraic objects (e.g.~groups) to topological spaces, and variants of this idea (`(co)homology theories') are pervasive throughout modern pure mathematics.
\end{itemize}

\paragraph{Functors `of several variables'} We can generalise the notion of a functor to a mapping from several domain categories to a codomain
category. For this we need the following definition.

\begin{mydefinition}
For categories $\CC,\DD$ define the \boldemph{product category} $\CC\times\DD$ as follows. An object in $\CC\times\DD$ is a pair of objects from $\CC$
and $\DD$, and an arrow in $\CC\times\DD$ is a pair of arrows from $\CC$ and $\DD$. Identities and arrow composition are defined componentwise:
\[ \id[(A,B)]:=(\id[A],\id[B])\,,\qquad (f,g)\circ(f',g'):=(f\circ f',g\circ g')\,. \]\deq[-1]
\end{mydefinition}
A functor `of two variables', with domains $\CC$ and $\DD$, to $\EE$ is simply a functor:
\[ F:\CC\times\DD\lrarr \EE\,. \]
For example, there are evident projection functors
\[ \CC \longleftarrow \CC \times \DD \lrarr \DD \, . \]

\subsection{Further Examples}
\paragraph{Set-valued functors} Many important constructions arise as functors $F:\CC\rarr\Set$. For example:
\begin{itemize}
\item If $G$ is a group, a functor $F:G\rarr\Set$ is an \emph{action of $G$ on a set}.
\item If $P$ is a poset representing time, a functor $F:P\rarr\Set$ is a notion of \emph{set varying through time}.
    This is related to \emph{Kripke semantics}, and to \emph{forcing arguments} in set theory.
\item Recall that $\TWO$ is the category $\xymatrix@=8mm{\bullet\ar@/^/[r]\ar@/_/[r] & \bullet}$\,. Then, functors $F:\TWO\rarr\Set$ correspond to    \emph{directed graphs} understood as in Definition~\ref{diagdef}, \ie as structures $(V, E, s, t)$, where $V$ is a set of vertices, $E$ is a set of edges, and $s, t : E \rarr V$ specify the source and target vertices for each edge.
\end{itemize}
Let us examine the first example in more detail. For a group $(G,\cdot,1)$, a functor $F:G\rarr\Set$ is specified by a set $X$ (to which the unique object
of $G$ is mapped), and by an arrow-map sending each element $m$ of $G$ to an endofunction on $X$, say $\actn{m}\uscore:X\rarr X$. Then, functoriality
amounts to the conditions
\[ \forall m_1,m_2\in G.\; F(m_1\cdot m_2)=F(m_1)\circ F(m_2)\,,\qquad F(1)=\id[X]\,,\]
that is, for all $m_1,m_2\in G$ and all $x\in X$,
\[ \actn{(m_1\cdot m_2)}x=\actn{m_1}\actn{m_2}x\,,\qquad \actn{1}x=x\,.\]
We therefore see that $F$ defines an action of $G$ on $X$.

\begin{myexercise}
Verify that functors $F:\TWO\rarr\Set$ correspond to {directed graphs}.
\end{myexercise}

\paragraph{Example: Lists}

Data-type constructors are functors. As a basic example, we consider
lists. There is a functor
\[ \List : \Set \longrightarrow \Set \]
which takes a set $X$ to the set of all finite lists (sequences) of
elements of $X$. $\List$ is functorial: its action on morphisms (\ie
functions, \ie (functional) programs) is given by \emph{maplist}:
\[ \frac{f : X \longrightarrow Y}{\List (f) : \List (X)
  \longrightarrow \List (Y)} \]
\[ \List (f) [ x_1 , \ldots , x_n ] := [ f(x_1 ), \ldots , f(x_n )] \]
We can upgrade $\List$ to a functor $\MList:\Set\rarr\Mon$ by mapping each set $X$ to the monoid $(\List(X),*,\epsilon)$ and $f:X\rarr Y$ to
$\List(f)$, as above. The monoid operation $*:\List(X)\times\List(X)\rarr\List(X)$ is list concatenation, and $\epsilon$ is the empty list. We call
$\MList(X)$ the \boldemph{free monoid} over $X$. This terminology will be justified in Chapter~5.

\paragraph{Products as functors} If a category $\CC$ has binary products, then there is automatically a functor
\[ \uscore\times\uscore : \CC \times \CC \lrarr \CC \]
which takes each pair $(A,B)$ to the product $A\times B$, and each $(f, g)$ to
\[ f \times g := \langle f \circ \pi_1 , g \circ \pi_2 \rangle\,. \]
Functoriality is shown as follows, using proposition~\ref{p:times} and uniqueness of pairings in its equational form.
\begin{align*}
  (f\times g)\circ(f'\times g') &= (f\times g)\circ\ang{f'\circ\pi_1,g'\circ\pi_2}=\ang{f\circ f'\circ\pi_1,g\circ g'\circ\pi_2} \\
    &= (f\circ f')\times(g\circ g')\,,\\
  \id[A]\times\id[B] &= \ang{\id[A]\circ\pi_1,\id[B]\circ\pi_2} = \ang{\pi_1\circ\id[A\times B],\pi_2\circ\id[A\times B]} = \id[A\times B]\,.
\end{align*}

\paragraph{The category of categories}
There is a category \textbf{Cat} whose objects are categories, and whose arrows are functors. Identities in $\Cat$ are given by identity functors:
\[ \Id_\CC:\CC\lrarr\CC:=A\mapsto A,\, f\mapsto f. \]
Composition of functors is defined in the evident
fashion. Note that if $F:\CC \rarr \DD$ and $G : \DD \rarr \EE$ then, for $f : A \rarr B$ in $\CC$,
\[ G \circ F (f) := G(F(f)) : G(F(A))  \lrarr G(F(B)) \]
so the types work out.
A category of categories sounds (and is) circular, but in practice is harmless: one usually makes some size restriction on the categories, and then \tbf{Cat} will be too `big' to be an object of itself. See Appendix~A.

Note that product categories are products in \textbf{Cat}! For any pair of categories $\CC,\DD$, set
\[ \CC\overset{\ppi_1}{\llarr}\CC\times\DD\overset{\ppi_2}{\lrarr}\DD \]
where $\CC\times\DD$ the product category (defined previously) and $\ppi_i$'s the obvious projection functors. For any pair of functors
$\CC\overset{F}{\llarr}\EE\overset{G}{\lrarr}\DD$, set
\[ \ang{F,G}:\EE\lrarr\CC\times \DD:=A\mapsto(FA,GA),\;f\mapsto(Ff,Gf)\,. \]
It is easy to see that $\ang{F,G}$ is indeed a functor. Moreover, satisfaction of the product diagram and uniqueness are shown exactly as in $\Set$.

\subsection{Contravariance}
By definition, the arrow-map of a functor $F$ is \emph{covariant}: it preserves the direction of arrows, so if $f:A\rarr B$ then $Ff:FA\rarr FB$. A
\emph{contravariant} functor $G$ does exactly the opposite: it reverses arrow-direction, so if $f:A\rarr B$ then $Gf:GB\rarr GA$. A concise way to
express contravariance is as follows.

\begin{mydefinition}
Let $\CC,\DD$ be categories. A \boldemph{contravariant} functor $G$ from $\CC$ to $\DD$ is a functor $G:\CCop\rarr\DD$.
(Equivalently, a functor $G : \CC \rarr \op{\DD}$.)\deq
\end{mydefinition}
Explicitly, a contravariant functor $G$ is given by an assignment of:
\begin{itemize}
\item an object $GA$ in $\DD$ to every object $A$ in $\CC$,
\item an arrow $Gf:GB\rarr GA$ in $\DD$ to every arrow $f:A\rarr B$ in $\CC$, such that (notice the change of order in composition):
\[ G(g \circ f) = Gf \circ Gg\,,\qquad G \id[A] = \id[GA] . \]
\end{itemize}
Note that functors of several variables  can be covariant in some variables and contravariant in others, e.g.
\[ F : \op{\CC} \times \DD \lrarr \EE\, . \]

\paragraph{Examples of Contravariant Functors}
\begin{itemize}
\item The contravariant powerset functor, $\op{\PP}: \op{\Set}\rarr \Set$\,, is given by:
\begin{align*}
\op{\PP}(X)&:=\PP(X)\,. \\
\op{\PP}(f:X\rarr Y) &: \PP(Y) \lrarr \PP(X) := T \mapsto \{ x \in X \mid f(x) \in T \}\,.
\end{align*}
\item The dual space functor on vector spaces:
\[  (\uscore)^{*}:\op{\Vect_k} \lrarr \Vect_k := V \mapsto V^{*} . \]
\end{itemize}
Note that these are both examples of the following idea: send an object $A$ into functions from $A$ into some fixed object. For example, the powerset can be written as $\PP(X) = 2^{X}$, where we think of a subset in terms of its characteristic function.

\paragraph{Hom-functors} We now consider some fundamental examples of $\Set$-valued functors. Given a category $\CC$ and an object $A$ of $\CC$, two  functors to $\Set$ can be defined:
\begin{itemize}
\item The covariant Hom-functor at $A$,
\[ \CC (A,\uscore) : \CC \lrarr \Set\,, \]
which is given by (recall that each $\CC(A,B)$ is a set):
\[  \CC (A,\uscore)(B) :=  \CC (A, B)\,, \qquad \CC (A,\uscore)(f : B \rarr C) := g \mapsto f \circ g\,. \]
We usually write $\CC(A,\uscore)(f)$ as $\CC(A,f)$. Functoriality reduces directly to the basic category axioms: associativity of composition and the unit laws for the
identity.
\item There is also a contravariant Hom-functor,
\[ \CC (\uscore\,,A) : \op{\CC} \lrarr \Set\,, \]
given by:
\[ \CC (\uscore\,,A)(B) := \CC (B,A)\,, \qquad \CC (\uscore\,,A)(h : C \rarr B) := g \mapsto g \circ h\,. \]
\end{itemize}
Generalising both of the above, we obtain a \boldemph{bivariant} Hom-functor,
\[ \CC(\uscore\,,\uscore) : \op{\CC} \times \CC \lrarr \Set\,. \]
\begin{myexercise}
Spell out  the definition of $\CC(\uscore\,,\uscore) : \op{\CC} \times \CC \lrarr \Set$. Verify carefully that it is a functor.
\end{myexercise}

\subsection{Properties of Functors}
\begin{mydefinition}\label{d:equiv1}
A functor $F:\CC\rarr\DD$ is said to be:
\begin{itemize}
\item \boldemph{faithful} if each map $F_{A,B}: \CC(A,B)\rarr\DD(FA,FB)$ is injective;
\item \boldemph{full} if each map $F_{A,B} :\CC (A, B) \rarr\DD (FA, FB)$ is surjective;
\item an \boldemph{embedding} if $F$ is full, faithful, and injective on objects;
\item an \boldemph{equivalence} if $F$ is  full, faithful, and \emph{essentially surjective}: \ie for every object $B$ of $\DD$ there is
    an object $A$ of $\CC$ such that $F(A) \cong B$;
\item an \boldemph{isomorphism} if there is a functor $G:\DD\rarr\CC$ such that
\[ G \circ F = \Id_{\CC}\,, \qquad F \circ G = \Id_{\DD}\,. \]
\deq[-1]
\end{itemize}
\end{mydefinition}
We say that categories $\CC$ and $\DD$ are isomorphic, $\CC\cong\DD$, if there is an isomorphism between them.
Note that this is just the usual notion of isomorphism applied to $\Cat$.
Examples:
\begin{itemize}
\item The forgetful functor $U : \Mon\rarr\Set$  is faithful, but not full. For the latter, note that not all functions $f:M\rarr N$ yield an arrow
    $f:(M,\cdot,1)\rarr(N,\cdot,1)$. Similar properties hold for other forgetful functors.
\item The free monoid functor $\MList:\Set\rarr\Mon$ is faithful, but not full.
\item The product functor $\uscore\times\uscore:\CC\times\CC\lrarr\CC$ is generally neither faithful nor full.
    For the latter, e.g.~in \Set, the function $f:\mathbb{N}^2\rarr\mathbb{N}^2:=(m,n)\mapsto(n,n)$ cannot be expressed in the form $f_1\times f_2$. Faithfulness of the functor is examined in exercise~\ref{ex:Funct}(2).
\item There is an equivalence between $\text{\textbf{FDVect}}_{k}$ the category of finite dimensional vector spaces over the field $k$, and
$\text{\textbf{Mat}}_{k}$, the category of matrices with entries in  $k$. Note that these categories are  very far from isomorphic! This example is elaborated in  exercise~\ref{ex:Funct}(1).
\end{itemize}

\paragraph{Preservation and Reflection}
Let $P$ be a property of arrows. We say that a functor $F :\CC \rarr \DD$ \emph{preserves} $P$ if whenever $f$ satisfies $P$, so does $F(f)$. We say
that $F$ \emph{reflects} $P$ if whenever $F(f)$ satisfies $P$, so does $f$. For example:
\begin{enumerate}\renewcommand{\theenumi}{\normalfont\alph{enumi}}
\item All functors preserve isomorphisms, split monics and split epics.
\item Faithful functors reflect monics and epics.
\item Full and faithful functors reflect isomorphisms.
\item Equivalences preserve monics and epics.
\item[$\bullet$\;] The forgetful functor $U:\Mon\rarr\Set$ preserves products.
\end{enumerate}
Let us show c; the rest are given as exercises below. So let $f:A\rarr B$ in $\CC$ be such that $Ff$ is an iso, that is, it has an inverse $g':FB\rarr
FA$. Then, by fullness, there exists some $g:B\rarr A$ so that $g'=Fg$. Thus,
\[ F(g\circ f)=Fg\circ Ff=g'\circ Ff=\id[FA]=F(\id[A])\,. \]
By faithfulness we obtain $g\circ f=\id[A]$\,. Similarly, $f\circ g=\id[B]$ and therefore $f$ is an isomorphism.

\begin{myexercise}
Show items a, b and d above.
\end{myexercise}
\begin{myexercise}
Show the following.
\begin{itemize}
  \item Functors do not in general reflect monics or epics.
  \item Faithful functors do not in general reflect isomorphisms.
  \item Full and faithful functors do not in general preserve monics or epics.
\end{itemize}
\end{myexercise}

\subsection{Exercises}\label{ex:Funct}
\begin{enumerate}\renewcommand{\theenumi}{\textbf{\arabic{enumi}}}
  \item Consider the category \textbf{FDVect}$_\R$ of finite dimensional vector spaces over $\R$, and \textbf{Mat}$_\R$ of matrices
        over $\R$. Concretely, \textbf{Mat}$_\R$ is defined as follows:
        \begin{align*}
            Ob(\mathbf{Mat}_\R) &:=\mathbb{N}\,,\\
            \mathbf{Mat}_\R(n,m)&:=\{M\ |\ M\text{ is an $n\times m$ matrix with entries in $\R$}\}\,.
        \end{align*}
        Thus, objects are natural numbers, and arrows $n\rarr m$ are $n\times m$ real matrices.
        Composition is matrix multiplication, and the identity on $n$ is the $n \times n$ identity matrix.
        \\
        Now let $F:\mathbf{Mat}_\R\rarr\mathbf{FDVect}_\R$ be the functor taking each $n$ to the vector space $\R^n$ and each $M:n\rarr m$ to the linear
        function
        \[ FM:\R^n\lrarr\R^m:=(x_1,...,x_n)\mapsto [x_1,...,x_n]M \]
        with the $1\times m$ matrix $[x_1,...,x_n]M$ considered as a vector in $\R^{m}$. Show that $F$ is full, faithful and essentially surjective,
        and hence  that \textbf{FDVect}$_\R$ and \textbf{Mat}$_\R$ are equivalent categories. Are they isomorphic?
  \item Let $\CC$ be a category with binary products such that, for each pair of objects $A,B$,
    \begin{equation}\label{ad}\tag{$*$}
    \CC(A,B)\neq\keno.
    \end{equation}
    Show that the product functor $F:\CC\times\CC\rarr\CC$ is faithful.
    \\ Would $F$ still be faithful in the absence of condition \eqref{ad}?
\end{enumerate}

\section{Natural Transformations}
\begin{quote}\it
``Categories were only introduced to allow functors to be defined;  functors were only introduced to allow natural transformations to be defined.''
\end{quote}
Just as categories have morphisms between them, namely functors, so functors have morphisms between them too\,---\,\emph{natural transformations}.

\subsection{Basics}
\begin{mydefinition}
Let $F, G : \CC \rarr\DD$ be functors. A \boldemph{natural transformation}
\[ t : F \lrarr G \]
is a family of morphisms in $\DD$ indexed by objects $A$ of $\CC$,
\[ \{\ t_A : FA \longrightarrow GA\ \}_{A\in Ob(\CC)} \]
such that, for all $f : A \rarr B$, the following diagram commutes.
\[ \xymatrix@=13mm{FA\ar[r]^{Ff}\ar[d]_{t_A} & FB\ar[d]^{t_B} \\ GA\ar[r]_{Gf} & GB} \]
This condition is known as \emph{naturality}.\\
If each $t_A$ is an isomorphism, we say that $t$ is a \boldemph{natural isomorphism}:
\[ t:F\overset{\cong}{\lrarr}G\,. \]\deq[-1]
\end{mydefinition}
Examples:
\begin{itemize}
\item Let $\Id$ be the identity functor on $\Set$, and $\times\circ\ang{\Id,\Id}$ be the functor taking each set $X$ to $X\times X$ and each function $f$ to $f\times f$. Then, there is a natural transformation $\Delta: \Id\lrarr\times\circ\ang{\Id,\Id}$ given by:
\[ \Delta_{X} :X\lrarr X\times X:= x \mapsto (x, x)\,. \]
Naturality amounts to asserting that, for any function $f:X\rarr Y$, the following diagram commutes.
\[ \xymatrix@=13mm{X\ar[r]^{f}\ar[d]_{\Delta_X} & Y\ar[d]^{\Delta_Y} \\ X\times X\ar[r]_{f\times f} & Y\times Y} \]
We call $\Delta$ the \emph{diagonal} transformation on $\Set$. In fact, it is the \emph{only} natural transformation between these functors.
\item The diagonal transformation can be defined for any category $\CC$ with binary products by setting, for each object $A$ in $\CC$,
    \[ \Delta_A:A\lrarr A\times A:=\ang{\id[A],\id[A]}\,. \]
  Projections also yield natural transformations. For example the arrows
    \begin{align*}
        \pi_{1(A,B)}&:A\times B\lrarr A
    \end{align*}
  specify a natural transformation $\pi_1:\times\rarr\ppi_1$\,. Note that $\times,\ppi_1:\CC\times\CC\rarr\CC$ are the functors for product and first
  projection respectively.
\item  Let $\CC$ be a category with terminal object $T$, and let $K_T:\CC\rarr\CC$ be the functor mapping all objects to $T$ and all arrows to
    $\id[T]$. Then, the canonical arrows
    \[ \tau_A:A\lrarr T \]
  specify a natural transformation $\tau:\Id\rarr K_T$\, (where $\Id$ the identity functor on $\CC$).
\item Recall the functor $\List:\Set\rarr\Set$ which takes a set $X$ to the set of finite lists with elements in $X$. We can define (amongst others) the following
natural transformations,
\[ \mathsf{reverse}: \List\longrightarrow \List\,,\quad \mathsf{unit}: \Id\longrightarrow \List\,,\quad
\mathsf{flatten}: \List\circ\List\longrightarrow \List\,, \]
by setting, for each set $X$,
\begin{align*}
 \mathsf{reverse}_X &: \List (X) \longrightarrow \List (X) := [x_1 , \ldots , x_n ] \mapsto [ x_n , \ldots , x_1 ]\,, \\[1mm]
 \mathsf{unit}_X &: X \longrightarrow \List (X) := x \mapsto [x]\,, \\[1mm]
 \mathsf{flatten}_X &: \List (\List (X)) \longrightarrow \List (X) \\
&:= [\,[x_1^1 , \ldots , x^1_{n_1}], \ldots , [x^k_1 , \ldots , x^k_{n_k}]\,] \mapsto [x_1^1 , \ldots \dots , x^k_{n_k}]\,.
\end{align*}
\item Consider the functor $P:=\times\circ \ang{U,U}$ with $U:\Mon\rarr\Set$, \ie
    \begin{align*}
      P:\Mon\lrar\Set:=(M,\cdot,1)\mapsto M\times M, \ f\mapsto f\times f\,.
    \end{align*}
    Then, the monoid operation yields a natural transformation $t:P\rarr U$ defined by:
    \[ t_{(M,\cdot,1)}:M\times M\lrar M:=(m,m')\mapsto m\cdot m'\,. \]
    Naturality corresponds to asserting that, for any $f:(M,\cdot,1)\rarr(N,\cdot,1)$, the following diagram commutes,
    \[ \xymatrix@=13mm{M\times M\ar[r]^{f\times f}\ar[d]_{t_M} & N\times N\ar[d]^{t_N} \\ M\ar[r]_{f} & N} \]
    that is, for any $m_1,m_2\in M$, \ $f(m_1)\cdot f(m_2)=f(m_1\cdot m_2)$.
\item If $V$ is a finite dimensional vector space, then $V$ is  isomorphic to both its first dual $V^{\ast}$ and to its second dual $V^{\ast \ast}$.\\
    However, while it is naturally isomorphic to its second dual, there is no natural isomorphism to the first dual. This was actually the original example which motivated Eilenberg and Mac Lane to define the concept of natural transformation; here naturality captures  \emph{basis independence}.
\end{itemize}
\begin{myexercise}
Verify naturality of diagonal transformations, projections and terminals for a category $\CC$ with finite products.
\end{myexercise}
\begin{myexercise}
Prove that the diagonal is the only natural transformation $\Id\lrarr\times\circ\ang{\Id,\Id}$ on $\Set$.
Similarly, prove that the first projection is the only natural transformation  $\times\rarr\ppi_1$ on $\Set$.
\end{myexercise}

\subsection{Further Examples}

\paragraph{Natural isomorphisms for products} Let $\CC$ be a category with finite products, \ie binary products and a terminal object $\ena$.
Then, we have the following canonical natural isomorphisms.
\begin{align*}
a_{A, B, C} &: A \times (B \times C) \iso (A \times B) \times C\,, \\
s_{A, B} &: A \times B \iso B \times A\,, \\
l_A &: \mathbf{1} \times A \iso A\,, \\
r_A &: A \times \mathbf{1} \iso A\,.
\end{align*}
The first two isomorphisms are meant to assert that the product is \emph{associative} and \emph{symmetric}, and the last two that $\ena$ is its \emph{unit}.
In later sections we will see that these conditions form part of the definition of  \emph{symmetric monoidal categories}.

These natural isomorphisms are defined explicitly by:
\begin{align*}
a_{A, B, C} &:= \langle \langle \pi_{1}, \pi_{1} \circ \pi_{2} \rangle , \pi_{2} \circ \pi_{2} \rangle\,, \\
s_{A,B} &:= \ang{\pi_2,\pi_1}\,, \\
l_A &:= \pi_2\,,\\
r_A &:= \pi_1\,.
\end{align*}
Since natural isomorphisms are a \emph{self-dual} notion, similar natural isomorphisms can be defined if $\CC$ has binary coproducts and an initial object.
\begin{myexercise}
Verify that these families of arrows are natural isomorphisms.
\end{myexercise}
\paragraph{Natural transformations between Hom-functors} Let $f : A \rarr B$ in a category $\CC$. Then, this induces a natural transformation
\begin{align*}
\CC(f,\uscore) &: \CC(B,\uscore) \lrarr \CC(A,\uscore)\,,\\[1mm]
\CC(f,\uscore)_{C} &: \CC(B,C) \lrarr \CC(A,C) := (g : B \rarr C) \mapsto (g \circ f : A \rarr C)\,.
\end{align*}
Note that  $\CC(f,\uscore)_{C}$ is the same as $\CC(f,C)$, the result of applying the contravariant functor $\CC(\uscore\,,C)$ to $f$. Hence,
naturality amounts to asserting that, for each $h:C\rarr D$, the following diagram commutes.
\[
\xymatrix@=13mm{
\CC(B,C)\ar[d]_{\CC(f,C)}\ar[r]^{\CC(B, h)} & \CC(B,D)\ar[d]^{\CC(f,D)} \\
\CC(A,C)\ar[r]_{\CC(A,h)} & \CC(A,D)}
\]
Starting from a $g : B \rarr C$, we compute:
\[ \CC(A,h)(\CC(f,C)(g)) = h \circ (g \circ f)  = (h \circ g) \circ f = \CC(f,D)(\CC(B,h)(g))\,. \]
The natural transformation $\CC(\uscore\,,f):\CC(\uscore\,,A)\rarr\CC(\uscore\,,B)$ is defined similarly.

\begin{myexercise}
Define the natural transformation $\CC(\uscore\,,f)$ and verify its naturality.
\end{myexercise}
There is a remarkable result, the \boldemph{Yoneda Lemma}, which says that \emph{every} natural transformation between Hom-functors comes from a
(unique) arrow in $\CC$ in the fashion described above.

\begin{lemma}
Let $A,B$ be objects in a category $\CC$. For each natural transformation $t:\CC(A,\uscore)\rarr\CC(B,\uscore)$, there is a unique arrow $f:B\rarr A$
such that
\[ t = \CC(f,\uscore)\,. \]
\end{lemma}
\proof Take any such $A,B$ and $t$ and let
\[ f:B\lrarr A:= t_A(\id[A])\,. \]
We want to show that $t=\CC(f,\uscore)$. For any object $C$ and any arrow $g:A\rarr C$, naturality of $t$ means that the following commutes.
\[ \xymatrix@=10mm{\CC(A,A)\ar[r]^{\CC(A,g)}\ar[d]_{t_A} & \CC(A,C)\ar[d]^{t_C}\\ \CC(B,A)\ar[r]_{\CC(B,g)} & \CC(B,C)} \]
Starting from $\id[A]$ we have that:
\[ t_C(\CC(A,g)(\id[A]))=\CC(B,g)(t_A(\id[A]))\,,\; \text{ \ie } \; t_C(g)=g\circ f\,. \]
Hence, noting that  $\CC(f,C)(g)=g\circ f$, we obtain $t=\CC(f,\uscore)$.\\
For uniqueness we have that, for any $f,f':B\rarr A$, if $\CC(f,\uscore)=\CC(f',\uscore)$ then
\[ f=\id[A]\circ f=\CC(f,A)(\id[A])=\CC(f',A)(\id[A])=\id[A]\circ f'=f'. \]\qed[-1]

\begin{myexercise}
\label{contrayonex}
Prove a similar result for contravariant hom-functors.
\end{myexercise}

\paragraph{Alternative definition of equivalence} Another way of defining equivalence of categories is as follows.
\begin{mydefinition}\label{d:equiv2}
We say that categories $\CC$ and $\DD$ are \boldemph{equivalent}, $\CC \simeq \DD$, if there are functors $F:\CC\rarr\DD$, $G:\DD\rarr\CC$ and natural
isomorphisms
\[ G \circ F \cong \Id_{\CC}\,, \quad F \circ G \cong \Id_{\DD}\,. \]\deq[-1]
\end{mydefinition}
%

\subsection{Functor Categories}
Suppose we have functors $F, G, H : \CC \rarr \DD$ and natural transformations
\[ t : F \lrarr G\,, \quad u : G \lrarr H\,. \]
Then, we can compose these natural transformations, yielding $u \circ t : F \rarr H$:
\[ (u \circ t)_{A} := FA \overset{t_{A}}{\lrarr}GA \overset{u_{A}}{\lrarr}HA . \]
Composition is associative, and has as identity the natural transformation
\[ I_{F} : F \lrarr F := \{\ (I_{F})_{A} := \id[A] : FA\lrarr FA\ \}_{A}\,. \]
These observations lead us to the following.

\begin{mydefinition}
For categories $\CC,\DD$ define the \boldemph{functor category} $\mathsf{Func}(\CC, \DD)$ by taking:
\begin{itemize}
\item Objects: functors $F : \CC \rarr \DD$.
\item Arrows: natural transformations $t : F \rarr G$.
\end{itemize}
Composition and identities are given as above.\deq
\end{mydefinition}
\begin{myremark}
We see that in the category $\Cat$ of categories and functors, each hom-set $\Cat (\CC, \DD)$ \emph{itself has the structure of a category}. In fact,
$\Cat$ is the basic example of a ``2-category'', \ie of a category where hom-sets are themselves categories.
\end{myremark}
Note that a natural isomorphism is precisely an isomorphism in the functor category.
Let us proceed to some examples of functor categories.
\begin{itemize}
\item Recall that, for any group $G$, functors from $G$ to $\Set$ are $G$-actions on sets. Then, $\mathsf{Func}(G, \Set)$ is the category
of $G$-actions on sets and \emph{equivariant functions}: $f : X \rarr Y$ such that $f(\actn m x) = \actn m f(x)$.

\item $\mathsf{Func}(\TWO, \Set)$: Graphs and graph homomorphisms.

\item If $F, G : P\rarr Q$ are monotone maps between posets, then $t : F \rarr G$ means that
\[ \forall x \in P. \; Fx \leq Gx\,. \]
Note that in this case naturality is trivial (hom-sets are singletons in $Q$).
\end{itemize}
\begin{myexercise}
Verify the above descriptions of $\mathsf{Func}(G, \Set)$ and $\mathsf{Func}(\TWO, \Set)$.
\end{myexercise}

\begin{myremark}
The composition of natural transformations defined above is called \emph{vertical composition}. The reason for this terminology is depicted below.
\[
\xymatrix@C=10mm{\CC\ar@/^10mm/[rr]^{F}_{}="F"\ar@/_10mm/[rr]_{H}^{}="H"\ar[rr]^{}="GU"_{}="GD"^<<<<<<<G && \DD\quad\ar@{=>}[r]&
\quad\CC\ar@/^5mm/[rr]^{F}_{}="F2"\ar@/_5mm/[rr]_{H}^{}="H2" && \DD \ar"F";"GU"^t \ar"GD";"H"^u \ar"F2";"H2"^{u\,\circ\,t}}
\]
As expected, there is also a \emph{horizontal composition}, which is given as follows.
\[
\xymatrix@C=10mm{\CC\ar@/^5mm/[rr]^F_{}="F"\ar@/_5mm/[rr]_G^{}="G" && \DD\ar@/^5mm/[rr]^{F'}_{}="F'"\ar@/_5mm/[rr]_{G'}^{}="G'" && \EE\quad\ar@{=>}[r]&
\quad\CC\ar@/^5mm/[rr]^-{F'\circ\,F}_{}="FF"\ar@/_5mm/[rr]_{G'\circ\,G}^{}="GG" && \EE \ar"F";"G"^{t}\ar"F'";"G'"^{t'}\ar"FF";"GG"^{t'\bullet\,t}}
\]
\end{myremark}

\subsection{Exercises}
\begin{enumerate}\renewcommand{\theenumi}{\textbf{\arabic{enumi}}}
\item By identifying the relevant functors, express pairing $\langle\uscore,\uscore\rangle$
 as a natural transformation. What does naturality correspond to explicitly?
  \item Show that the two definitions of equivalence of categories, namely
  \begin{enumerate}
    \item $\CC$ and $\DD$ are equivalent if there is an equivalence $F:\CC\rarr\DD$ (definition~\ref{d:equiv1}),
    \item $\CC$ and $\DD$ are equivalent if there are $F:\CC\rarr\DD$, $G:\DD\rarr\CC$, and isomorphisms $F\circ G\cong\Id_{\DD}$, $G\circ
    F\cong\Id_{\CC}$ (definition~\ref{d:equiv2}),
  \end{enumerate}
  are: equivalent!
  Note that this will need the Axiom of Choice.
  \item Define a relation on objects in a category $\CC$ by: $A \cong B$ iff $A$ and $B$ are isomorphic.
    \begin{enumerate}
      \item Show that this relation is an equivalence relation.
    \end{enumerate}
    Define a \emph{skeleton} of $\CC$ to be the (full) subcategory obtained by choosing one object from each equivalence class of $\cong$\, (note that
    this involves choices, and is not uniquely defined).
    \begin{enumerate}\setcounter{enumii}{1}
      \item Show that $\CC$ is equivalent to any skeleton.
      \item Show that any two skeletons of $\CC$ are isomorphic.
      \item Give an example of a category whose objects form a proper class, but whose skeleton is finite.
    \end{enumerate}

\item   Given a category $\CC$, we can define a functor
\[ y : \CC \lrarr \mathsf{Func}(\CCop, \Set) := A \mapsto \CC (\uscore\,,A), \;\; f \mapsto  \CC (\uscore\,,f) \, . \]
Prove carefully that this is indeed a functor. Use exercise~\ref{contrayonex} to  conclude that $y$ is full and faithful. Prove that it is also injective on objects, and hence an embedding. It is known as the \emph{Yoneda embedding}.

\item Define the horizontal composition $u \bullet t$ of natural transformations explicitly. Prove that it is associative.
\end{enumerate}

\section{Universality and Adjoints}\label{s:Univ}

There is a fundamental triad of categorical notions:
\[ \text{\emph{Functoriality, Naturality, Universality.}} \]
We have studied the first two notions explicitly. We have also seen many \emph{examples} of universal definitions, notably the various notions of limits and colimits considered in section~\ref{s:BasConstr}. It is now time to consider universality in general; the proper formulation of this fundamental and pervasive notion is one of the major achievements of basic category theory.

Universality arises when we are interested in finding \emph{canonical solutions} to problems of construction: that is, we are interested not just in the \emph{existence} of a solution but in its \emph{canonicity}. This canonicity should guarantee uniqueness, in the sense we have become familiar with: a canonical solution should be \emph{unique up to (unique) isomorphism}.

The notion of canonicity has a simple interpretation in the case of posets, as an \emph{extremal solution}: one that is the least or the greatest among all solutions. Such an extremal solution is obviously unique.
For example, consider the problem of finding a lower bound of a pair of elements $A$, $B$ in a poset $P$: a \emph{greatest lower bound} of $A$ and $B$ is an extremal solution to this problem. As we have seen, this is the specialisation to posets of the problem of constructing a product:
\begin{itemize}
\item[$\rightsquigarrow$] A product of $A$, $B$ in a poset is an element $C$ such that $C \leq A$ and $C\leq B$, ($C$ is a lower bound);
\item[$\rightsquigarrow$] and for any other solution $C'$, \ie $C'$ such that $C'\leq A$ and $C'\leq B$, we have $C'\leq C$.
    ( $C$ is a greatest lower bound.)
\end{itemize}
Because the ideas of universality and adjunctions have an appealingly simple form in posets, which is, moreover,  useful in its own right, we will develop the ideas in that special case first, as a prelude to the general discussion for categories.

\subsection{Adjunctions for Posets}

Suppose $g : Q \rarr P$ is a monotone map between posets. Given $x \in P$,
a \boldemph{$g$-approximation of $x$} (from above) is an element $y \in Q$ such that $x \leq g(y)$.\\
A \boldemph{best $g$-approximation of $x$} is an element $y \in Q$ such that
\[ x \leq g(y)\ \land\ \forall z \in Q.\,(\, x \leq g(z) \implies y \leq z \,)\,. \]
If a best $g$-approximation exists then it is clearly unique.

\paragraph{Discussion} It is worth clarifying the notion of best $g$-approximation. If $y$ is a best $g$-approximation to $x$, then in particular, by monotonicity of $g$, $g(y)$ is the least element of the set of all $g(z)$ where $z \in Q$ and $x \leq g(z)$. However, the property of being a best approximation is much \emph{stronger} than the mere existence of a least element of this set. We are asking for \emph{$y$ itself} to be the least, \emph{in $Q$},  among all elements $z$ such that $x \leq g(z)$. Thus, even if $g$ is \emph{surjective}, so that for every $x$ there is a $y \in Q$ such that $g(y) = x$, there need not exist a \emph{best} $g$-approximation to $x$. This is exactly the issue of having a \emph{canonical choice} of solution.
\begin{myexercise}
Give an example of a surjective monotone map $g : Q \rarr P$ and an element $x \in P$ such that there is no best $g$-approximation to $x$ in $Q$.
\end{myexercise}

\noindent If such a best $g$-approximation $f(x)$ exists for all $x \in P$ then we have a function $f : P \rarr Q$ such that, for all $x \in P$, $z \in Q$:
\begin{equation}\label{adjeq}
x \leq g(z) \;\; \Longleftrightarrow \;\; f(x) \leq z\, .
\end{equation}
We say that $f$ is the \boldemph{left adjoint} of $g$, and $g$ is the \boldemph{right adjoint} of $f$.
It is immediate from the definitions that the left adjoint of $g$, if it exists,  is uniquely determined by $g$.

\begin{myproposition}
If such a function $f$ exists, then it is  monotone. Moreover,
\[ \id[P]\leq g\circ f\,, \qquad f \circ g \leq \id[Q]\,, \qquad f \circ g \circ f = f\,, \qquad g \circ f \circ g = g\,. \]
\end{myproposition}
\proof  If we take $z = f(x)$ in equation (\ref{adjeq}), then since $f(x) \leq f(x)$, $x \leq g
\circ f(x)$. Similarly, taking $x = g(z)$ we obtain  $f \circ g(z) \leq z$.
Now, the ordering on functions $h, k : P \lrarr Q$  is the \emph{pointwise order}:
\[ h \leq k \iff \forall x \in P. \, h(x) \leq k(x) . \]
This gives the first two equations.

Now, if $x \leq_{P} x'$ then $x \leq x' \leq g \circ f(x')$, so $f(x')$ is a $g$-approximation of $x$, and hence $f(x) \leq f(x')$. Thus,
$f$ is monotone.

Finally, using the fact that composition is monotone with
respect to the pointwise order on functions, and the first two equations:
\[ g = \id[P] \circ g \leq g \circ f \circ g \leq g \circ \id[Q] = g, \]
and hence $g = g \circ f \circ g$. The other equation is proved similarly. \qed[2]
Examples:
\begin{itemize}
\item Consider the inclusion map
\[ i : \ZZ \hookrightarrow \RR\,. \]
This has both a left adjoint $f^{L}$ and a right adjoint $f^{R}$, where $f^{L}, f^{R} : \RR \rarr \ZZ$. For all $z \in \ZZ$, $r \in \RR$:
\[ z \leq f^{R}(r)  \;\; \Longleftrightarrow \;\;  i(z) \leq r \,,\qquad f^{L}(r) \leq z  \;\; \Longleftrightarrow \;\;  r \leq i(z)\,. \]
We see from these defining properties that the right adjoint maps a real $r$ to the \emph{greatest integer below it} (the extremal solution to finding
an integer below a given real). This is the standard \emph{floor function}.\\ Similarly, the left adjoint maps a real to the least integer above it
yielding the \emph{ceiling function}. Thus:
\[ f^{R}(r) = \lfloor r \rfloor\,, \qquad f^{L}(r) = \lceil r \rceil\,. \]
\item Consider a relation $R \subseteq  X \times Y$. $R$ induces a function:
\[ f_{R} : \pow(X) \lrarr \pow(Y):= S \mapsto \{ y \in Y \mid \exists x \in S.\, xRy \}\,. \]
This has a right adjoint $[R] : \pow(Y) \lrarr \pow(X)$:
\[ S \subseteq [R]T \;\; \Longleftrightarrow \;\; f_{R}(S) \subseteq T\,. \]
The definition of $[R]$ which satisfies this condition is:
\[ [R]T := \{ x \in X \mid \forall y \in Y. \, xRy \; \Rightarrow \; y \in T \}\,. \]
If we consider a set of \emph{worlds} $W$ with an \emph{accessibility relation} $R \subseteq W \times W$ as in Kripke semantics for modal logic, we see
that $[R]$ gives the usual Kripke semantics for the modal operator $\Box$, seen as a propositional operator mapping the set of worlds satisfied by a
formula $\phi$ to the set of worlds satisfied by $\Box \phi$.
\\
On the other hand, if we think of the relation $R$ as the denotation of a (possibly non-deterministic) program, and $T$ as a predicate on
\emph{states}, then $[R]T$ is exactly the \emph{weakest precondition} $\mathbf{wp}(R, T)$. In \emph{Dynamic Logic}, the two settings are combined, and
we can write expressions such as $[R]T$ directly, where $T$ will be (the denotation of) some formula, and $R$ the relation corresponding to a program.
\item Consider a function $f : X \rarr Y$. This induces a function:
\[ f^{-1} : \PP(Y) \lrarr \PP(X) := T \mapsto \{ x \in X \mid f(x) \in T \}\,. \]
This function $f^{-1}$ has both a left adjoint $\exists (f) : \PP(X) \lrarr \PP(Y)$, and a right adjoint $\forall (f) : \PP(X) \lrarr \PP(Y)$. For all
$S \subseteq X$, $T \subseteq Y$:
\[ \exists(f)(S) \subseteq T \iff S \subseteq f^{-1}(T)\,,\qquad  f^{-1}(T) \subseteq S \iff  T \subseteq \forall(f)(S)\,. \]
How can we define $\forall(f)$ and $\exists(f)$ explicitly so as to fulfil these defining conditions? --\,As follows:
\begin{align*}
\exists(f)(S) &:= \{ y \in Y \mid \exists x \in X. \, f(x) = y \; \wedge \; x \in S \}\,,  \\
\forall(f)(S) &:= \{ y \in Y \mid \forall x \in X. \, f(x) = y \; \Rightarrow \; x \in S \}\,.
\end{align*}
If $R \subseteq X \times Y$, which we write in logical notation as $R(x, y)$, and we take the projection function $\pi_{1} : X \times Y \lrarr X$,
then:
\[ \forall(\pi_{1})(R) \equiv \forall y. \, R(x, y)\,, \qquad \exists(\pi_{1})(R) \equiv \exists y. \, R(x, y)\,. \]
This extends to an algebraic form of the usual Tarski model-theoretic semantics for first-order logic, in which:
\begin{center}
\fbox{\emph{Quantifiers are Adjoints}}
\end{center}
\end{itemize}

\paragraph{Couniversality}
We can dualise the discussion, so that starting with a monotone map $f : P \rarr Q$ and $y \in Q$, we can ask for the best $P$-approximation to $y$ from below: $x \in P$ such that $f(x) \leq y$, and for all $z \in P$:
\[ f(z) \leq y \;\; \Longleftrightarrow \;\; z \leq x . \]
If such a best approximation $g(y)$ exists for all $y \in Q$, we obtain a monotone map $g : Q \rarr P$ such that $g$ is right adjoint to $f$.
From the symmetry of the definition, it is clear that:
\begin{center}
 $f$ is the left adjoint of $g$ $\;\; \Longleftrightarrow \;\;$ $g$ is the right adjoint of $f$
\end{center}
and each determines the other uniquely.

\subsection{Universal Arrows and Adjoints}
Our discussion of best approximations for posets is lifted to general categories as follows.

\begin{mydefinition}
Let $G : \DD\rarr\CC$ be a functor, and $C$ an object of $\CC$. A \boldemph{universal arrow from $C$ to $G$} is a pair $(D,\eta)$ where $D$ is an object of $\DD$ and
\[ \eta : C \longrightarrow G(D)\,, \]
such that, for any object $D'$ of $\DD$ and morphism $f : C\rarr G(D')$, \emph{there exists a unique} morphism $\hat{f} : D\rarr D'$ in $\DD$ such that
$f=G(\hat{f})\circ\eta$\,.
\\
Diagrammatically:
\[
\xymatrix@=13mm{C\ar[r]^{\eta}\ar[dr]_f & G(D)\ar@{-->}[d]^{G(\hat{f})} & D\ar@{-->}[d]^{\hat{f}}\\ & G(D') & D'}
\]
\deq[-1]
\end{mydefinition}
As in previous cases, uniqueness can be given a purely equational specification:
\begin{equation}
\label{univeq}
\forall h:D\lrarr D'.\;\widehat{G(h)\circ\eta}=h\,.
\end{equation}
\begin{myexercise}
Show that if $(D,\eta)$ and $(D',\eta')$ are universal arrows from $C$ to $G$ then there is a unique isomorphism $D\cong D'$.
\end{myexercise}
\begin{myexercise}
Check that the equational specification of uniqueness (\ref{univeq}) is valid.
\end{myexercise}
Examples:
\begin{itemize}
  \item Take $U : \Mon\rarr\Set$. Given a set $X$, the universal arrow is
    \[
    \eta_{X} : X \lrarr U(\MList(X)) := x \mapsto [x]\,.
    \]
    Indeed, for any monoid $(M,\cdot,1)$ and any function $f:X\rarr M$, set
    \[ \hat{f}:\MList(X)\lrarr(M,\cdot,1):= [x_1,\dots,x_n] \mapsto f(x_1)\cdot\,\cdots\,\cdot f(x_n)\,. \]
    It is easy to see that $\hat{f}$ is a monoid homomorphism, and that $U(\hat{f})\circ\eta_X=f$. Moreover, for uniqueness we have that, for any
    $h:\MList(X)\rarr(M,\cdot,1)$,
    \begin{align*}
        \widehat{U(h)\circ\eta_X} &=\widehat{x\mapsto h([x])}=[x_1,\dots,x_n]\mapsto h([x_1])\cdot\,\cdots\,\cdot h([x_n]) \\
        &=[x_1,\dots,x_n]\mapsto h([x_1]*\cdots*[x_n]) \\
        &=[x_1,\dots,x_n]\mapsto h([x_1,\dots,x_n]) = h\,.
    \end{align*}
  \item Let $K:\CC\rarr\1$ be the unique functor to the one-object/one-arrow category. A universal arrow from the object of $\1$ to $K$ corresponds
    to an initial object in $\CC$. \\
    Indeed, such a universal arrow is given by an object $I$ of $\CC$ (and a trivial arrow in $\1$), such that for any $A$ in $\CC$ (and relevant arrow
    in $\1$) there exists a unique arrow from $I$ to $A$ (such that a trivial condition holds).
  \item Consider the functor $\ang{\Id_\CC,\Id_\CC}:\CC\rarr\CC\times\CC$, taking each object $A$ to $(A,A)$ and each arrow $f$ to $(f,f)$.
    A universal arrow from an object $(A,B)$ of $\CC\times\CC$ to $\ang{\Id_\CC,\Id_\CC}$ corresponds to a coproduct of $A$ and $B$.
\end{itemize}
\begin{myexercise}
Verify the description of coproducts as universal arrows.
\end{myexercise}
As in the case of posets, a related notion to universal arrows is that of \emph{adjunction}.

\begin{mydefinition}
Let $\CC,\DD$ be categories. An \boldemph{adjunction} from $\CC$ to $\DD$ is a triple $(F,G,\theta)$, where $F$ and $G$ are functors
\[ \xymatrix@=13mm{{\CC\ }\ar@<.8mm>[r]^F & {\ \DD}\ar@<.8mm>[l]^G} \]
and $\theta$ is a family of bijections
\[ \theta_{A,B}:\CC(A,G(B))\overset{\cong}{\lrarr}\DD(F(A),B)\,, \]
for each $A\in Ob(\CC)$ and $B\in Ob(\DD)$, \textbf{natural} in $A$ and $B$.\\
We say that $F$ is \boldemph{left adjoint} to $G$, and $G$ is \boldemph{right adjoint} to $F$.\deq
\end{mydefinition}
Note that $\theta$ should be understood as the ``witnessed'' form --- \ie arrows instead of mere relations --- of the defining condition for adjunctions in the case of posets:
\[ x \leq g(y) \iff f(x) \leq y . \]
This is often displayed as a two-way `inference rule':
\[ \frac{A\lrarr GB}{FA\lrarr B} \]
Naturality of $\theta$ is expressed as follows: for any $f:A\rarr G(B)$ and any $g:A'\rarr A$, $h:B\rarr B'$,
\begin{align*}
  \theta_{A',B}(f\circ g) &= \theta_{A,B}(f)\circ F(g)\,, \\
  \theta_{A,B'}(G(h)\circ f) &= h\circ\theta_{A,B}(f)\,.
\end{align*}
Note that $f,g$ are in $\CC$, and $h$ is in $\DD$.
In one line:
\[ \theta_{A',B'}(G(h)\circ f\circ g) = h\circ\theta_{A,B}(f)\circ F(g)\,.\]
Diagrammatically:
\[ \xymatrix@=10mm{
    \CC(A,GB')\ar[d]_{\theta_{A,B'}} & \CC(A,GB)\ar[l]_-{\CC(A,Gh)}\ar[r]^-{\CC(g,GB)}\ar[d]_{\theta_{A,B}} & \CC(A',GB)\ar[d]^{\theta_{A',B}}\\
    \DD(FA,B') & \DD(FA,B)\ar[l]^-{\DD(FA,h)}\ar[r]_-{\DD(Fg,B)} & \DD(FA',B)}\quad
\xymatrix@=10mm{\CC(A,GB)\ar[d]_{\theta_{A,B}}\ar[r]^-{\CC(g,Gh)} & \CC(A',GB')\ar[d]^{\theta_{A',B'}} \\
    \DD(FA,B)\ar[r]_-{\DD(Fg,h)} & \DD(FA',B')}
\]
Thus, $\theta$ is in fact a natural isomorphism
\[ \theta:\CC(\uscore\,,G(\uscore))\xrightarrow{\cong}\DD(F(\uscore),\uscore)\,, \]
where $\CC(\uscore\,,G(\uscore)):\CCop\times\DD\rarr\Set$ is the result of composing the bivariant hom-functor $\CC(\uscore\,,\uscore)$ with
$\Id_{\CCop}\times G$, and $\DD(F(\uscore),\uscore)$ is similar.

In the next propositions we show that universal arrows and adjunctions are equivalent notions.

\begin{myproposition}[Universals define adjunctions]\label{p:UnivtoAdj}
Let $G:\DD\rarr\CC$. If for every object $C$ of $\CC$ there exists a universal arrow $\eta_{C} : C \rarr G(F(C))$, then:
\begin{enumerate}\renewcommand{\theenumi}{\rm\arabic{enumi}}
\item $F$ uniquely extends to a functor $F:\CC\rarr\DD$ such that $\eta : \Id_{\CC} \rarr G \circ F$ is a natural transformation.
\item $F$ is uniquely determined by $G$ (up to unique natural isomorphism), and vice versa.
\item For each pair of objects $C$ of $\CC$ and $D$ of $\DD$, there is a natural bijection:
\[ \theta_{C,D} : \CC(C, G(D))\cong\DD(F(C), D) \,. \]
\end{enumerate}
\end{myproposition}
\proof %
For 1, we extend $F$ to a functor as follows. Given $f : C \rarr C'$ in $\CC$, we consider the composition
\[ \eta_{C'} \circ f : C \lrarr GFC' . \]
By the universal property of $\eta_{C}$, there exists a unique arrow $Ff : FC \rarr FC'$ such that the following diagram commutes.
\[
\xymatrix{C\ar[r]^-{\eta_{C}}\ar[d]_{f} & GFC\ar[d]^{GFf} \\  C'\ar[r]_-{\eta_{C'}} & GFC'}
\]
Note that the above is the naturality diagram for $\eta$ on $C$, hence the arrow-map thus defined for $F$ is the unique candidate that makes
$\eta$ a natural transformation.
\\
It remains to verify the functoriality of $F$. To show that $F$ preserves composition, consider $g : C' \rarr C''$. We have the following commutative
diagram,
\[
\xymatrix{C\ar[r]^f\ar[d]_{\eta_C} & C'\ar[r]^g\ar[d]_{\eta_{C'}} & C''\ar[d]^{\eta_{C''}}\\ GFC\ar[r]_{GFf} & GFC'\ar[r]_{GFg} & GFC''}
\]
from which it follows that
\[ G(Fg \circ Ff) \circ \eta_{C} = GFg \circ GFf \circ \eta_{C} = \eta_{C''} \circ g \circ f\,, \]
\[\therefore\ F(g\circ f)=\widehat{\eta_{C''}\circ g\circ f}=\widehat{G(Fg\circ Ff)\circ\eta_C}=Fg\circ Ff\,, \]
where the last equality above holds because of~\eqref{univeq}. The verification that $F$ preserves identities is similar.
\\
For 2, we have that each $FC$ is determined uniquely up to unique isomorphism, by the universal property, and once the object part of $F$ is fixed, the
arrow part is uniquely determined.
\\
For 3, we need to define a natural isomorphism $\theta_{C,D} : \CC(C, G(D))\cong\DD(F(C), D)$.
Given $f : C \rarr GD$, $\theta_{C,D}(f)$ is defined to be the unique arrow $FC \rarr D$ such that the following commutes, as dictated by universality.
\[
\xymatrix@=10mm{C\ar[r]^{\eta_C}\ar[dr]_f & GFC\ar[d]^{G(\theta_{C,D}(f))}\\ & GD}
\]
Suppose that $\theta_{C,D}(f) = \theta_{C,D}(g)$. Then
\[ f = G(\theta_{C,D}(f)) \circ \eta_{C} =  G(\theta_{C,D}(g)) \circ \eta_{C} = g\,. \]
Thus $\theta_{C,D}$ is injective. Moreover, given $h : FC \rarr D$, by the equational formulation of uniqueness~\eqref{univeq} we have:
\[ h=\theta_{C,D}(Gh\circ\eta_C)\,. \]
Thus $\theta_{C,D}$ is surjective. We are left to show naturality, \ie that the following diagram commutes, for all $h:C'\rarr C$ and $g:D\rarr D'$.
\[\xymatrix@C=13mm{%
\CC(C, GD)\ar[d]_{\theta_{C, D}}\ar[r]^{\CC(h,Gg)} & \CC(C', GD')\ar[d]^{\theta_{C', D'}} \\
\DD(FC, D)\ar[r]_{\DD(Fh, g)} & \DD(FC', D')} \]
We chase around the diagram, starting from $f:C\rarr GD$.
\[ \DD(Fh, g) \circ \theta_{C,D}(f) = g\circ\theta_{C,D}(f)\circ Fh \]
\[ \theta_{C',D'}\circ\CC(h, Gg)(f) = \theta_{C',D'}(Gg \circ f \circ h) \]
Now:
\begin{align*}
g\circ\theta_{C,D}(f)\circ Fh &= \theta_{C',D'}(G(g\circ\theta_{C,D}(f)\circ Fh)\circ\eta_{C'}) && \text{by (\ref{univeq})} \\
    &= \theta_{C',D'}(Gg\circ G(\theta_{C,D}(f))\circ GFh\circ\eta_{C'}) && \text{functoriality of $G$}\\
    &= \theta_{C',D'}(Gg\circ G(\theta_{C,D}(f))\circ \eta_C\circ h) && \text{naturality of $\eta$}\\
    &= \theta_{C',D'}(Gg\circ f\circ h) && \text{by (\ref{univeq}).}
\end{align*}
\qed[-1]

\begin{myproposition}[Adjunctions define universals]
Let $G:\DD\rarr\CC$ be a functor, $D\in Ob(\DD)$ and $C\in Ob(\CC)$. If, for any $D'\in Ob(\DD)$, there is a bijection
\[ \phi_{D'} : \CC(C, G(D'))\cong\DD(D, D')  \]
natural in $D'$ then there is a universal arrow $\eta:C\rarr G(D)$.
\end{myproposition}
\proof Take $\eta:C\rarr G(D) := \phi_D^{-1}(\id[D])$ and, for any $g:C\rarr G(D')$, take $\hat{g}:D\rarr D':=\phi_{D'}(g)$.\\
We have that
\[ G(\hat{g})\circ\eta= G(\hat{g})\circ\phi_D^{-1}(\id[D])\overset{\text{nat}}{=}\phi_{D'}^{-1}(\hat{g})=g\,. \]
Moreover, for any $h:D\rarr D'$,
\[ \phi_{D'}(Gh\circ\eta)=\phi_{D'}(Gh\circ\phi_D^{-1}(\id[D]))\overset{\text{nat}}{=}\phi_{D'}(\phi_{D'}^{-1}(h))=h\,, \]
where equalities labelled with ``nat" hold because of naturality of $\phi$. \qed

\begin{mycorollary}
Let $(F,G,\theta)$ be an adjunction with $F:\CC\rarr\DD$. Then, for each $C\in Ob(\CC)$ there is a universal arrow $\eta:C\rarr G(F(C))$. \qed
\end{mycorollary}

\paragraph{Equivalence of Universals and Adjoints}
Thus we see that the following two situations are equivalent, in the sense that each determines the other uniquely.
\begin{itemize}
\item We are given a functor $G : \DD \rarr \CC$, and for each object $C$ of $\CC$ a universal arrow from $C$ to $G$.
\item We are given functors $F : \CC \rarr \DD$ and $G : \DD \rarr \CC$, and a natural bijection
\[  \theta_{C,D} : \CC(C, G(D)) \cong \DD(F(C), D)\,.  \]
\end{itemize}

\paragraph{Couniversal Arrows}
Let $F : \CC \rarr \DD$ be a functor, and $D$ an object of $\DD$. A \boldemph{couniversal arrow from $F$ to $D$} is an object $C$ of $\CC$ and a  morphism
\[ \epsilon : F(C) \longrightarrow D \]
such that, for every object $C'$ of $\CC$ and morphism $g:F(C')\rarr D$, there exists a unique morphism $\bar{g}:C'\rarr C$ in $\CC$ such that
$g=\epsilon\circ F(\bar{g})$.\\
Diagrammatically:
\[
\xymatrix@=13mm{C & F(C)\ar[r]^{\epsilon} & D \\ C'\ar@{-->}[u]^{\bar{g}} & F(C')\ar@{-->}[u]^{F(\bar{g})}\ar[ur]_g}
\]
By exactly similar (but dual) reasoning to the previous propositions,  an adjunction implies the existence of couniversal arrows, and the
existence of the latter implies the existence of the adjunction. Hence,
\[ \text{Universality $\equiv$ Adjunctions $\equiv$ Couniversality}\,. \]
Some examples of couniversal arrows:
\begin{itemize}
\item A terminal object in a category $\CC$ is a couniversal arrow from the unique functor $K:\CC\rarr\1$ to the unique object in $\1$.
\item Let $A$, $B$ be objects of $\CC$. A product of $A$ and $B$ is a couniversal arrow from $\ang{\Id_\CC,\Id_\CC}:\CC\rarr\CC\times\CC$
to $(A,B)$.
\end{itemize}

\subsection{Limits and Colimits}
In the previous paragraph we described products $A\times B$ as couniversal arrows from the \emph{diagonal functor} $\Delta:\CC\rarr\CC\times\CC$ to
$(A,B)$. $\Delta$ is the functor assigning $(A,A)$ to each object $A$, and $(f,f)$ to each arrow $f$. Noting that $\CC\times\CC=\CC^{\2}$, where
$\CC^{\2}$ is a \emph{functor category}, this suggests an important generalisation.

\begin{mydefinition}
Let $\CC$ be a category and $\II$ be another category, thought of as an ``index category". A \emph{diagram of shape $\II$} in $\CC$ is just a functor $F:\II\rarr\CC$. Consider the functor category $\CC^{\II}$ with objects the functors from $\II$ to $\CC$, and
natural transformations as morphisms. There is a \emph{diagonal functor}
\[ \Delta : \CC \longrightarrow \CC^{\II} , \]
taking each object $C$ of $\CC$ to the constant functor $K_C:\II\rarr\CC$, which maps every object of $\II$ to $C$. A \boldemph{limit} for the diagram $F$ is a couniversal arrow from $\Delta$ to $F$. \deq
\end{mydefinition}
This concept of limit subsumes products (including
infinite products), pullbacks, inverse limits, etc.
\\
For example, take $\II:=\2_{\scriptscriptstyle\rightrightarrows}$  (we have seen this before:
$\2_{\scriptscriptstyle\rightrightarrows}=\xymatrix@1{\bullet\ar@/^/[r]\ar@/_/[r]&\bullet}$\ ). A functor $F$ from $\II$ to $\CC$ corresponds to a
diagram:
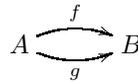
\[ \xymatrix@=10mm{A\ar@/^/[r]^f\ar@/_/[r]_g & B} \]
A couniversal arrow from $\Delta$ to $F$ corresponds to the following situation,
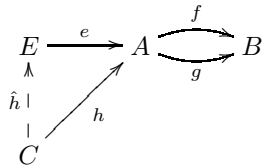
\[
\xymatrix@=10mm{E\ar[r]^e & A\ar@/^/[r]^f\ar@/_/[r]_g & B\\ C\ar@{-->}[u]^{\hat{h}}\ar[ur]_h}
\]
\ie to an equaliser!

By dualising limits we obtain \emph{colimits}. Some important examples are coproducts, coequalisers, pushouts and $\omega$-colimits.

\begin{myexercise}
Verify that pullbacks are limits by taking: \[ \II \; := \; \bullet \longrightarrow \bullet \longleftarrow \bullet\]
\end{myexercise}

\paragraph{Limits as terminal objects} Consider $\Delta : \CC \rarr\CC^{\II}$ and $F:\II\rarr\CC$.
A \emph{cone to  $F$} is an object $C$ of $\CC$ and family of arrows $\gamma$,
\[ \{\ \gamma_I:C\lrarr FI\ \}_{I\in Ob(\II)}\,, \]
such that, for any $f:I\rarr J$, the following triangle commutes.
\[ \xymatrix@=10mm{FI\ar[rr]^{Ff} && FJ\\ & C\ar[ul]^{\gamma_I}\ar[ur]_{\gamma_J}} \]
Thus a cone is exactly a natural transformation $\gamma : \Delta C \rarr F$.
A morphism of cones (`mediating morphism') $(C,\gamma)\lrarr(D,\delta)$ is an arrow $g:C\rarr D$ such that each of the following triangles commutes.
\[ \xymatrix@=10mm{& FI\\ C\ar[ur]^{\gamma_I}\ar[rr]_g && D\ar[ul]_{\delta_I}} \]
We obtain a category $\mathbf{Cone}(F)$ whose objects are cones to $F$ and whose arrows are mediating morphisms. Then, a {limit} of $F$ is a terminal
object in $\mathbf{Cone}(F)$.

\subsection{Exponentials}
In $\Set$, given sets $A$, $B$, we can form the set of functions $B^A:=\Set(A,B)$, \emph{which is again a set, \ie an object of $\Set$}. This closure of $\Set$ under forming
``function spaces" is one of its most important properties.

How can we axiomatise this situation? Once again, rather than asking what the elements of a function space \emph{are}, we ask  instead what  we can \emph{do} with them operationally. The answer is simple: apply functions to their arguments. That is, there is a map
\[ \ev{A}{B} : B^A \times A \longrightarrow B\ \text{ such that }\ \ev{A}{B}(f, a) = f(a)\,. \]
We can think of the function as a `black box': we can feed it inputs and observe the outputs.

Evaluation has the following couniversal property. For any $g : C \times A \rarr B$, there is a unique map $\Lambda(g):C\rarr B^A$ such that the
following diagram commutes.
\[
\xymatrix@=13mm{B^A\times A\ar[r]^-{\ev{A}{B}} & B\\ C\times A\ar[u]^{\Lambda(g)\times\id[A]}\ar[ur]_g}
\]
In $\Set$, this is defined by:
\[ \Lambda(g)(c): A \longrightarrow B := a\mapsto g(c, a)\,. \]
This process of transforming a function of two arguments into a function-valued function of one argument is known as \emph{currying}, after
H.~B.~Curry. It is an algebraic form of \emph{$\lambda$-abstraction}.

We are now led to the general definition of exponentials. Note that, for each object $A$ of a category $\CC$ with products, we can define a functor
\[ \uscore\times A : \CC \longrightarrow \CC\,. \]

\begin{mydefinition}\label{d:exp}
Let $\CC$ be a category with binary products. We say that $\CC$ \boldemph{has exponentials} if for all objects $A$ and $B$ of $\CC$ there is a
couniversal arrow from $\uscore\times A$ to $B$, {\ie}~an object $B^A$ of $\CC$ and a  morphism
\[ \ev{A}{B} : B^A \times A \longrightarrow B \]
with the couniversal property: for every $g : C \times A \rarr B$, there is a unique morphism $\Lambda (g) : C \rarr B^A$ such that the
following diagram commutes.
\[
\xymatrix@=13mm{B^A\times A\ar[r]^-{\ev{A}{B}} & B\\ C\times A\ar[u]^{\Lambda(g)\times\id[A]}\ar[ur]_g}
\]\deq[-1]
\end{mydefinition}
As before, the couniversal property can be given in purely equational terms, as follows. For every $h:C\rarr B^A$,
\[ \Lambda(\ev{A}{B}\circ h\times\id[A])=h\,.\]
Equivalently, $\CC$ has exponentials if, for every object $A$, the functor $\uscore\times A$ has a right adjoint, that is, there exists a functor
$\uscore^A:\CC\rarr\CC$ and a bijection
\[ \Lambda_{B,C}:\CC(C\times A,B)\overset{\cong}{\lrarr}\CC(C,B^A) \]
natural in $B,C$. In that case, $\ev{A}{B}:=\Lambda^{-1}(\id[B^A])$.

\begin{myexercise}
Derive $\uscore^A$ and $\Lambda^{-1}$ of the above description from $\mathsf{ev}$ and $\Lambda$ of definition~\ref{d:exp}.
\end{myexercise}
\begin{myexercise}\label{ex:expon}
Show that $\CC$ has exponentials iff, for every $A,B,C\in\Ob(\CC)$, there is an object $B^A$ and a bijection
\[ \theta_C:\CC(C\times A,B)\overset{\cong}{\lrarr}\CC(C,B^A) \]
natural in $C$.
\end{myexercise}

\begin{notation} The notation  $B^A$ for exponential objects is standard in the category theory literature. For our purposes, however, it will be more convenient to write $A \Rightarrow B$.
\end{notation}
Exponentials bring us to another fundamental notion, this time for understanding functional types, models of $\lambda$-calculus, and the structure of
proofs.

\begin{mydefinition}
A category with a terminal object, products and exponentials is called a \boldemph{Cartesian Closed Category (CCC)}. \deq
\end{mydefinition}
For example, $\Set$ is a CCC. Another class of examples are \emph{Boolean algebras}, seen as categories:
\begin{itemize}
  \item Products are given by conjunctions $A \wedge B$. We define exponentials as \emph{implications}:
    \[ A \Rightarrow B \; := \; \neg A \vee B\,. \]
  \item Evaluation is just \emph{Modus Ponens},
    \[ (A \Rightarrow B) \wedge A \;\leq\; B \]
    while couniversality is the \emph{Deduction Theorem},
    \[ C \wedge A \;\leq\; B \;\; \Longleftrightarrow \;\; C \;\leq\; A \Rightarrow B\,. \]
\end{itemize}

\subsection{Exercises}
\begin{enumerate}\renewcommand{\theenumi}{\textbf{\arabic{enumi}}}
  \item Suppose that $U:\CC\rarr\DD$ has a left adjoint $F_{1}$, and $V:\DD\rarr\EE$ has a left adjoint $F_{2}$. Show that
    $V\circ U:\CC\rarr \EE$ has a left adjoint.
  \item A \emph{sup-lattice} is a poset $P$ in which every subset $S \subseteq P$ has a supremum (least upper bound) $\bigvee S$.
  Let $P$, $Q$ be sup-lattices, and $f : P \rarr Q$ be a monotone map.
    \begin{enumerate}
    \item Show that if $f$ has a right adjoint then $f$ preserves least upper bounds:
    \[ f(\bigvee S) = \bigvee \{ f(x) \mid x \in S \}\,. \]
    \item Show that if $f$ preserves least upper bounds then it has a right adjoint $g$, given by:
    \[ g(y) = \bigvee \{ x \in P \mid f(x) \leq y \}\,. \]
    \item Dualise to get a necessary and sufficient condition for the existence of left adjoints.
    \end{enumerate}
  \item Let $F : \CC \rarr \DD$, $G : \DD \rarr \CC$ be functors such that $F$ is left adjoint to $G$, with natural bijection
    $\theta_{C,D}:\CC(C,GD)\overset{\cong}{\lrar}\DD(FC,D)$.
    Show that there is a natural transformation $\varepsilon : F \circ G \rarr \Id_{\DD}$, the \boldemph{counit} of the adjunction.
    \\
    Describe this counit explicitly in the case where the right adjoint is the forgetful functor $U : \Mon \rarr \Set$.
  \item Let $F:\CC\rarr\DD$ and $G:\DD\rarr\CC$ be functors, and assume $F$ is left adjoint to $G$ with natural bijection $\theta$.
    \begin{enumerate}
      \item Show that $F$ preserves epimorphisms.
      \item Show that $F$ is faithful if and only if, for every object $A$ of $\CC$, $\eta_{A}:A\rarr GF(A)$ is monic.
      \item Show that if, for each object $A$ of $\CC$, there is a morphism $s_{A}:GF(A)\rarr A$ such that
      $\eta_{A} \circ s_{A} = \id[GF(A)]$ then $F$ is full.
    \end{enumerate}
\end{enumerate}

\section{The Curry-Howard Correspondence}
We shall now study a beautiful three-way connection between logic, computation and categories:
\begin{center}
\rnode{L}{\psshadowbox{\textbf{Logic}}} \hspace{30mm}%
\rnode{C}{\psshadowbox{\textbf{Computation}}} \\[10mm]
\rnode{X}{\psshadowbox{\textbf{Categories}}}\hspace{13mm} %
\ncline[linewidth=.6mm,nodesep=1pt]{<->}{L}{C}\ncline[linewidth=.6mm,nodesep=1pt]{<->}{L}{X}\ncline[linewidth=.6mm,nodesep=1pt]{<->}{C}{X}
\captionof{table}{The Curry-Howard Correspondence.}
\end{center}
%
This connection has been known since the 1970's, and is widely used in Computer Science\,---\,it is also beginning to be used in Quantum Informatics!
It is the upper link (Logic\,--\,Computation) that is usually attributed to Haskell B.~Curry and William A.~Howard, although the idea is related to  the  operational interpretation of intuitionistic logic  given in various formulations by Brouwer, Heyting and Kolmogorov. The link to Categories is mainly due to the pioneering
work of Joachim Lambek.

\subsection{Logic}
Suppose we ask ourselves the question: \emph{What is Logic about?} There are  two main kinds of answer:  one focuses on \emph{Truth}, and the other
on \emph{Proof}. We focus on the latter, that is, on:
\[ \fbox{What follows from what} \]
Traditional introductions to logic focus on Hilbert-style proof systems, that is, on generating the set of \emph{theorems} of a system from a set of
\emph{axioms} by applying rules of inference (e.g.~Modus Ponens).

A key step in logic took place in the 1930's with the advent of \emph{Gentzen-style systems}. Instead of focusing on theorems, we look more generally
and symmetrically at \emph{What follows from what}: in these systems the primary focus is on \emph{proofs from assumptions}.
We will examine two such kinds of systems: Natural Deduction systems and Gentzen sequent calculi.

\begin{mydefinition}
Consider the fragment of propositional logic with logical connectives $\wedge$ and $\supset$. The assertion that  a formula $A$ can be proved from {assumptions}
$A_{1},...,A_{n}$ is expressed by a \boldemph{sequent}:
\[ A_{1}, \ldots , A_{n} \vdash A \]
We use $\Gamma$, $\Delta$ to range over finite sets of formulas, and write $\Gamma,A$ for $\Gamma\cup\{A\}$. Proofs are built using the proof rules of table~\ref{t:NatDed}; the resulting proof system is called \emph{the} \boldemph{Natural Deduction system for} $\wedge$,$\supset$. \deq
\end{mydefinition}

\begin{center}
\renewcommand{\arraystretch}{.5}
\fbox{$\begin{array}{@{\;\;}c@{\;\;}|@{\;\;}c@{\;\;}|@{\;\;}c@{\;\;}} &&\\ \textbf{Identity} & \textbf{Conjunction} &
\textbf{Implication}
\\&&\\\hline\hline&&\\
\infer[\Ax]{\Gamma , A \vdash  A}{} & %
\infer[\someQWEqwe{\land}{intro}]{\Gamma \vdash A \wedge B}{\Gamma \vdash A \qquad \Gamma \vdash B} &%
\infer[\someQWEqwe{\mathnormal\supset}{intro}]{\Gamma \vdash A \supset B}{\Gamma , A \vdash B}
\\&&\\&&\\&
\infer[\someQWEqwe{\land}{elim_1}]{\Gamma \vdash A}{\Gamma \vdash A \wedge B} &%
\infer[\someQWEqwe{\mathnormal\supset}{elim}]{\Gamma \vdash B}{\Gamma \vdash A \supset B \qquad \Gamma \vdash A}
\\&&\\&&\\&
\infer[\someQWEqwe{\land}{elim_2}]{\Gamma \vdash B}{\Gamma \vdash A \wedge B}\\&&\\
\end{array}$} \captionof{table}{Natural Deduction System for $\wedge$,$\supset$.}\label{t:NatDed}
\end{center}
For example, we have the following proof of $\supset$-transitivity.
\[\hspace{-24mm}\text{\small
\infer[\ImpI]{A\supset B,B\supset C\vdash A\supset C}{\infer[\ImpE]{A\supset B,B\supset C,A\vdash C}{
    \hspace{23mm}\infer[\Ax]{A\supset B,B\supset C,A\vdash B\supset C}{}\
    \infer[\ImpE]{A\supset B,B\supset C,A\vdash B}{\infer[\Ax]{A\supset B,B\supset C,A\vdash A\supset B}{}\quad
    \infer[\Ax]{A\supset B,B\supset C,A\vdash A}{}}}}}
\]
An important feature of Natural Deduction is the systematic pattern it exhibits in the structure of the inference rules. For each connective $\Box$ there are \emph{introduction rules}, which show how formulas $A \Box B$ can be derived, and \emph{elimination rules}, which show how such formulas can be used to derive other formulas.

\paragraph{Admissibility}
We say that a proof rule
\[ \infer{\Delta\vdash B}{\Gamma_1\vdash A_1\quad\cdots\quad\Gamma_n\vdash A_n} \]
is \boldemph{admissible} in Natural Deduction if, whenever there are proofs of $\Gamma_i\vdash A_i$ then there is also a proof of $\Delta\vdash B$. For example, the following \emph{Cut rule} is admissible.
\[ \infer[\Cut]{\Gamma,\Delta\vdash B}{\Gamma\vdash A\quad A,\Delta\vdash B} \]

\begin{myexercise}
Show that the following rules are admissible in Natural Deduction.
\begin{enumerate}
\item The Weakening rule:
\[ \infer{\Gamma,A\vdash B}{\Gamma\vdash B} \]
\item The Cut rule.
\end{enumerate}
\end{myexercise}
Our focus will be on \boldemph{Structural Proof Theory}, that is studying the ``space of formal proofs" as a mathematical structure in its own right,
rather than  focussing only on
\begin{center}
Provability $\longleftrightarrow$ Truth
\end{center}
(\ie the usual notions of ``soundness and completeness"). One motivation for this approach comes from trying to understand and use the
\emph{computational content} of proofs, epitomised in the ``Curry-Howard correspondence".

\subsection{Computation}
Our starting point in computation is the pure calculus of functions called the \emph{$\lambda$-calculus}.

\begin{mydefinition}
Assume a countably infinite set of \boldemph{variables}, ranged over by $x,y,z$ and variants.
\boldemph{$\lambda$-calculus terms}, ranged over by $t,u,v$ etc, are constructed from the following inductive definition.
\begin{itemize}
  \item Every variable $x$ is a term.
  \item If $t$ and $u$ are terms, then $t\,u$ is a term (\boldemph{application}).
  \item If $x$ is a variable and $t$ is a term, then $\lambda x.\,t$ is a term (\boldemph{$\lambda$-abstraction}). \deq
\end{itemize}
\end{mydefinition}
The above definition can be given in the following compact form, which will be followed in similar definitions in the sequel.
\begin{align*}
    \mathrm{VA} &\ni x,y,z,\dots \\
    \mathrm{TE} &\ni t,u,v \;\; ::= \;\; x \; \mid \; {t\,u} \; \mid \; {\lambda x . \, t}
\end{align*}
The computational content of the calculus is exhibited in the following examples.
Note that the first example is not part of our formal syntax: it presupposes some \emph{encoding} of numerals and successors.
\begin{center}
\begin{tabular}{l@{\quad}l}
$\lambda x. \, x+1$ & successor function \\
$\lambda x . \, x$ & identity function \\
$ \lambda f . \, \lambda x. \, fx$ & application \\
$\lambda f. \, \lambda x. \, f(fx)$ & double application \\
$\lambda f. \, \lambda g. \, \lambda x. \, g(f(x))$ & composition and application
\end{tabular}\captionof{table}{Examples of $\lambda$-terms.}
\end{center}
What we also note above is the use of parentheses in order to disambiguate the structure of terms (\ie the precedence of term constructors). To avoid notational clutter we also use the following conventions.
\begin{itemize}
  \item Applications associate to the left. For example, $f\,x\,y$ stands for $(fx)\,y$\,.
  \item The scope of an abstractions goes as far to the right as possible. For example,
  \[ \lambda f.(\lambda x.f(xx))\,\lambda x.f(xx)\text{ \ stands for \ }\lambda f.((\lambda x.(f(xx)))(\lambda x.(f(xx))))\,. \]
\end{itemize}
The \emph{free variables} of a term are those that are not bound by any $\lambda$; they can be seen as the \emph{assumptions} of the term.

\begin{mydefinition}\label{d:fv}
The set of \boldemph{free variables} of a term $t$, $\fv(t)$, is given by:
\begin{align*}
\fv(x)  &:= \{ x \}\,, \\
\fv(t\,u) &:= \fv(t)\cup\fv(u)\,, \\
\fv(\lambda x.t) &:= \fv(t)\setminus\{x\}\,.
\end{align*}\deq[-1]
\end{mydefinition}
The notation $\lambda x.t$ is meant to serve the purpose of expressing formally
\[ \emph{the function that returns $t$ on input $x$.} \]
Thus, $\lambda$ is a \emph{binder}, that is, it binds the variable $x$ in the `function' $\lambda x.t$, in the same way that e.g.~$\int$ binds $x$ in
$\int f(x)\,dx$\,.
%
This means that there should be no difference between $\lambda x.t$ and $\lambda x'.t'$, where $t'$ is obtained from $t$ by swapping $x$ with some
\emph{fresh} variable $x'$ (\ie with some $x'$ not appearing free in $t$). For example, the terms
\[ \lambda x.x \text{\;and\;} \lambda x'.x' \]
should be `equal', as they both stand for the identity function. We formalise this by stipulating that
\[\fbox{\text{Terms are identified up to $\alpha$-equivalence}} \]
where we say that two terms are $\alpha$-equivalent iff they differ solely in the choice of variables appearing in \emph{binding positions}.
This is formally defined in two steps, as follows.

\begin{mydefinition}
We define \emph{variable-swapping} on terms recursively as follows.
\begin{align*}
\sw{z}{y}{x}  \ &:=\ \begin{cases}y & \text{if }z=x \\ x & \text{if }z=y \\ z & \text{otherwise}\end{cases} \\
\sw{t\,u}{y}{x} \ &:=\ (\sw{t}{y}{x})(\sw{u}{y}{x}) \\
\sw{\lambda z. t}{y}{x} \ &:=\ \lambda(\sw{z}{y}{x}).(\sw{t}{y}{x})
\end{align*}
Then, \boldemph{$\alpha$-equivalence}, $=_\alpha$, is the relation on terms defined inductively by:\footnote{%
The last clause can be replaced by any of the following:
\begin{itemize}
    \item \dots if, for some $y$ not appearing in $t\,t'$, $\sw{t}{y}{x}=_\alpha \sw{t'}{y}{x'}$\,.
    \item \dots if, for all $y$ not appearing free in $t\,t'$, $\sw{t}{y}{x}=_\alpha \sw{t'}{y}{x'}$\,.
    \item \dots if, for some $y$ not appearing free in $t\,t'$, $\sw{t}{y}{x}=_\alpha \sw{t'}{y}{x'}$\,.
\end{itemize}}
\begin{itemize}
  \item $x=_{\alpha}x$,
  \item $t\,u=_{\alpha}t'\,u'$ if $t=_{\alpha}t'$ and $u=_{\alpha}u'$,
  \item $\lambda x.t=_\alpha\lambda x'.t'$ if,  for all $y$ not appearing in $t\,t'$, $\sw{t}{y}{x}=_\alpha \sw{t'}{y}{x'}$\,. \deq
\end{itemize}
\end{mydefinition}
Equating terms modulo $\alpha$-equivalence means that we work with $\mathrm{TE}/\mathnormal{=_\alpha}$ instead of $\mathrm{TE}$.
Henceforth, we will refer to elements of $\mathrm{TE}/\mathnormal{=_\alpha}$ as \emph{terms}, and to elements of $\mathrm{TE}$ as \boldemph{raw terms}. Note that $\alpha$-equivalence is meaningful only on raw terms.
\begin{myexercise}
Prove the following $\alpha$-equivalences.
\[ \lambda x.x=_\alpha\lambda y.y\,,\quad \lambda x.\lambda y.\,xy=_\alpha\lambda y.\lambda x.\,yx\,,\quad x(\lambda x.x)=_\alpha x(\lambda y.y)\,. \]
\end{myexercise}
\begin{myexercise}
Show that, for all raw terms $t,t'$ and variables $x,x'$, if $t=_\alpha t'$ then $\fv(t)=\fv(t')$ and $\sw{t}{x}{x'}=_\alpha\sw{t'}{x}{x'}$.  \\
Moreover, show that, for any $x,x'\not\in\fv(t)$, \ $t=_\alpha\sw{t}{x}{x'}$\,. Hence infer that, for any $x'\not\in\fv(t)$, \ $\lambda
x.\,t=_\alpha\lambda x'.\sw{t}{x}{x'}$\,.
\end{myexercise}
From the above exercise we obtain that $\fv$ and variable-swapping extend to terms (\ie to $\mathrm{TE}/\mathnormal{=_\alpha}$) in a straightforward manner.
Moreover, we have that, for any term $t$ and any $x'\notin\fv(t)$,
\[ \lambda x.\,t= \lambda x'.\sw{t}{x}{x'}\,. \]
Since $\lambda$-abstractions stand for functions, an application of a $\lambda$-abstraction on
another term should result to a \emph{substitution} of the latter inside the body of the abstraction.
\begin{mydefinition}
Define the \boldemph{substitution} of a term $t$ for a variable $x$ inside a term inductively by:
\begin{align*}
y[t/x]  &:=\begin{cases}t & \text{if }y=x \\ y & \text{if }y \neq x\end{cases} \\
(uv)[t/x] &:= (u[t/x])(v[t/x]) \\
(\lambda z. u) [t/x] &:= \lambda z.\, (u[t/x]) \qquad (\ast )
\end{align*}
where  ($*$) indicates the condition that  $z\not\in\fv(x\,t)$.\deq
\end{mydefinition}
Note that, due to identification of $\alpha$-equivalent (raw) terms, it is always possible to rename bound variables so that condition ($*$) be satisfied: for example,
\[ (\lambda z. zx) [z/x] = (\lambda y. yx) [z/x] = \lambda y. yz \]
\begin{myexercise}\label{ex:SubSub}
Show that, for all $\lambda$-terms $u,t,t'$ and variables $x,x'$ such that $x'\notin\fv(u)\setminus\{x\}$,
    \[ u[t/x][t'/x'] = u[(t[t'/x'])/x]\,. \]
\end{myexercise}
%
%
%
%
We proceed to the definition of $\beta$-reduction and $\beta$-conversion. These are relations defined on pairs of terms and express the
computational content of the calculus.
\begin{mydefinition}
We take \boldemph{$\beta$-reduction}, $\longrightarrow_\beta$, to be the relation defined by: 
\[ (\lambda x. t)\,u \longrightarrow_\beta t[u/x]\,. \]
This extends to arbitrary terms as follows. If $t\brar t'$ then:
\[ t\,u\brar t'u\,,\quad u\,t\brar u\,t',\quad \lambda x.t\brar\lambda x.t'.\]
We take \boldemph{$\beta$-conversion}, $=_\beta$, to be the symmetric reflexive transitive closure of $\beta$-reduction, that is, the
equivalence relation induced by:
\[ (\lambda x.t)\,u =_\beta t[u/x]\,. \]\deq[-1]
\end{mydefinition}
With $\beta$-reduction we obtain a notion of ``computational dynamics". For example:
\begin{align*}
(\lambda f.\,f(f\,y))(\lambda x.\,x+1) &\brar (\lambda x.\,x+1)((\lambda x.\,x+1)\,y) \\
    &\brar((\lambda x.\,x+1)\,y)+1\brar(y+1)+1 \\[5pt]
(\lambda f.\,f(f\,y))(\lambda x.\,x+1) &\brar (\lambda x.\,x+1)((\lambda x.\,x+1)\,y) \\
    &\brar(\lambda x.\,x+1)(y+1)\brar(y+1)+1
\end{align*}
Note that in the sequel we will usually write $\beta$-reduction simply by ``$\lrarr$".

\subsection{Simply-Typed $\lambda$-calculus}

The `pure' $\lambda$-calculus we have discussed so far is \emph{very} unconstrained. For example, it allows \emph{self-application}, \ie terms like $xx$ are perfectly legal. On the one hand, this means that the calculus very expressive: for example, we can encode \emph{recursion} by setting
\[ \YY := \lambda f.  (\lambda x. \, f(x x))\,\lambda x.  f(x x)\,. \]
We have:
\[ \YY t \lrarr (\lambda x. \, t(x x)) \,\lambda x. \, t(x x) \lrarr t ((\lambda x. \, t(x x))\,\lambda x. \, t(x x)) \llarr t(\YY t) \]
However, self-application leads also to \emph{divergences}. The most characteristic example is the following. Setting $\Omega :=(\lambda x.xx)\,\lambda x.xx$, we have:
\[ \Omega \longrightarrow \Omega \longrightarrow \Omega \longrightarrow \cdots \]
Historically, Curry extracted $\YY$ from an analysis of Russell's Paradox, so it should come as no surprise that it too leads to divergences:
setting $t'$ to be $\lambda x. \, t(x x)$,
\[ \YY t \lrarr t't' \lrarr t(t't')\lrarr t(t(t't')) \lrarr \cdots \]
The solution is to introduce \emph{types}.  The original idea, due to Church following Russell, was that:
\[ \fbox{\textit{Types are there to stop you doing bad things}} \]
However, it has turned out that types constitute one of the most fruitful \emph{positive} ideas in Computer Science, and provide  one of the
key disciplines of programming.
\begin{mydefinition}
Let us assume a set of \emph{base types}, ranged over by $b$. The \boldemph{simply-typed $\lambda$-calculus} is defined as follows.
\begin{align*}
&\text{Type} &&
\mathrm{TY}\ni\, T,U \; ::= \; b \; \mid \; T \rightarrow U \; \mid \; T \times U 
\\[1.5mm]
&\text{Term} && \mathrm{TE}\ni\,t,u\; ::= \; x \; \mid \; {t\,u} \;\mid\; {\lambda x.\, t} \;\mid\; \ang{t,u} \;\mid\; \pi_1u \;\mid\; \pi_2u \\[1.5mm]
&\text{Typing context} && \Gamma \;::=\; \varnothing\; \mid x : T,\Gamma \qquad(x\text{ does not appear in }\Gamma)
\end{align*}
A \boldemph{typing judgement} is a triple of the form
\[ \Gamma \vdash t : T\,, \]
which is to be understood as the assertion that term $t$ has the type $T$ under the assumptions that $x_{1}$ has type $T_{1}$, \ldots , $x_{k}$ has type $T_{k}$, if $\Gamma = x_{1}:T_{1}, \ldots , x_{k} : T_{k}$. A \boldemph{typed term} is a term $t$ accompanied with a type $T$ and a context $\Gamma$, such that the judgement $\Gamma\vdash t:T$ is derivable by use of the \emph{typing rules} of table~\ref{t:SimTyp}.
\deq
\end{mydefinition}

\begin{center}
\renewcommand{\arraystretch}{.5}
\fbox{$\begin{array}{@{\;\;}c@{\;\;}|@{\;\;}c@{\;\;}|@{\;\;}c@{\;\;}} &&\\ \textbf{Variable} & \textbf{Product} &
\textbf{Function}
\\&&\\\hline\hline&&\\
\infer{\Gamma , x:T \vdash  x:T}{} & %
\infer{\Gamma \vdash \langle t, u \rangle : T \times U}{\Gamma \vdash t:T \qquad  \Gamma \vdash u:U} &%
\infer{\Gamma \vdash \lambda x. \, t: U \rightarrow T}{\Gamma , x:U \vdash t:T}
\\&&\\&&\\&
\infer{\Gamma \vdash \pi_1 v : T}{\Gamma \vdash v : T \times U} &%
\infer{\Gamma \vdash t\,u : T}{\Gamma \vdash t : U \rightarrow T \qquad \Gamma  \vdash u:U}
\\&&\\&&\\&
\infer{\Gamma \vdash \pi_2 v : U}{\Gamma \vdash v : T \times U}
\\&&\\
\end{array}$}\captionof{table}{Simply-Typed $\lambda$-calculus for $\times,\rarr$.}\label{t:SimTyp}
\end{center}%
Note that contexts are sets, and so $x : T,\Gamma$ stands for $\{x:T\}\cup\Gamma$ with $x$ not appearing in $\Gamma$. As before, terms are identified up to $\alpha$-equivalence.

From the definition of types we see that the simply-typed $\lambda$-calculus is a calculus of functions and products.
For example:
\[\begin{array}{c@{\quad}l}
 b \rightarrow b \rightarrow b & \mbox{first-order function type} \\[5pt]
(b \rightarrow b ) \rightarrow b & \mbox{second-order function type}
\end{array}\]


%
\begin{myexercise}
Can you type the following terms?
\[ \lambda x. \, xx\,,\quad \lambda f. \, (\lambda x.\, f(x x)) (\lambda x. \, f(x x))\,. \]
\end{myexercise}
\begin{myexercise}[Weakening \& Cut] Show that Weakening and Cut are admissible in the typing system of the simply-typed $\lambda$-calculus:
    \[ \infer[\Weaks]{\Gamma,x:U\vdash t:T}{\Gamma\vdash t:T}\qquad
    \infer[\Cut]{\Gamma\vdash u[t/x]:U}{\Gamma\vdash t:T\qd[2]\Gamma,x:T\vdash u:U} \]
\end{myexercise}
%
\newlines{1}
We proceed to the rules for reduction and conversion. These are given as in the untyped case, with the addition of $\eta$-rules, which are essentially
extensionality principles.
\begin{mydefinition}
We define $\beta$-reduction, $\longrightarrow_\beta$, by the following rules, and let $\beta$-conversion, $=_\beta$, be its symmetric reflexive
transitive closure.
\[ \begin{array}{lcl}
(\lambda x. \, t ) u & \longrightarrow_\beta & t[u/x] \\
\pi_1 \langle t, u \rangle & \longrightarrow_\beta & t \\
\pi_2 \langle t, u \rangle & \longrightarrow_\beta & u
\end{array} \]
Moreover, $\eta$-conversion, $=_\eta$, is the symmetric reflexive transitive relation obtained by the following rules,
\[ \begin{array}{lcl@{\quad}l}
t & =_\eta & \lambda x.\, tx & \mbox{$x\not\in\fv(t)$, at function  types} \\
v & =_\eta & \langle \pi_1 v, \pi_2 v \rangle & \mbox{at product types}
\end{array} \]
and $\lambda$-conversion, $=_\lambda$, is the transitive closure of $=_\beta\cup=_\eta$\,.\deq
\end{mydefinition}
Implicit in the above definition is the fact that $\eta$-rules relate typed terms. For example, $t=_\eta\lambda x.\,tx$ has as side condition that $t$ be of function type, \ie that $t$ be a typed term $\Gamma\vdash t:T\rarr U$.
Now, following our intuitive interpretation of arrows as functions, we can read this $\eta$-rule as:
\[ \textit{$t$ is the function returning $t(x)$ to every input $x$}\]
Note that the above statement is in fact the couniversal property of currying in $\Set$; we will see more on this in the next sections!
\begin{myexercise}[Subject Reduction]
    Show that, for any typed term $\Gamma\vdash t:T$, if $t\brar t'$ then $\Gamma\vdash t':T$ is derivable.
\end{myexercise}

\paragraph{Strong Normalisation} Term reduction results in a \emph{normal form}: an \emph{explicit but much longer expression} in which
no more reductions are applicable.
Formally, a $\lambda$-term is called a \boldemph{redex} if it is in one of forms of the left-hand-side of the $\beta$-reduction rules,
and therefore $\beta$-reduction can be applied to it. A term is in \boldemph{normal form} if it contains no redexes.
In the light of the correspondence presented in the next paragraph, a term in normal form corresponds to a proof in which all lemmas have been eliminated.
\begin{fact}[SN]
For every term $t$, there is no infinite sequence of $\beta$-reductions:
\[ t\lrarr t_0\lrarr t_1\lrarr t_2\lrarr \cdots \]
\end{fact}
The above result states that \emph{every} reduction sequence leads eventually to a term in normal form. Note, though,
that reduction to normal form has enormous (\emph{non-elementary}) complexity.

\paragraph{The correspondence between Logic and Computation} Comparing the following two systems,
\[
  \fbox{Natural Deduction System for $\wedge$,$\supset$}\text{ \ vs \ }
  \fbox{Simply-Typed $\lambda$-calculus for $\times$,$\rightarrow$}
\]
we notice that if we equate
\[ \begin{array}{lcl}
\wedge & \equiv & \times \\
\supset & \equiv & \rightarrow
\end{array} \]
then they are the same! This is the Logic--Computation part of the Curry-Howard correspondence (sometimes: ``Curry-Howard isomorphism"). It works on three levels:
\begin{center}
\fbox{\renewcommand{\arraystretch}{1.3}\begin{tabular}{@{\;\,}l@{\quad}|@{\;\,}l@{\quad}}
Natural Deduction System & Simply-Typed $\lambda$-calculus\\ \hline
Formulas & Types \\
Proofs & Terms \\
Proof transformations & Term reductions
\end{tabular}}\captionof{table}{Correspondence between Logic and Computation.}
\end{center}
%
%
%
The view of proofs as containing computational content can also be detected in the \emph{Brouwer-Heyting-Kolmogorov interpretation} of intuitionistic
logic:
\begin{itemize}
\item A proof of an implication $A \supset B$ is a procedure which transforms any proof of $A$ into a proof of $B$.
\item A proof of $A \wedge B$ is a pair consisting of a proof of $A$  and a proof of $B$.
\end{itemize}
These readings motivate identifying $A \wedge B$ with $A \times B$, and $A \supset B$ with $A \rightarrow B$.
Moreover, these ideas have strong connections to computing. The $\lambda$-calculus is a `pure' version of functional programming languages such as
Haskell and SML. So we get a reading of:
\begin{center}
\fbox{{\textit{Proofs as Programs}}}
\end{center}

\subsection{Categories}
We now have our link between Logic and Computation. We now proceed to complete the triangle of the Curry-Howard correspondence by showing the connection to Categories.

We establish the link from Logic (and Computation) to Categories. Let $\CC$ be a cartesian closed category. We shall interpret formulas (or types) as
objects of $\CC$. A morphism $f:A\rarr B$ will then correspond to a proof of $B$ from assumption $A$, {\ie}~a proof of $A\vdash B$ (a typed term
$x:A\vdash t:B$).
Note that the bare structure of a category only supports proofs from a single assumption. Since $\CC$ has finite products, a proof of
\[ A_1 , \ldots , A_k \vdash A \]
will correspond to a morphism
\[ f : A_1 \times \cdots \times A_k \longrightarrow A\,. \]
The correspondence is depicted in the next table.
\begin{center}
\renewcommand{\arraystretch}{0.5}\fbox{$\begin{array}{c||@{\;\;}c@{\;\;}|@{\;\;}c@{\;\;}} &&\\
\textbf{Axiom} & \infer[\Ax]{\Gamma,A\vdash A}{} & \infer{\pi_{2}:\Gamma\times A\longrightarrow A}{}
\\&&\\\hline&&\\
\textbf{Conjunction} & \infer[\AndI]{\Gamma\vdash A\wedge B}{\Gamma\vdash A\qquad\Gamma\vdash B} & \infer{\langle f,g\rangle:\Gamma\longrightarrow A
\times B}{f : \Gamma \longrightarrow A \qquad g : \Gamma \longrightarrow B}
\\&&\\
& \infer[\AndEl]{\Gamma \vdash A}{\Gamma \vdash A \wedge B} & \infer{\pi_{1} \circ f : \Gamma \longrightarrow A}{f : \Gamma \longrightarrow A \times B}
\\&&\\
& \infer[\AndEr]{\Gamma \vdash B}{\Gamma \vdash A \wedge B} & \infer{\pi_{2}\circ f : \Gamma \longrightarrow B}{f : \Gamma \longrightarrow A \times B}
\\&&\\\hline&&\\
\textbf{Implication} & \infer[\ImpI]{\Gamma\vdash A\supset B}{\Gamma,A\vdash B} & \infer{\Lambda(f):\Gamma\longrightarrow(A\Rightarrow B)}{f:\Gamma
\times A\longrightarrow B} \\&&\\
& \infer[\ImpE]{\Gamma\vdash B}{\Gamma\vdash A\supset B \quad\Gamma\vdash A} & \infer{\App[A,B]\circ\langle f,g\rangle:\Gamma\longrightarrow B}{f :
\Gamma \longrightarrow (A \Rightarrow B) \quad g : \Gamma \rarr A}\\&&\\
\end{array}$}\captionof{table}{Correspondence between Logic and Categories.}
\end{center}
Moreover, the rules for $\beta$- and $\eta$-conversion are all then \emph{derivable} from the equations of cartesian closed categories. So cartesian
closed categories are \emph{models} of $\wedge$,$\supset$-logic at the level of proofs and proof-transformations, and of simply typed
$\lambda$-calculus at the level of terms and term-conversions.
The connection to computation is examined in more detail below.

\begin{myremark}
In our translation of Logic sequents there is an implicit ordering of assumptions: a \emph{set} of assumptions is mapped to an assumption \emph{product},
\[ \{A_1,\dots,A_n\} \longmapsto A_1\times\cdots\times A_n\,. \]
In practice, since for any permutation $A_1',\dots,A_n'$ of $A_1,\dots,A_n$ we have
\[ A_1\times\cdots\times A_n \cong A_1'\times\cdots\times A_n'\,,\]
such an ordering is harmless.
\end{myremark}

\subsection{Categorical Semantics of Simply-Typed $\lambda$-calculus}

We translate the simply-typed $\lambda$-calculus into a cartesian closed category $\CC$, so that to each typed term $x_1:T_1,...,x_k:T_k\vdash t:T$
corresponds an arrow
\[ \trn{t}:\trn{T_1}\times\cdots\times\trn{T_k}\lrarr\trn{T}\,. \]
The translation if given by the function $\trn{\uscore\,}$ defined below (``semantic brackets").

\begin{mydefinition}[{Semantic translation}]\label{d:SemTr}
Let $\CC$ be a CCC and suppose we are given an assignment of an object $\tilde{b}$ to each base type $b$. Then, the translation is defined recursively
on types by:
\[ \trn{b}:=\tilde{b} \,,\qd \trn{T\times U}:=\trn{T}\times\trn{U} \,,\qd \trn{T\rarr U}:=\trn{T}\Rightarrow\trn{U}\,, \]
and on typed terms by:
\[\begin{array}{c}
\infer{\trn{\Gamma,x:T\vdash x:T} := \pi_2:\trn{\Gamma}\times\trn{T}\lrarr\trn{T}}{}
\\[3mm]
\infer{\trn{\Gamma\vdash\pi_1 t:T}:=\trn{\Gamma}\overset{f}{\lrar}\trn{T}\times\trn{U}\overset{\pi_1}{\lrar}\trn{T}}{%
    \trn{\Gamma\vdash t:T\times U}=f:\trn{\Gamma}\lrar\trn{T}\times\trn{U}}
\\[3mm]
\infer{\trn{\Gamma\vdash \ang{t,u}:T\times U}:=\trn{\Gamma}\overset{\ang{f,g}}{\lrar}\trn{T}\times\trn{U}}{%
    \trn{\Gamma\vdash t:T}=f:\trn{\Gamma}\lrar\trn{T}\qd[2]\trn{\Gamma\vdash u:U}=g:\trn{\Gamma}\lrar\trn{U}}
\\[3mm]
\infer{\trn{\Gamma\vdash\lambda x.\,t:T\rarr U}:=\Lambda(f):\trn{\Gamma}\lrarr\funsp{\trn{T}}{\trn{U}}}{%
    \trn{\Gamma,x:T\vdash t:U}=f:\trn{\Gamma}\times\trn{T}\lrarr\trn{U}}
\\[3mm]
\infer{\trn{\Gamma\vdash t\,u: U}:=\trn{\Gamma}\overset{\ang{f,g}}{\lrarr}\funsp{\trn{T}}{\trn{U}}\times\trn{T}
    \overset{\App}{\lrarr}\trn{U}}{%
        \trn{\Gamma\vdash t:T\rarr U}=f\qd[2]\trn{\Gamma\vdash u:T}=g}
\end{array}\]\deq[-1]
\end{mydefinition}
Our aim now is to verify that $\lambda$-conversion (induced by~$\beta$- and $\eta$-rules) is \emph{preserved} by the translation, \ie that, for any
$t,u$,
\[ t=_\lambda u \implies \trn{t}=\trn{u}\,. \]
This would mean that our categorical semantics is \emph{sound}.

Let us recall some structures from CCC's. Given $f_{1} : D_{1} \rarr E_{1}$,   $f_{2} : D_{2} \rarr E_{2}$, we defined
\[ f_{1} \times f_{2}  = \langle f_{1} \circ \pi_{1},  f_{2} \circ \pi_{2} \rangle : D_{1} \times D_{2} \lrarr E_{1} \times E_{2}\,, \]
and we showed that \ $(f_1\times f_2)\circ\langle h_1,h_2\rangle = \langle f_1\circ h_1, f_2\circ h_2\rangle$\,.
Moreover, exponentials are given by the following natural bijection.
\[ \frac{f : D \times E \lrarr F}{\Curry{f} : D \lrarr \funsp{E}{F}} \]
Equivalently, recall the basic equation:
\[ \App \circ (\Curry{f} \times \id[E] ) = f\,, \]
where $\Curry{f}$ is the unique arrow $D\rarr\funsp{E}{F}$ satisfying this equation, with uniqueness being specified by:
\[ \forall h : D \lrarr \funsp{E}{F}. \; \Curry{\App \circ (h \times \id[E])} = h\,. \]
Naturality of $\Lambda$ is then proven as follows.

\begin{myproposition}
\label{p:Curry}
For any $f:A\times B\rarr C$ and $g:A'\rarr A$\,,
\[ \Curry{f} \circ g = \Curry{f \circ (g \times \id[B])}\,. \]
\end{myproposition}
\proof
\[ \begin{array}{rcl}
 \Curry{f} \circ g & = & \Curry{\App \circ ((\Curry{f} \circ g) \times \id[B])} \\
 & = & \Curry{\App \circ (\Curry{f}\times \id[B]) \circ (g \times \id[B])) } = \Curry{f \circ (g \times \id[B])}\,.
\end{array}\]
\qed[-1]

\paragraph{Substitution Lemma} We consider a \boldemph{simultaneous substitution} for all the free variables in a term.
\begin{mydefinition}
Let $\Gamma=x_{1}:T_{1}\lldots x_{k}:T_{k}$\,. Given typed terms
\[ \Gamma \vdash t : T  \text{\qd and\qd} \Gamma \vdash t_{i} : T_{i}\,,\; 1 \leq i \leq k\,,\]
we define \ $t[\vec{t}/\vec{x}] \equiv t[t_{1}/x_{1}\lldots t_{k}/x_{k}]$ \ recursively by:
\[\begin{array}{rcl}
x_{i}[\vec{t}/\vec{x}]  & := & t_{i}\\[1.1mm]
(\pi_i\,t)[\vec{t}/\vec{x}] & := & \pi_i(t[\vec{t}/\vec{x}]) \\[1.1mm]
\ang{t,u}[\vec{t}/\vec{x}] & := & \ang{t[\vec{t}/\vec{x}],u[\vec{t}/\vec{x}]} \\[1.1mm]
(t\,u)[\vec{t}/\vec{x}]  & := & (t[\vec{t}/\vec{x}])(u[\vec{t}/\vec{x}]) \\[1.1mm]
(\lambda x.\,t)[\vec{t}/\vec{x}]  & := & \lambda x.\, t[\vec{t},x/\vec{x},x]\,.
\end{array}\]\deq[-1]
\end{mydefinition}
Note that, in contrast to ordinary substitution, simultaneous substitution can be defined directly on raw terms, that is, prior to equating them
modulo $\alpha\text{-equivalence}$. Moreover, we can show that:
\[ t[t_{1}/x_{1}\lldots t_{k}/x_{k}] = t[t_{1}/x_{1}]\cdots[t_{k}/x_{k}]\,. \]
We can now show the following \emph{Substitution Lemma}.
\begin{myproposition} For $t,t_1\lldots t_k$ as in the previous definition,
\[ \lsem t[t_{1}/x_{1}\lldots t_{k}/x_{k}] \rsem = \lsem t \rsem \circ \langle \lsem t_{1} \rsem \lldots \lsem t_{k} \rsem \rangle\,. \]
\end{myproposition}
\proof By induction on the structure of $t$.
\\
(1) If $t = x_{i}$:
\[ \lsem x_{i}[\vec{t}/\vec{x}]  \rsem = \lsem t_{i} \rsem = \pi_{i} \circ  \langle \lsem t_{1} \rsem \lldots \lsem t_{k}
\rsem \rangle = \lsem x_{i} \rsem \circ  \langle \lsem t_{1} \rsem \lldots \lsem t_{k} \rsem \rangle\,. \]
(2) If $t = uv$ then, abbreviating $\langle \lsem t_{1} \rsem \lldots \lsem t_{k} \rsem \rangle$ to $\langle \lsem \vec{t} \rsem \rangle$ we have:
\[ \begin{array}{rcl@{\quad}l}
\lsem uv[\vec{t}/\vec{x}] \rsem & = & \lsem  (u[\vec{t}/\vec{x}])(v[\vec{t}/\vec{x}])\rsem & \text{Defn of substitution} \\
& = & \App \circ \langle \lsem u[\vec{t}/\vec{x}] \rsem, \lsem v[\vec{t}/\vec{x}] \rsem \rangle & \text{Defn of semantic function} \\
& = &  \App \circ \langle \lsem u \rsem\circ\ang{\lsem\vec{t}\rsem}, \lsem v \rsem\circ\ang{\lsem\vec{t}\rsem}\rangle & \text{Induction hyp.} \\
& = &  \App \circ \langle \lsem u \rsem  , \lsem v \rsem \rangle \circ \langle \lsem \vec{t} \rsem \rangle & \text{Property of products} \\
& = & \lsem uv \rsem  \circ \langle \lsem \vec{t} \rsem \rangle & \text{Defn of semantic function}
\end{array}\]
(3) If $t = \lambda x.\,u$:
\[ \begin{array}{rcl@{\quad}l}
\lsem \lambda x.u [\vec{t}/\vec{x}] \rsem & = & \lsem \lambda x. (u[\vec{t},x/\vec{x},x]) \rsem &  \text{Defn. of substitution} \\
& = & \Curry{\lsem u[\vec{t},x/\vec{x},x]  \rsem} & \text{Defn. of semantic function} \\
& = &  \Curry{\lsem u \rsem \circ (\ang{\lsem\vec{t}\rsem}\times \mathsf{id})} &  \text{Induction hyp.} \\
& = & \Curry{\lsem u \rsem} \circ  \langle \lsem \vec{t} \rsem \rangle & \text{Prop.~\ref{p:Curry}} \\
& = & \lsem \lambda x.u \rsem \circ  \langle \lsem \vec{t} \rsem \rangle & \text{Defn. of semantic function}
\end{array}
\]
(4,5) The cases of projections and pairs are left as exercise. \qed

\begin{myexercise}
Complete the proof of the above proposition.
\end{myexercise}

\paragraph{Validating the conversion rules}
We can now show that the conversion rules of the $\lambda$-calculus are {preserved} by the translation, and hence the interpretation is
{sound}. Observe the correspondence between $\eta$-rules and uniqueness (couniversality) principles.
\begin{itemize}
    \item For $\beta$-conversion:\qd $\big[\text{\small$(\lambda x. \, t)u = t[u/x]\,,\ \pi_1\ang{t,u}=t\,,\ \pi_2\ang{t,u}=u$}\big]$
\begin{align*}
\lsem (\lambda x.\,t)u \rsem &=  \App \circ \langle \Curry{\lsem t \rsem}, \lsem u \rsem \rangle & \text{Defn. of semantics}& \\
& =  \App \circ (\Curry{\lsem t \rsem} \times \id) \circ \langle \id[\lsem \Gamma \rsem], \lsem u \rsem \rangle & \text{Property of $\times$}& \\
& = \lsem t \rsem \circ \langle \id[\lsem \Gamma \rsem], \lsem u \rsem \rangle & \text{Defn. of }\Lambda& \\
& = \lsem t[\vec{x},u/\vec{x},x] \rsem & \text{Substitution lemma}.&
\\[3mm]
\lsem \pi_1\ang{t,u} \rsem &= \pi_1\circ \lsem\ang{t,u}\rsem = \pi_1\circ \ang{\lsem t\rsem,\lsem u\rsem} = \lsem t \rsem\,.
\end{align*}
\item For
$\eta$-conversion:\qd $\big[\text{\small$t = \lambda x.\,tx \,,\ \ang{\pi_1t,\pi_2t}=t$}\big]$
\begin{align*}
    \lsem  \lambda x.\,tx \rsem &= \Curry{\App \circ (\lsem t \rsem \times \id)} = \lsem t \rsem & \text{Uniqueness equation }&(\impl)\\[2mm]
    \lsem \ang{\pi_1t,\pi_2t} \rsem &= \ang{\pi_1\circ\lsem t\rsem, \pi_2\circ\lsem t\rsem} = \lsem t \rsem & \text{Uniqueness equation }&(\times)
\end{align*}
\end{itemize}

\subsection{Completeness?}
It is the case that, in a general CCC $\CC$, there may be equalities which are \emph{not} reflected by the semantic translation, \ie
\[ \lsem t \rsem = \lsem u \rsem \text{\qd yet \qd} t\neq_\lambda u\,. \]
In the rest of this section,
we show how to \emph{construct} a CCC $\CC_\lambda$ in which equalities between arrows correspond precisely to $\lambda$-conversions between terms.
We call $\CC_\lambda$ a \emph{term model}, due to its dependence on the syntax.

\begin{mydefinition}\label{d:EqT}
We define a family of relations on variable-term pairs by setting \hbox{$(x,t)\sim_{T,U}(y,u)$} if $x:T\vdash t:U$ and $y:T\vdash u:U$ are derivable and
\[ t=_\lambda u[x/y]\,. \]
These are equivalence relations, so we set:
\[ [(x,t)]_{T,U} := \{\ (y,u)\ |\ (x,t)\sim_{T,U}(y,u)\ \}. \]
Similarly, $(\trt,t)\sim_{\trt,U}(\trt,u)$ if $\vdash t:U$ and $\vdash u:U$ are derivable and $t=_\lambda u$. Moreover,
\[ [(\trt,t)]_{\trt,U} := \{\ (\trt,u)\ |\ (\trt,t)\sim_{\trt,U}(\trt,u)\ \}\,. \]\deq[-1]
\end{mydefinition}
We denote $[(x,t)]_{T,U}$ simply by $[x,t]$\,, and $[(\trt,t)]_{\trt,U}$ simply by $[\trt,t]$ (these are \emph{not} to be confused with copairings!). We proceed with $\CC_\lambda$.

\begin{mydefinition}
The category $\CC_\lambda$ is defined as follows. We take as set of objects the set of $\lambda$-types augmented with a terminal object $\ena$:
\begin{align*}
Ob(\CC_\lambda)&:= \{\ena\}\cup\{\ \objk{T}\ |\ T\text{ a $\lambda$-type}\ \}
\myintertext{The homsets of $\CC_\lambda$ contain equivalence relations on typed terms (definition~\ref{d:EqT}), or terminal arrows $\tau$:}
        \CCl(\objk{T},\objk{U}) &:= \{\ [x,t]\ |\ x:T\vdash t:U\ \text{ is derivable}\ \} \\
        \CCl(\ena,\objk{U})     &:= \{\ [\trt,t]\ |\ \vdash t:U\ \text{ is derivable}\ \} \\
        \CCl(A,\ena)            &:= \{\ \tau_A \ \}
\myintertext{The identities are:}
\id[\objk{T}] &:= [x,x] \,,\quad \id[\ena] := \tau_\ena\,,
\myintertext{and arrow composition is defined by:}
        [x,t]\circ[y,u] &:= [y,t[u/x]] \\ 
        [x,t]\circ[\trt,u] &:= [\trt,t[u/x]] \\
        [\trt,t]\circ\tau_A &:=\begin{cases} [y,t] &\text{if }A=\objk{U}\\ [\trt,t] &\text{if }A=\ena \end{cases} \\
        \tau_B\circ h &:= \tau_A \qd[7]\;(h\in\CCl(A,B))
\end{align*}
\deq[-1]
\end{mydefinition}
Note that, for each variable $x'$, any arrow $[x,t]:\objk{T}\rarr\objk{U}$ can be written in the form $[x',t']$, since $t=(t[x'/x])[x/x']$ and therefore $[x,t]=[x',t[x'/x]]$.

\begin{myproposition} $\CCl$ is a category.
\end{myproposition}
\proof It is not difficult to see that $\id[]$'s are identities. For associativity, we show the most interesting case (and leave the rest as an
exercise):
\begin{align*}
  [x,t]\circ([y,u]\circ[z,v]) &= [x,t]\circ[z,u[v/y]] = [z,t[(u[v/y])/x]]\,, \\[3pt]
  ([x,t]\circ[y,u])\circ[z,v] &= [y,t[u/x]]\circ[z,v] = [z,t[u/x][v/y]]\,.
\end{align*}
By exercise~\ref{ex:SubSub}, the above are equal.
\qed

\begin{myproposition} $\CCl$ has finite products.
\end{myproposition}
\proof Clearly, $\ena$ is terminal with canonical arrows $\tau_A:A\rarr\ena$. For (binary) products, $\ena\times A=A\times\ena=A$. Otherwise, define
$\objk{T}\overset{\pi_1}{\llarr}\objk{T}\times\objk{U}\overset{\pi_2}{\lrarr}\objk{U}$ by:
\begin{align*}
    \objk{T}\times\objk{U} &:= \objk{T\times U} \\
    \pi_i &:= [x,\pi_ix] \qd i=1,2\,.
\end{align*}
Given \, $\objk{T}\overset{[x,t]}{\llarr}\objk{V}\overset{[x,u]}{\lrarr}\objk{U}$\,, take \; $\ang{[x,t],[x,u]}:\objk{V}\rarr\objk{T}\times\objk{U}:=
[x,\ang{t,u}]$\,. Then:
\begin{align*}
  \pi_1\circ \ang{[x,t],[x,u]} &= [y,\pi_1 y]\circ[x,\ang{t,u}] &&\text{Definitions}\\
        &=[x,\pi_1\ang{t,u}] &&\text{Defn of composition}\\
        &=[x,t] &&\beta\text{-conversion}
\end{align*}
Uniqueness left as exercise. The case of $\objk{T}\overset{[\trt,t]}{\llarr}\ena\overset{[\trt,u]}{\lrarr}\objk{U}$ is similar. \qed

\begin{myproposition} $\CCl$ has exponentials.
\end{myproposition}
\proof We have that $\ena\impl A=A$ and $A\impl\ena=\ena$, with obvious evaluation arrows. Otherwise,
\begin{align*}
  \objk{U}\impl\objk{T} &:= \objk{U\rarr T} \\
  \ev{\objk{U}}{\objk{T}} &:(\objk{U}\impl\objk{T})\times\objk{U}\lrar\objk{T}:=[x,(\pi_1x)(\pi_2x)]
\end{align*}
Given $[x,t]:\objk{V}\times\objk{U}\rarr\objk{T}$\,, take \, $\Lambda([x,t]):=[x_1,\lambda x_2.t[\ang{x_1,x_2}/x]]$\,. \\
Then,
\begin{align*}
\qd[3]  \evv\circ\Lambda([x,t])\times\id &= \evv\circ\ang{\Lambda([x,t])\circ\pi_1,\id\circ\pi_2} \\
    &= \evv\circ\ang{[x_1,\lambda x_2.t[\ang{x_1,x_2}/x]]\circ[y,\pi_1y],[y,\pi_2y]} \\
    &= \evv\circ\ang{[y,\lambda x_2.t[\ang{\pi_1y,x_2}/x]],[y,\pi_2y]} \\
    &= [z,(\pi_1z)(\pi_2z)]\circ[y,\rnode{A}{}\ang{\lambda x_2.t[\ang{\pi_1y,x_2}/x],\pi_2y}\rnode{B}{}] \\
    \ncbar[angle=-90,linewidth=.1mm,linecolor=darkgray,armB=1mm,armA=1mm,nodesep=1mm]{A}{B}\nbput[npos=2]{\color{darkgray}u}
    &= [y,(\pi_1u)(\pi_2u)]
    \overset{\beta}{=} [y,(\lambda x_2.t[\ang{\pi_1y,x_2}/x])(\pi_2y)] \\
    &\overset{\beta}{=} [y,t[\ang{\pi_1y,\pi_2y}/x]] \overset{\eta}{=} [y,t[y/x]] = [x,t]\,.
\end{align*}
Uniqueness left as exercise. The case of $[x,t]:\ena\times\objk{U}\rarr\objk{T}$ is similar. \qed

\begin{myexercise}
Complete the proof of the previous propositions.
\end{myexercise}
Hence, $\CC_\lambda$ is a CCC and a sound model of the simply-typed $\lambda$-calculus. Moreover, applying our translation from the $\lambda$-calculus to a CCC (definition~\ref{d:SemTr}) we can show that we have
\[ \lsem\Gamma\vdash t:T\rsem = [x,\, t[\pi_ix/x_i]_{i=1..n}] \]
where $\Gamma=\{x_1:T_1,...,x_n:T_n\}$, $x\notin\Gamma$ and $x:\prod_{i=1}^n T_i$\,. Then,
\[ t=_\lambda u \ \iff \ \lsem\Gamma\vdash t:T\rsem = \lsem\Gamma\vdash u:T\rsem\,.  \]
This means that our term model is \emph{complete}.

\subsection{Exercises}\label{ex:CH}
\begin{enumerate}\renewcommand{\theenumi}{\textbf{\arabic{enumi}}}
  \item Give Natural Deduction proofs of the following sequents.
    \begin{itemize}
    \item $\vdash (A \supset B) \supset ((B \supset C) \supset (A \supset C))$
    \item $\vdash (A \supset (A \supset B)) \supset (A \supset B)$
    \item $\vdash (C \supset A) \supset ((C \supset B) \supset (C \supset (A \wedge B)))$
    \item $\vdash (A \supset (B \supset C)) \supset ((A \supset B) \supset (A \supset C))$
    \end{itemize}
    In each case, give the corresponding $\lambda$-term and the corresponding arrow in a CCC $\CC$.
  \item For each of the following $\lambda$-terms, find a type for it. Try to find the `most general' type, built from `type variables' $\alpha$,
    $\beta$ etc. For example, the most general type for the identity $\lambda x. \, x$ is $\alpha \rightarrow \alpha$.
    In each case, give the derivation of the type for this term (where you may assume that types can be
    built up from type variables as well as base types).
    \begin{itemize}
      \item $\lambda f.\, \lambda x. \, fx$
      \item $ \lambda x. \, \lambda y. \, \lambda z. \, x (yz)$
      \item $ \lambda x . \, \lambda y. \, \lambda z. \, x z y$
      \item $ \lambda x.\, \lambda y.\, x y y$
      \item $ \lambda x.\, \lambda y.\, x$
      \item $ \lambda x. \, \lambda y.\, \lambda z. \, xz (yz)$
    \end{itemize}
    Reflect a little on the methods you used to do this exercise. Could they be made algorithmic?
\end{enumerate}

\section{Linearity}
In the system of Natural Deduction, implicit in our treatment of assumptions in sequents
\[ A_{1}, \ldots , A_{n} \vdash A \]
is that we can use them as many times as we want (including not at all). In this section we will explore the field that is opened once we apply
restrictions to this approach, and thus render our treatment of assumptions more \emph{linear} (or \emph{resource sensitive}).

\subsection{Gentzen Sequent Calculus}
In order to make the manipulation of assumptions more visible, we now represent the assumptions as a \emph{list} (possibly
with repetitions) rather than a set, and use explicit structural rules to control copying, deletion and interchange of assumptions.

\begin{mydefinition}\label{d:Struct} The \boldemph{structural rules} for Logic are given in the following table. \deq
\end{mydefinition}

\begin{center}\renewcommand{\arraystretch}{.5}
\fbox{$\begin{array}{@{\quad}c@{\quad}|@{\quad}c@{\quad}} &\\ %
\infer[\Ax]{A\vdash A}{} &  \infer[\Exch]{\Gamma,B,A,\Delta\vdash C}{\Gamma,A,B,\Delta\vdash C} \\&\\\hline&\\
\infer[\Con]{\Gamma,A\vdash B}{\Gamma,A,A\vdash B} & \infer[\Weak]{\Gamma,A\vdash B}{\Gamma\vdash B}\\ \end{array}$}
\captionof{table}{Structural Rules for Logic.}
\end{center}
If we think of using proof rules \emph{backwards} to reduce the task of proving a given sequent to various sub-tasks, then we see that the Contraction rule lets us  \emph{duplicate} premises, and the Weakening rule lets us \emph{discard} them, while the Exchange rule merely lets us re-order them.
The Identity axiom as given here is equivalent to the one with auxiliary premises given previously in the presence of Weakening.

The structural rules have clear categorical meanings in a category $\CC$ with products. Recalling the \emph{diagonal transformation}
$\Delta_A:=\ang{\id[A],\id[A]}$ and the \emph{symmetry transformation} $s_{A,B}:=\ang{\pi_2,\pi_1}$, the meanings are as follows.
\[\begin{array}{c@{\qd[2]}c}
\infer[\Exchs]{\Gamma , B, A, \Delta \vdash C}{\Gamma , A, B, \Delta  \vdash C} &
\infer{f \circ (\id[\Gamma]\times s_{A, B}\times\id[\Delta]):\Gamma\times B\times A\times\Delta\longrightarrow C}%
            {f:\Gamma\times A \times B \times \Delta \longrightarrow C} \\\\
\infer[\Cons]{\Gamma , A \vdash B}{\Gamma , A, A \vdash B}
& \infer{f\circ(\id[\Gamma]\times\Delta_{A}) : \Gamma \times A \longrightarrow B}{f : \Gamma \times A \times A \longrightarrow B} \\\\
\infer[\Weaks]{\Gamma , A \vdash B}{\Gamma \vdash B} & \infer{f\circ\pi_1:\Gamma \times A \longrightarrow B}{f : \Gamma \longrightarrow B}
\end{array}\]
In order to analyse Natural Deduction, Gentzen introduced \textbf{sequent calculi} based on \textbf{Left} and \textbf{Right} rules, instead of Introduction and Elimination rules. These kind of systems are more adequate for our discussion on linearity.

\begin{mydefinition}
We define the \boldemph{Gentzen sequent calculus} for $\wedge$,$\supset$ as the proof system obtained by the structural rules (def.~\ref{d:Struct}) and the rules in table~\ref{t:Gentz} for connectives.
\deq
\end{mydefinition}
\begin{center}
\renewcommand{\arraystretch}{.5}
\fbox{$\begin{array}{@{\;\;}c@{\;\;}|@{\;\;}c@{\;\;}|@{\;\;}c@{\;\;}} &&\\ \textbf{Conjunction} & \textbf{Implication} &
\textbf{Cut} \\&&\\\hline\hline&&\\
\infer[\AndR]{\Gamma , \Delta \vdash A \wedge B}{\Gamma \vdash A \qquad  \Delta \vdash B} &%
\infer[\ImpR]{\Gamma \vdash A \supset B}{\Gamma , A \vdash B} &%
\infer[\Cut]{\Gamma , \Delta \vdash B}{\Gamma \vdash A \qquad A, \Delta  \vdash B}
\\&&\\&&\\
\infer[\AndL]{\Gamma , A \wedge B \vdash C}{\Gamma , A, B \vdash C} &%
\infer[\ImpL]{\Gamma , A \supset B, \Delta \vdash C}{\Gamma \vdash A  \qquad B, \Delta \vdash C} &
\\&&\\
\end{array}$}\captionof{table}{Gentzen Sequent Calculus for $\wedge$,$\supset$.}\label{t:Gentz}
\end{center}
For example, the proof of $\supset$-transitivity is now given as follows.
\[ \infer[\ImpR]{A\supset B,B\supset C\vdash A\supset C}{\infer[\Exch]{A\supset B,B\supset C,A\vdash C}{
    \infer[\ImpL]{A\supset B,A,B\supset C\vdash C}{\infer[\Exch]{A\supset B,A\vdash B}{\infer[\ImpL]{A,A\supset B\vdash B}{
        \infer[\Ax]{A\vdash A}{}\quad\infer[\Ax]{B\vdash B}{}}}\quad\infer[\Ax]{C\vdash C}{}}}}
\]
%
%
\begin{myexercise}
Show that the Gentzen-rules are admissible in Natural Deduction. Moreover, show that the Natural Deduction rules are admissible in the Gentzen sequent
calculus.
\end{myexercise}
The $\Cut$ rule allows the use of \emph{lemmas} in
proofs. It also yields a \emph{dynamics of proofs} via \emph{Cut Elimination}, that is, a dynamics of proof transformations towards the goal of eliminating the uses of the Cut rule in a proof, \ie removing all lemmas and making the proof completely `explicit', meaning Cut-free.
Such transformations are always possible, as is shown in the following seminal result of Gentzen (\emph{Hauptsatz}).
\begin{fact}[Cut Elimination] The $\Cut$ rule is admissible in the Gentzen sequent calculus without $\Cut$.
\end{fact}

\subsection{Linear Logic}
In the presence of the structural rules, the Gentzen sequent calculus is entirely equivalent to the Natural Deduction system we studied earlier.
Nevertheless:
\begin{quote}
What happens if we \textbf{drop} the Contraction and Weakening rules (but keep the Exchange rule)?
\end{quote}
It turns out we can still make good sense of the resulting proofs, terms and categories, but now in the setting of a different, `resource-sensitive'
logic.

\begin{mydefinition}
\boldemph{Multiplicative Linear logic} is a variant of standard logic with \emph{linear} logical connectives. The multiplicative connectives for
conjunction and implication are $\otimes$ and $\linimpl$. {Proof sequents} are of the form $\Gamma \vdash A$, where $\Gamma$ is now a \emph{multiset}.
The proof rules for $\otimes$,$\linimpl$-Linear Logic, given in table~\ref{t:LL}, are the multiplicative versions of the Gentzen rules.
\deq\end{mydefinition}
\begin{center}
\renewcommand{\arraystretch}{.5}
\fbox{$\begin{array}{@{\;\;}c@{\;\;}|@{\;\;}c@{\;\;}|@{\;\;}c@{\;\;}} &&\\ \textbf{Conjunction} & \textbf{Implication} &
\textbf{Cut} \\&&\\\hline\hline&&\\
\infer[\TenR]{\Gamma , \Delta \vdash A \otimes B}{\Gamma \vdash A \qquad  \Delta \vdash B} &%
\infer[\LimpR]{\Gamma \vdash A \linimpl B}{\Gamma , A \vdash B} &%
\infer[\Cut]{\Gamma , \Delta \vdash B}{\Gamma \vdash A \qquad A, \Delta  \vdash B}
\\&&\\&&\\
\infer[\TenL]{\Gamma , A \otimes B \vdash C}{\Gamma , A, B \vdash C} &%
\infer[\LimpL]{\Gamma , A \linimpl B, \Delta \vdash C}{\Gamma \vdash A  \qquad B, \Delta \vdash C} &
\\&&\\
\end{array}$}\captionof{table}{Rules for $\otimes$,$\linimpl$-Linear Logic.}\label{t:LL}
\end{center}
Multiplicativity here means the use of \emph{disjoint} (\ie non-overlapping) contexts.
The use of multisets allows us to omit explicit use of the Exchange rule in our proof system.

Note that the given system satisfies Cut-elimination, and this leans heavily on the $\LimpL$ rule. We could have used instead the following rule,
\[ \infer[\LimpE]{\Gamma,\Delta\vdash B}{\Gamma\vdash A \linimpl B \qquad \Delta\vdash A} \]
which is more intuitive computationally, but then cut-elimination would fail. Note, though, that:
\[\label{e:qwe}
\LimpL\,,\,\Cut\,,\,\Ax \;\equiv\; \LimpE\,,\,\Cut\,,\,\Ax\,.
\]
This is shown as follows.
\begin{center}\small
$\hspace{-5mm}\infer[\Cut]{\Ga,A\linimpl B,\Delta\vdash C}%
    {\infer[\LimpE]{\Ga,A\linimpl B\vdash B}{\Ga\vdash A\qd\infer[\Ax]{A\linimpl B\vdash A\linimpl B}{}}\qd B,\Delta\vdash C\qd[4]\;\;}\hspace{-7mm}
\infer[\Cut]{\Ga,\Delta\vdash B}{\Ga\vdash A\linimpl B\qd \infer[\LimpL]{A\linimpl B,\Delta\vdash B}{\Delta\vdash A\qd\infer[\Ax]{B\vdash B}{}}}$
\end{center}
The resource-sensitive nature of Linear Logic is reflected in the following exercise.

\begin{myexercise} Can you construct proofs in Linear Logic of the following sequents? (Hint: Use the Cut Elimination property.)
\begin{itemize}
\item $A \vdash A \otimes A$
\item  $(A \otimes A) \linimpl B \vdash A \linimpl B$
\item $\vdash A \linimpl (B \linimpl A)$
\end{itemize}
\end{myexercise}
\newlines{1}
Related to linear logic is the \emph{linear $\lambda$-calculus}, which is a linear version of the simply-typed $\lambda$-calculus.

\begin{mydefinition}
The \boldemph{linear $\lambda$-calculus} is defined as follows.
\begin{align*}
&\text{Type} && \mathrm{TY}\ni\, T,U \; ::= \; b \; \mid \; T \linimpl U \; \mid \; T \otimes U
\\[1.5mm]
&\text{Term} && \mathrm{TE}\ni\,t,u\; ::= \; x \; \mid \; {t\,u} \;\mid\; {\lambda x.\, t} \;\mid\; t\otimes u \;\mid\; \olet{z}{x}{y}{t} \\[1.5mm]
&\text{Typing context} && \Gamma \;::=\; \varnothing\; \mid x : T,\Gamma
\end{align*}
Terms are typed by use of the typing rules of table~\ref{t:LinLambda}. Finally,
the rules for $\beta$-reduction are:
\[ \begin{array}{rcl}
(\lambda x. \, t)u & \longrightarrow_\beta & t[u/x] \\
\olet{t \otimes u}{x}{y}{v} & \longrightarrow_\beta & v[t/x,u/y]\,.
\end{array} \]
\deq[-1]
\end{mydefinition}
Note here that, again,  $x : T,\Gamma$ stands for $\{x:T\}\cup\Gamma$ with $x$ not appearing in $\Gamma$. Note also that Cut-free proofs \emph{always yield terms in normal form}.
\begin{center}
\renewcommand{\arraystretch}{.3}
{$\begin{array}{|@{\;}l@{\;}||@{\;\;}c@{\;\;}|@{\;\;}c@{\;\;}|} \hline&&\\&&\\&&\\
\parbox{15mm}{\bfseries Variable\\/\,Cut\\} & \infer{x:T \vdash  x:T}{} & \infer{\Gamma,\Delta\vdash u[t/x] : U}{\Gamma \vdash t : T  \qquad x : T, \Delta  \vdash u : U}
\\\hline&&\\&&\\&&\\
\parbox{15mm}{\bfseries Linear\\ Tensor\\} &\infer{\Gamma,\Delta\vdash t\otimes u:T\otimes U}{\Gamma\vdash t:T\qquad\Delta\vdash u:U} &%
\infer{\Gamma,z:T\otimes U\vdash\olet{z}{x}{y}{v}:V}{\Gamma,x:T,y:U\vdash v:V}
\\\hline&&\\&&\\&&\\
\parbox{15mm}{\bfseries Linear\\ Function\\} &\infer{\Gamma\vdash\lambda x.\,t:U\linimpl T}{\Gamma,x:U\vdash t:T} &%
\infer{\Gamma,f:T\linimpl U,\Delta\vdash u[ft/x]:V}{\Gamma\vdash t:T\qquad x:U,\Delta\vdash u:V}
\\\hline
\end{array}$}\captionof{table}{Linear $\lambda$-calculus for $\otimes,\linimpl$.}\label{t:LinLambda}
\end{center}

\noindent
Term formation is now highly constrained by the form of the typing judgements. In particular,
\[ x_1 : A_1 , \ldots , x_k : A_k \vdash t : A \]
now implies that each $x_i$ occurs \textbf{exactly once} (free) in $t$.

Moreover, note that, for function application, instead of the rule on the LHS
below, we could have used the more intuitive rule on the RHS.
\[ \infer{\Gamma,f:T\linimpl U,\Delta\vdash u[ft/x]:V}{\Gamma\vdash t:T\qquad x:U,\Delta\vdash u:V} \qquad\qquad
\infer{\Gamma,\Delta\vdash t\,u:B}{\Gamma\vdash t:A\linimpl B \qquad \Delta\vdash u:A}
\]
As we did in the logic, we can show that the typing systems with one or the other rule are equivalent.

\subsection{Linear Logic in Monoidal Categories}
We proceed to give a categorical counterpart to linearity by providing a categorical interpretation of linear logic.
Note that CCC's are no longer adequate for this task as they contain arrows
\[ \Delta_A:A\lrarr A\times A\,,\quad \pi_1:A\times B\lrarr A \]
which violate linearity. It turns out that the right setting is that of \emph{symmetric monoidal closed
categories}.

\begin{mydefinition}
A \boldemph{monoidal category} is a structure $(\CC , \otimes , I, a, l, r)$ where:
\begin{itemize}
\item $\CC$ is a category,
\item $\otimes : \CC \times \CC \rarr\CC$ is a functor (\emph{tensor}),
\item $I$ is a distinguished object of $\CC$ (\emph{unit}),
\item $a$, $l$, $r$ are natural isomorphisms (\emph{structural isos}) with components:
\[ a_{A, B, C} : A \otimes (B \otimes C) \iso (A \otimes B) \otimes C \]
\[ l_A : I \otimes A \iso A \qquad \quad r_A : A \otimes I \iso A \]
such that\, $l_{I} = r_{I} : I \otimes I \rarr I$\, and the following diagrams commute.
\[
\xymatrix@R=10mm@C=7.5mm{A\otimes(I\otimes B)\ar[r]^a\ar[dd]_{\id\otimes l} & (A\otimes I)\otimes B\ar[ddl]^{r\otimes\id} \\\\ A\otimes B}
\hspace{-10mm}
\xymatrix@R=8.2mm@C=-11mm{& (A \otimes B) \otimes (C \otimes D)\ar[dr]^a \\
    A\otimes(B\otimes(C\otimes D))\ar[ur]^a\ar[d]_{\id\otimes a}  && ((A\otimes B)\otimes C)\otimes D \\
    A\otimes((B\otimes C)\otimes D)\ar[rr]^a && (A\otimes(B\otimes C))\otimes D\ar[u]_{a\otimes\id}} \]
\end{itemize}\deq
\end{mydefinition}
The monoidal diagrams ensure \boldemph{coherence}, described by the slogan:
\begin{quote}\it
``\dots`all' diagrams involving $a,l$ and $r$ must commute."
\end{quote}
Examples:
\begin{itemize}
\item Both products and coproducts give rise to monoidal structures\,---\,which are the common denominator between them.
(But in addition, products have \emph{diagonals} and \emph{projections}, and coproducts have \emph{codiagonals} and \emph{injections}.)
\item $(\mathbb{N}, \leq , + , 0)$ is a monoidal category.
\item $\Rel$, the category of sets and relations, with cartesian  product (which is \emph{not} the categorical product).
\item $\Vect_k$ with the tensor product.
\end{itemize}
Let us examine the example of $\Rel$ in some detail. We take $\otimes$ to be the cartesian product, which is defined on relations $R:X\rarr X'$ and
$S:Y\rarr Y'$ as follows.
\[ \forall(x,y)\in X\times Y,(x',y')\in X'\times Y'.\; (x,y)R\otimes S(x',y')\iff xRx'\land ySy'\,. \]
It is not difficult to show that this is  indeed a functor. Note that, in the case that $R,S$ are \emph{functions}, $R\otimes S$ is the same as $R\times S$ in
\Set. Moreover, we take each $a_{A,B,C}$ to be the associativity function for products (in \Set), which is an iso in $\Set$ and hence also in \Rel.
Finally, we take $I$ to be the one-element set, and $l_A,r_A$ to be the projection functions: their relational converses are their inverses in \Rel.
The monoidal diagrams commute simply because they commute in $\Set$.
\begin{myexercise}
Verify that $(\mathbb{N}, \leq , + , 0)$ and $\Vect_k$ are monoidal categories.
\end{myexercise}

\paragraph{Tensors and products}
As we mentioned earlier, products are tensors with extra structure: natural diagonals and projections. This fact, which reflects \emph{no-cloning} and
\emph{no-deleting} of Linear Logic, is shown as follows.

\begin{myproposition}
Let $\CC$ be a monoidal category $(\CC,\otimes,I,a,l,r)$. $\otimes$ induces a product structure iff there exist natural diagonals and projections,
\ie natural transformations given by arrows
\[ d_A:A\lrar A\otimes A\,,\qd[2] p_{A,B}:A\times B\lrar A\,,\qd[2]q_{A,B}:A\times B\lrar B\,, \]
such that the following diagrams commute.
\[\xymatrix@=12mm{& A\ar[d]^{d_A}\ar[dl]_{\id[A]}\ar[dr]^{\id[A]}\\ A& A\otimes A\ar[l]^{p_{A,A}}\ar[r]_{q_{A,A}} & A}\qd[2]
  \xymatrix@=12mm{A\otimes B\ar[r]^-{d_{A,B}}\ar[dr]_{\id[A\otimes B]} & (A\otimes B)\otimes(A\otimes B)\ar[d]^{p_{A,B}\otimes q_{A,B}}\\ & A\otimes B}\]
\end{myproposition}
\proof The ``only if" direction is straightforward. For the converse, let $\CC$ be monoidal with natural projections and diagonals. Then, we take
product pairs to be pairs of the form
\[ A\overset{p_{A,B}}{\llarr}A\otimes B\overset{q_{A,B}}{\lrarr}B\,. \]
Moreover, for any pair of arrows $B\overset{f}{\llarr}A\overset{g}{\lrarr}C$\,, define
\[ \ang{f,g}:= A\overset{d_A}{\lrarr}A\otimes A\overset{f\otimes g}{\lrarr}B\otimes C\,. \]
Then the product diagram commutes. For example:
\[
\xymatrix@R=13mm@C=19.9mm{A\ar@{}[ddr]_{(1)}\ar[r]^{d_A}\ar[dd]_f& A\otimes A\ar[r]^{f\otimes g}\ar[d]^{f\otimes\id[A]}\ar@/_7mm/[dd]_{f\otimes f}
& B\otimes C\ar[dd]^{p_{B,C}}\\
& B\otimes A\ar[ur]_{\id[B]\otimes g}\ar[dr]^{p_{B,A}}_{(*)}\ar[d]^{\id[A]\otimes f}\ar@{}[r]|>>>>>>>>>{(*)} &{} \\
B\ar@/_7mm/[rr]_{\id[B]}^{(2)} \ar[r]^{d_B} & B\otimes B\ar[r]^{p_{B,B}} & B} %
\qd[3]\parbox[t]{.3\linewidth}{\newlines{3}\small$(*)$ naturality of $p$\\ $(1)$ naturality of $d$\\ $(2)$ hypothesis}
\]
For uniqueness, if $h:A\rarr B\otimes C$ then the following diagram commutes,
\[
\xymatrix@=13mm{A\ar[r]^h\ar[d]_{d_A}\ar@{}[dr]|{(1)} & B\otimes C\ar[d]_{d_{B\otimes C}}\ar[dr]^{\id[B\otimes C]}_{(2)} \\
A\otimes A\ar[r]_-{h\otimes h} & (B\otimes C)\otimes(B\otimes C)\ar[r]_-{p_{B,C}\otimes q_{B,C}} & B\otimes C}
\qd[3]\parbox[t]{.3\linewidth}{\newlines{3}\small$(1)$ naturality of $d$\\ $(2)$ hypothesis}
\]
so $h=\ang{\pi_1\circ h,\pi_2\circ h}$. \qed[1]

\paragraph{SMCC's} Linear Logic is interpreted in monoidal categories with two more pieces of structure: monoidal symmetry and closure.
The former allows the Exchange rule to be interpreted, while the latter realises
linear implication.

\begin{mydefinition}
A \boldemph{symmetric monoidal category} is a monoidal category $(\CC,\otimes,I,a,l,r)$ with an additional natural isomorphism (\emph{symmetry}),
\[ s_{A, B} : A \otimes B \iso B \otimes A \]
such that $s_{B,A}=s_{A,B}^{-1}$ and the following diagrams commute.
\[
\xymatrix@=9mm{ A\otimes I\ar[r]^{s}\ar[dr]_{r} & I\otimes A\ar[d]^{l} \\ & A}\qquad
\xymatrix@=9mm{ A\otimes(B\otimes C)\ar[d]_a\ar[r]^{\id\otimes s} & A\otimes(C\otimes B)\ar[r]^a & (A\otimes C)\otimes B\ar[d]^{s\otimes\id} \\
(A\otimes B)\otimes C\ar[r]_s & C\otimes(A\otimes B)\ar[r]_a & (C\otimes A)\otimes B}
\]\deq
\end{mydefinition}

\begin{mydefinition}
A \boldemph{symmetric monoidal closed category (SMCC)} is a symmetric monoidal category $(\CC , \otimes , I, a, l, r, s)$ such that, for each object
$A$, there is a couniversal arrow to the functor
\[ \uscore \otimes A : \CC \longrightarrow \CC\,. \]
That is, for all pairs $A,B$, there is an object $A \multimap B$ and a morphism
\[ \App[A, B] : (A \multimap B) \otimes A \longrightarrow B \]
such that, for every morphism $f : C \otimes A \rarr B$, there is a \emph{unique} morphism $\Lambda(f) : C \rarr(A \multimap B)$ such that
\[ \App[A, B] \circ (\Lambda (f) \otimes \id[A]) = f\,. \]\deq[-1]
\end{mydefinition}
Note that, although we use notation borrowed from CCC's ($\App,\Lambda$), these are \emph{different} structures! Examples of symmetric monoidal closed
categories are $\Rel$, $\Vect_k$, and (\emph{a fortiori}) cartesian closed categories.

\begin{myexercise}
Show that $\Rel$ is a symmetric monoidal closed category.
\end{myexercise}

\paragraph{Linear logic in SMCC's}
Just as cartesian closed categories correspond to $\wedge$,$\supset$-logic (and simply-typed $\lambda\text{-calculus}$), so
do symmetric monoidal closed categories correspond to \hbox{$\otimes$,$\linimpl$-logic} (and linear $\lambda$-calculus).

So let $\CC $ be a symmetric monoidal closed category. The interpretation of a linear sequent
\[ A_1 , \ldots , A_k \vdash A \]
will be a morphism
\[ f : A_1 \otimes \cdots \otimes A_k \longrightarrow A\,. \]
To be precise in our interpretation, we will again treat contexts as {lists} of formulas, and explicitly interpret the Exchange rule by:
\[ \infer{\Gamma , B, A, \Delta \vdash C}{\Gamma , A, B, \Delta  \vdash C}
\qquad \infer{f \circ (\id[\Gamma] \otimes s_{B, A} \otimes \id[\Delta])   : \Gamma
  \otimes B \otimes A \otimes \Delta \longrightarrow C}{f : \Gamma \otimes  A \otimes B \otimes \Delta \longrightarrow C}\]
The rest of the rules are translated as follows.
\begin{center}\renewcommand{\arraystretch}{0.6}
\fbox{$\begin{array}{@{\qd}c@{\qd}|@{\qd}c@{\qd}}
&\\
\infer{A \vdash A}{} & \infer{\id[A] : A \longrightarrow A}{}
\\&\\\hline&\\
\infer{\Gamma , \Delta \vdash B}{\Gamma \vdash A \qquad A, \Delta  \vdash B} & \infer{g\circ (f \otimes \id[\Delta]) : \Gamma \otimes \Delta
\longrightarrow B}{f:\Gamma\longrightarrow A\qquad g:A\otimes\Delta \longrightarrow B}
\\&\\\hline&\\
\infer{\Gamma , \Delta \vdash A \otimes B}{\Gamma \vdash A \qquad  \Delta \vdash  B} & \infer{f \otimes g : \Gamma \otimes \Delta \longrightarrow A
\otimes  B}{f : \Gamma \longrightarrow A \qquad g : \Delta \longrightarrow B}
\\&\\
\infer{\Gamma , A \otimes B \vdash C}{\Gamma , A, B \vdash C} & \infer{f \circ a_{\Gamma,A, B} : \Gamma \otimes (A \otimes B)  \longrightarrow C}{f :
(\Gamma \otimes A) \otimes B \longrightarrow C}
\\&\\\hline&\\
\infer{\Gamma \vdash A \linimpl B}{\Gamma , A \vdash B} & \infer{\Lambda (f) : \Gamma \longrightarrow (A \linimpl B)}{f : \Gamma  \otimes A
\longrightarrow B}
\\&\\
\infer{\Gamma,\Delta \vdash B}{\Gamma \vdash A \linimpl B \qquad \Delta \vdash A} & \infer{\App[A,B]\circ (f \otimes g) :\Gamma\otimes\Delta
\longrightarrow B}{f :\Gamma \longrightarrow (A \linimpl B)\qquad g : \Delta \longrightarrow A}
\\&\\
\end{array}$} \captionof{table}{Categorical Translation of $\otimes$,$\linimpl$-Linear Logic.}
\end{center}
Note that, because of coherence in monoidal categories, we will not be scholastic with associativity arrows $a$ in our translations and will usually omit them. For the same reason, consecutive applications of tensor will be written without specifying associativity, e.g.~$A_1\otimes\cdots\otimes A_n$.

\begin{myexercise} Let $\CC$ be a symmetric monoidal closed category. Give the interpretation of the $\linimpl$-left rule in $\CC$:
        \[ \infer[\LimpL]{\Gamma , A \linimpl B, \Delta \vdash C}{\Gamma \vdash A \qquad B, \Delta \vdash C} \]
\end{myexercise}

\begin{myexercise}
Is it possible to translate $\otimes$,$\linimpl$-logic into a CCC $\CC$? Is this in accordance with linearity of $\otimes$,$\linimpl$-logic?
\end{myexercise}

\subsection{Beyond the Multiplicatives}
Linear Logic has three `levels' of connectives, each describing a different aspect of standard logic:
\begin{itemize}
\item The \textbf{multiplicatives}: e.g. $\otimes$, $\linimpl$,
\item The \textbf{additives}: \emph{additive conjunction} $\With$ and \emph{disjunction} $\oplus$,
\item The \textbf{exponentials}, allowing controlled access to copying and  discarding.
\end{itemize}
We focus on additive conjunction and the exponential ``\,$\Bang$\,", which will allow us to recover the `expressive power' of standard
$\wedge$,$\supset$-logic.

\begin{mydefinition}\label{d:With} The logical connective for \boldemph{additive disjunction} is $\With$, and the related proof rules are the following.
\[\infer[\WithR]{\Gamma \vdash A \With B}{\Gamma \vdash A \qquad \Gamma \vdash B}
\qquad \infer[\WithL]{\Gamma , A \With B \vdash C}{\Gamma , A \vdash C} \qquad \infer[\WithL]{\Gamma , A \With B \vdash C}{\Gamma , B \vdash C}
\]\deq[-1]
\end{mydefinition}
So additive conjunction has proof rules that are identical to those of standard conjunction ($\wedge$). Note though that, since by linearity an
argument of type $A \With B$ can only be used once, each use of a left rule for $\With$ makes a once-and-for-all \emph{choice} of a projection. On the
other hand, $A\otimes B$ represents a conjunction where \emph{both} projections must be available.

Additive conjunction can be interpreted in any symmetric monoidal category with products, \ie a category $\CC$ with structure
$(\otimes,\times)$ where $\otimes$ is a symmetric monoidal tensor  and $\times$ is a product.
\[\infer{\ang{f,g}:\Gamma\lrar A\times  B}{f:\Gamma\lrar A \qquad g:\Gamma\lrar B}
\qd[2] \infer{f\circ(\id\otimes\pi_1):\Gamma\otimes(A\times B)\lrar C}{f:\Gamma\otimes A\lrar C} \]
Moreover, we can extend the linear $\lambda$-calculus with term constructors for additive conjunction as follows.
\[\begin{array}{cc}
\infer{\Gamma \vdash \langle t, u \rangle : A \With B}{\Gamma \vdash t : A \qquad \Gamma \vdash u : B}\qquad &%
\infer{\Gamma , z : A \With B \vdash \mathbf{let} \; z = \langle x,- \rangle \; \mathbf{in} \; t : C}{\Gamma , x : A \vdash t : C} \\[2mm]
& \infer{\Gamma , z : A \With B \vdash \mathbf{let} \; z = \langle - , y \rangle \; \mathbf{in} \; t : C}{\Gamma , B \vdash C}
\end{array}\]
The
$\beta$-reduction rules related to these constructs are:
\[ \begin{array}{lcl}
\mathbf{let} \; \langle t, u \rangle = \langle x, - \rangle \; \mathbf{in} \; v & \longrightarrow_\beta & v[t/x] \\
\mathbf{let} \; \langle t, u \rangle = \langle - , y \rangle \; \mathbf{in} \; v & \longrightarrow_\beta & v[u/y]\,.
\end{array} \]
Finally, we can gain back the lost structural rules, in \emph{disciplined} versions, by introducing an exponential \emph{bang} operator $!$ which
is a kind of \emph{modality} enabling formulas to participate in structural rules.
\begin{mydefinition}\label{d:Bang}
The logical connective for \boldemph{bang} is $\Bang$, and the related proof rules are the following.
\[
\infer[\someQWE{\Bang L}]{\Gamma,\Bang A\vdash B}{\Gamma,A\vdash B}\qd[2]\infer[\someQWE{\Bang R}]{\Bang\Gamma\vdash\Bang A}{\Bang\Gamma\vdash A}
\qd[2] \infer[\someQWE{Weak}]{\Gamma,\Bang A\vdash B}{\Gamma\vdash B} \qd[2] \infer[\someQWE{Contr}]{\Gamma,\Bang A\vdash B}{\Gamma,\Bang A,\Bang A
\vdash B}
\]
Note that $\Bang\{A_1,...,A_n\}:=\Bang A_1,...,\Bang A_n$\,.
\deq
\end{mydefinition}
We can now see the discipline imposed on structural rules: in order for the rules to be applied, the participating formulas need to be tagged
with a bang.

\paragraph{Interpreting standard Logic}
We are now in position to recover the standard logical connectives $\wedge$, $\supset$ within Linear Logic. If we interpret
\[ \begin{array}{lcl}
A \supset B & := & \Bang A \linimpl B \\
A \wedge B & := & A \With B
\end{array} \]
and each $\wedge$,$\supset$-sequent $\Gamma \vdash A$ as $\Bang\Gamma\vdash A$\,, then each proof rule of the Gentzen system for $\wedge$,$\supset$ is
admissible in the proof system of Linear Logic for $\otimes$,$\linimpl$,$\With$,$\Bang$\,.

Note in particular that the interpretation
\[ A \supset B \;\; := \;\; !A \linimpl B \]
\emph{decomposes} the fundamental notion of {implication} into finer notions\HY like `splitting the  atom of logic'!

\subsection{Exercises}
\begin{enumerate}\renewcommand{\theenumi}{\textbf{\arabic{enumi}}}
  \item Give proofs of the following sequents in Linear Logic.
    \begin{enumerate}
    \item $\vdash A \linimpl A$
    \item $ A \linimpl B , B \linimpl C \vdash A \linimpl C$
    \item $\vdash (A \linimpl B \linimpl C) \linimpl (B \linimpl A \linimpl C)$
    \item $A \otimes (B \otimes C) \vdash (A \otimes B) \otimes C$
    \item $A \otimes B \vdash B \otimes A$
    \end{enumerate}
    For each of the proofs constructed give:
    \begin{itemize}
    \item the corresponding linear $\lambda$-term,
    \item its interpretation in \Rel.
    \end{itemize}
  \item Consider a symmetric monoidal closed category $\CC$.
    \begin{enumerate}
        \item
        Suppose the sequents $\Gamma_1\vdash A$, $\Gamma_2\vdash B$ and $A,B,\Delta\vdash C$ are provable and let their interpretations (\ie the
        interpretations of their proofs) in $\CC$ be $f_1:\Gamma_1\rarr A$, $f_2:\Gamma_2\rarr B$ and $g:A\otimes B\otimes\Delta\rarr C$ respectively.
        Find then the interpretations $h_1,h_2$ of the following proofs.
        \begin{center}\small$\hspace{-6mm}
        \infer[\Cut]{\Gamma_1,\Gamma_2,\Delta\vdash C}{\infer[\TenR]{\Gamma_1,\Gamma_2\vdash A\otimes B}{\infer{\Gamma_1\vdash A}{\vdots}\qd%
        \infer{\Gamma_2\vdash B}{\vdots}}\qd
        \infer[\TenL]{A\otimes B,\Delta\vdash C}{\infer{A,B,\Delta\vdash C}{\vdots}}}
        \infer[\Cut]{\Gamma_1,\Gamma_2,\Delta\vdash C}{\qd[2]\infer{\Gamma_2\vdash B}{\vdots}
        \qd\infer[\Cut]{\Gamma_1,B,\Delta\vdash C}{\infer{\Gamma_1\vdash A}{\vdots}\qd\infer{A,B,\Delta\vdash C}{\vdots}}}$
        \end{center}
        and show that $h_1=h_2$.
        \item
        Suppose now $\CC$ has \emph{also} binary products, given by $\times$. Given that the sequents $\Gamma\vdash A$, $\Gamma\vdash B$ and
        $A,\Delta\vdash C$
        are provable, and that their interpretations in $\CC$ are $f_1:\Gamma\rarr A$, $f_2:\Gamma\rarr B$
        and $g:A\otimes\Delta\rarr C$ respectively, find the interpretations $h_1,h_2$ of the following proofs.
        \[ \infer[\Cut]{\Gamma,\Delta\vdash C}{\infer[\WithR]{\Gamma\vdash A\With B}{\infer{\Gamma\vdash A}{\vdots}\qd\infer{\Gamma\vdash B}{\vdots}}\qd
            \infer[\WithL]{A\With B,\Delta\vdash C}{\infer{A,\Delta\vdash C}{\vdots}}}
        \qd[3] \infer[\Cut]{\Gamma,\Delta\vdash C}{\infer{\Gamma\vdash A}{\vdots}\qd\infer{A,\Delta\vdash C}{\vdots}}
        \]
        and show that $h_1=h_2$.
    \end{enumerate}
  \item Show that the condition $l_I=r_I$ in the definition of monoidal categories is redundant.\\
        Moreover, show that the condition\, $\id[A]\otimes l_B=a_{A,I,B}\circ r_A\otimes\id[B]$\, in the definition of symmetric monoidal categories
        is redundant.
\end{enumerate}

\section{Monads and Comonads}
Recall that an adjunction is given by a triple $\ang{F,G,\theta}$, with $F:\CC\rarr\DD$ and $G:\DD\rarr\CC$ being functors, and $\theta$ a natural bijection between homsets. By composing the two functors we obtain endofunctors
\[ G\circ F:\CC\lrarr\CC\,,\quad F\circ G:\DD\lrarr\DD\,. \]
These can be seen as \emph{encapsulating the effect} of the adjunction inside their domain category. For example, if we consider the functors
\[ \MList:\Set\lrarr\Mon\,,\quad U:\Mon\lrarr\Set\,,\]
then $U\circ\MList$ encodes the free monoid construction inside $\Set$.

The study of such endofunctors on their own right gave rise to the notions of \emph{monad} and \emph{comonad}, which we examine in this section.

\subsection{Basics}
\begin{mydefinition}
A \boldemph{monad} over a category $\CC$ is a triple $(T,\eta,\mu)$ where $T$ is an endofunctor on $\CC$ and $\eta:\Id_{\CC}\rarr T$, $\mu:\TT\rarr T$
are natural transformations such that the following diagrams commute. (Note that $\TT:=T\circ T$, etc.).
  \[
  \xymatrix@=12mm{\TTT A\ar[r]^-{\mu_{TA}}\ar[d]_{T\mu_{A}} & \TT A\ar[d]^{\mu_A}\\ \TT A\ar[r]_-{\mu_A} & TA}\qd[3]
  \xymatrix@=12mm{TA\ar[r]^-{\eta_{TA}}\ar[dr]_{\id[TA]}\ar[d]_{T\eta_A} & \TT A\ar[d]^{\mu_A} \\
                  \TT A\ar[r]_{\mu_A} & TA}
  \]\deq[-1]
\end{mydefinition}
We call $\eta$ the \emph{unit} of the monad, and $\mu$ its \emph{multiplication}; the whole terminology comes from {monoids}. Let us now
proceed to some examples.
\begin{itemize}
\item Let $\CC$ be a category with coproducts and let $E$ be an object in $\CC$. We can define a monad $(T,\eta,\mu)$ of $E$-coproducts
    (computationally, $E$-exceptions) by taking $T:\CC\rarr\CC$ to be the functor $\uscore+E$\,, and $\eta,\mu$ as follows.
    \begin{align*}
      T &:= A\mapsto A+E\,,\, f\mapsto f+\id[E] \\
      \eta_A &:= A\lred{\incl[1]}A+E \\
      \mu_A &:=  (A+E)+E\lred{[\id[A+E],\,\incl[2]]}A+E
    \end{align*}
    As an injection, $\eta$ is a natural transformation. For $\mu$, we can use the properties of the coproduct. For $f:A\rarr B$,
   \begin{align*}
     Tf\circ\mu_A &= Tf\circ[\id[A+E],\incl[2]]=[Tf\circ\id[A+E],Tf\circ\incl[2]]=[Tf,Tf\circ\incl[2]] \\
     &= [Tf,(f+\id[E])\circ\incl[2]] = [Tf,\incl[2]] \\
     &= [\id[B+E]\circ Tf,\incl[2]\circ\id[E]]=[\id[B+E],\incl[2]]\circ(Tf+\id[E]) \\
     &= \mu_B\circ\TT f\,.
   \end{align*}
   The monadic diagrams follow in a similar manner. For example,
   \begin{align*}
     \mu_A\circ\mu_{TA} &= \mu_A\circ[\id[TA+E],\incl[2]] = [\mu_A\circ\id[TA+E],\mu_A\circ\incl[2]]=[\mu_A,\mu_A\circ\incl[2]] \\
       &= [\mu_A,[\id[A+E],\incl[2]]\circ\incl[2]] = [\mu_A,\incl[2]] \\
       &= [\id[A+E]\circ\mu_A, \incl[2]\circ\id[E]] = [\id[A+E], \incl[2]]\circ(\mu_A+\id[E]) \\
       &= \mu_A\circ T\mu_A\,.
   \end{align*}
    \item Now let $\CC$ be a cartesian closed category and let $\xi$ be some object in $\CC$. We can define a monad of $\xi$-side-effects by taking
    $T$ to be the functor $\xi\impl(\uscore\times\xi)$, and $\eta,\mu$ as follows.
    \begin{align*}
      T      &:= A\mapsto\xi\impl(A\times\xi)\,,\,f\mapsto\xi\impl(f\times\id[\xi]) \\
      \eta_A &:= \Lambda(\,A\times\xi\lred{\id[A\times\xi]}A\times\xi\,) \\
      \mu_A  &:= \Lambda(\, T(TA)\times\xi \lred{\ev\xi{TA\times\xi}} TA\times\xi\lred{\ev\xi{A\times\xi}}A\times\xi\,)
    \end{align*}
    Naturality of $\eta,\mu$ follows from naturality of $\Lambda$: for any $f:A\rarr A'$,
    \begin{align*}
      Tf\circ\eta_A &= (\xi\impl f\times\id[\xi])\circ\Lambda(\id[A\times\xi]) = \Lambda(f\times\id[\xi]\circ\id[A\times\xi]) \\
        &= \Lambda(\id[A'\times\xi]\circ f\times\id[\xi]) = \Lambda(\id[A'\times\xi])\circ f = \eta_{A'}\circ f\,, \\[1mm]
      \mu_{A'}\circ\TT f &= \Lambda(\ev\xi{A'\times\xi}\circ\ev\xi{TA'\times\xi})\circ\TT f
        = \Lambda(\ev\xi{A'\times\xi}\circ\ev\xi{TA'\times\xi}\circ\TT f\times\id[\xi]) \\
        &= \Lambda(\ev\xi{A'\times\xi}\circ Tf\times\id[\xi]\circ\ev\xi{TA\times\xi})
        = \Lambda(f\times\id[\xi]\circ\ev\xi{A\times\xi}\circ\ev\xi{TA\times\xi}) \\
        &= (\xi\impl f\times\id[\xi])\circ\Lambda(\ev\xi{A\times\xi}\circ\ev\xi{TA\times\xi}) = Tf\circ\mu_A\,.
    \end{align*}
    The monadic diagrams are shown in a similar manner.
\item Our third example employs the functor $U:\Mon\rarr\Set$. In particular, we take $T:=U\circ\MList$ and $\eta,\mu$ as follows.
    \begin{align*}
      T &:= X\mapsto \bigcup\nolimits_{n\in\omega}\{[x_1,\dots,x_n]\ |\ x_1,\dots,x_n\in X\}\,,\,\\
        &\qd[2]\! f\mapsto (\, [x_1,\dots,x_n]\mapsto [f(x_1),\dots,f(x_n)]\,)\,. \\[.75mm]
      \eta_X &:= x\mapsto [x] \\
      \mu_X  &:= [[x_{11},\dots,x_{1n_1}],\dots,[x_{k1},\dots,x_{kn_k}]]\mapsto[x_{11},\dots,x_{1n_1},\dots,x_{k1},\dots,x_{kn_k}]
    \end{align*}
    Naturality of $\eta,\mu$ is obvious\HY besides, $\eta$ is the unit of the corresponding adjunction. The monadic diagrams are also straightforward:
    they correspond to the following equalities of mappings (we use $\vec x$ for $x_1,\dots,x_n$).
    \[\xymatrix{[[[\vec x_{11}],...,[\vec x_{1n_1}]],...,[[\vec x_{k1}],...,[\vec x_{kn_k}]]]\ar@{|->}[r]^-{\mu}\ar@{|->}[d]_{T\mu}
        & [[\vec x_{11}],...,[\vec x_{1n_1}],...,[\vec x_{k1}],...,[\vec x_{kn_k}]]\ar@{|->}[d]^{\mu} \\
        [[\vec x_{11},...,\vec x_{1n_1}],\dots,[\vec x_{k1},...,\vec x_{kn_k}]] \ar@{|->}[r]_-{\mu}
        & [\vec x_{11},...,\vec x_{1n_1},\dots,\vec x_{k1},...,\vec x_{kn_k}]} \]
    \[\xymatrix{[x_1,\dots,x_n]\ar@{|->}[r]^-{\eta}\ar@{|->}[dr]^{\id}\ar@{|->}[d]_{T\eta} & [[x_1,\dots,x_n]]\ar@{|->}[d]^{\mu} \\
        [[x_1],\dots,[x_n]]\ar@{|->}[r]_-{\mu} & [x_1,\dots,x_n]} \]
\end{itemize}
\begin{myexercise}
Show that the $E$-coproduct monad and the $\xi$-side-effect monads are indeed monads.
\end{myexercise}
Our discussion on monads can be dualised, leading us to \emph{comonads}.
\begin{mydefinition}
A \boldemph{comonad} over a category $\CC$ is a triple $(Q,\eps,\delta)$ where $Q$ is an endofunctor on $\CC$ and $\eps:Q\rarr\Id_\CC$,
$\delta:Q\rarr\QQ$ are natural transformations such that the following diagrams commute.
  \[
  \xymatrix@=12mm{QA\ar[r]^-{\delta_{A}}\ar[d]_{\delta_{A}} & \QQ A\ar[d]^{\delta_{QA}}\\ \QQ A\ar[r]_-{Q\delta_A} & \QQQ A}\qd[3]
  \xymatrix@=12mm{QA\ar[r]^-{\delta_{A}}\ar[dr]_{\id[QA]}\ar[d]_{\delta_A} & \QQ A\ar[d]^{\eps_{QA}} \\
                  \QQ A\ar[r]_{Q\eps_A} & QA}
  \]\deq[-1]
\end{mydefinition}
$\eps$ is the \emph{counit} of the comonad, and $\delta$ its \emph{comultiplication}. Two of our examples from monads dualise to comonads.
\begin{itemize}
\item If $\CC$ has finite products then, for any object $S$, we can define the $S$-product comonad with functor $Q:=S\times\uscore$\,.
\item We can form a comonad on $\Mon$ with functor $Q:=\MList\circ U$ (and counit that of the corresponding adjunction).
\end{itemize}
\begin{myexercise}
Give an explicit description of the comonad on $\Mon$ with functor $Q:=\MList\circ U$ described above. Verify it is a comonad.
\end{myexercise}
\subsection{(Co)Monads of an Adjunction}
In the previous section, we saw that an adjunction between $\Mon$ and $\Set$ yielded a monad on $\Set$ (and a comonad on $\Mon$), with its unit being
the unit of the adjunction. We now show that this observation generalises to any adjunction. Recall that an adjunction is specified by:
\begin{itemize}
  \item a pair of functors\, $\xymatrix@1{{\CC\ }\ar@<.7mm>[r]^F & {\ \DD}\ar@<.5mm>[l]^G}$,
  \item for each $A\in Ob(\CC),B\in Ob(\DD)$, a bijection\, $\theta_{A,B} : \CC(A, GB)\cong\DD(FA, B)$\,
  natural in $A,B$.
\end{itemize}
For such an adjunction we build a monad on $\CC$: the functor of the monad is simply $T:=G\circ F$, and unit and multiplication are defined by setting
\begin{align*}
    \eta_A &: A\lrar GFA := \theta_{A,FA}^{-1}(\id[FA])\,,\\
    \mu_A  &: GFGFA\lrarr GFA := G(\theta_{GFA,FA}(\id[GFA]))\,.
\end{align*}
Observe that $\eta$ is the {unit} of the adjunction.
%
%
%
\begin{myproposition}
Let $(F,G,\eta)$ be an adjunction. Then the triple $(T,\eta,\mu)$ defined above is a monad on $\CC$.
\end{myproposition}
\proof
%
Recall that naturality of $\theta$ means concretely that, for any $f:A\rarr GB$, $g:A'\rarr A$ and $h:B\rarr B'$,
\[ \theta_{A',B'}(Gh\circ f\circ g) = h\circ\theta_{A,B}(f)\circ Fg\,. \]
$\eta$ is the unit of the adjunction and hence natural. We show naturality of $\mu$:
\begin{align*}
GFGFf\circ\mu_B &= G\theta_{GFB,FB}(\id[GFB])\circ GFGFf=G(\theta_{GFB,FB}(\id[GFB])\circ FGFf) \\
    &\overset{\text{nat.}\theta}{=} G\theta_{GFA,FB}(\id[GFB]\circ GFf) = G\theta_{GFA,FB}(GFf\circ\id[GFA]) \\
    &\overset{\text{nat.}\theta}{=} G(Ff\circ\theta_{GFA,FA}(\id[GFA])) = GFf\circ\mu_A\,.
\end{align*}
The monoidal condition for $\mu$ also follows from naturality of $\theta$:
\begin{align*}
\mu_A\circ\mu_{GFA} &= G(\theta(\id[GFA])\circ\theta(\id[GFGFA])) \overset{\text{nat}}{=} G\theta(G\theta(\id[GFA])\circ\id[GFGFA]) \\
    &= G\theta(\id[GFA]\circ G\theta(\id[GFA])) \overset{\text{nat}}{=} G(\theta(\id[GFA])\circ FG\theta(\id[GFA])) \\
    &= \mu_A\circ GF\mu_A\,.
\end{align*}
Finally, for the $\eta$-$\mu$ conditions we also use the universality diagram for $\eta$ and the uniqueness property (in equational form).
\begin{align*}
\mu_A\circ\eta_{GFA} &= G\theta_{GFA,FA}(\id[GFA])\circ\eta_{GFA} = \id[GFA]\,, \\
\mu_A\circ GF\eta_{GFA} &= G\theta_{GFA,FA}(\id[GFA])\circ GF\eta_{GFA} = G(\theta_{GFA,FA}(\id[GFA])\circ F\eta_{GFA}) \\
    &\overset{\text{nat}}{=} G\theta_{}(\id[GFA]\circ\eta_{GFA}) = G\theta_{}(G\id[FA]\circ\eta_{GFA}) = G\id[FA]  = GFA\,.
\end{align*}
\qed[-3]
Hence, every adjunction gives rise to a monad. It turns out that the converse is also true: \emph{every} monad is described by means of an adjunction
in this way. In particular, there are two canonical constructions of adjunctions from a given monad: the \emph{Kleisli construction}, and the
\emph{Eilenberg-Moore construction}. These are in a sense minimal and maximal solutions to describing a monad via an adjunction. We describe the former
one in the next section. 

Finally, note that\HY because of the symmetric definition of adjunctions\HY the whole discussion can be dualised to comonads. That is, every adjunction
gives rise to a comonad with counit that of the adjunction, and also every comonad can be derived from an adjunction in this manner.

\subsection{The Kleisli Construction}
The Kleisli construction starts from a monad $(T,\eta,\mu)$ on a category $\CC$ and builds a category $\CC_T$ of $T$-computations, as follows.
\begin{mydefinition}
Let $(T,\eta,\mu)$ be a monad on a category $\CC$. Construct the \boldemph{Kleisli category} $\CC_T$ by taking the same objects as $\CC$, and by
including an arrow $\kl{f}:A\rarr B$ in $\CC_T$ for each $f:A\rarr TB$ in $\CC$. That is,
\begin{align*}
  Ob(\CC_T) &:= Ob(\CC)\,, \\
  \CC_T(A,B) &:= \{ \kl{f}\ |\ f\in\CC(A,TB)\}\,.
\end{align*}
The identity arrow for $A$ in $\CC_T$ is $\eta_{A.T}$\,, while the composite of $\kl{f}:A\rarr B$ and $\kl{g}:B\rarr C$ is $\kl{h}$\,, where:
\[ h := A\lred{f}TB\lred{Tg}\TT C\lred{\mu_C}TC\,. \]
\deq[-1]
\end{mydefinition}
The conditions for $\CC_T$ being a category follow from the monadic conditions. For composition with identity, for any $f:A\rarr TB$,
\begin{align*}
    \kl{f}\circ\eta_{A.T} &= \kl{(\mu_B\circ Tf\circ\eta_A)} = \kl{(\mu_B\circ\eta_B\circ f)}=\kl{f}\,,  \\
    \eta_{B.T}\circ\kl{f} &= \kl{(\mu_B\circ T\eta_B\circ f)} = \kl{f}\,.
\end{align*}
For associativity of composition, for any $f:A\rarr TB$, $g:B\rarr TC$ and $h:C\rarr TD$,
\begin{align*}
  (\kl{h}\circ\kl{g})\circ\kl{f} &= \kl{(\mu_D\circ Th\circ g)}\circ\kl{f} = \kl{(\mu_D\circ T(\mu_D\circ Th\circ g)\circ f)} \\
    &= \kl{(\mu_D\circ T\mu_D\circ \TT h\circ Tg\circ f)} = \kl{(\mu_D\circ \mu_{TD}\circ \TT h\circ Tg\circ f)}\\
    &= \kl{(\mu_D\circ Th\circ\mu_C\circ Tg\circ f)} = \kl{h}\circ(\kl{g}\circ\kl{f})\,.
\end{align*}
Let us now proceed to build the adjunction between $\CC$ and $\CC_T$ that will eventually give us back the monad $T$. Construct the functors
$F:\CC\rarr\CC_T$ and $G:\CC_T\rarr\CC$ as follows.
\begin{align*}
  F&:= A\mapsto A\,,\, (f:A\rarr B) \mapsto (\kl{(\eta_B\circ f)}:A\rarr B)\,, \\[1mm]
  G&:= A\mapsto TA\,,\, (\kl{f}:A\rarr B) \mapsto (\mu_B\circ Tf:TA\rarr TB)\,.
\end{align*}
Functoriality of $F,G$ follows from the monad laws and the definition of $\CC_T$.
Moreover, for each $A,B\in Ob(\CC)$, construct the following bijection of arrows.
\[ \theta_{A,B}:\CC(A,TB)\iso\CC_T(A,B):= f\mapsto\kl{f} \]
To establish that $(F,G,\theta)$ is an adjunction we need only show that $\theta$ is natural in $A,B$. So take $f:A\rarr TB$, $g:A'\rarr A$ and
$\kl{h}:B\rarr B'$. We then have:
\begin{align*}
  \theta_{A',B'}(G(\kl{h})\circ f\circ g) &= \theta_{A',B'}(\mu_{B'}\circ Th\circ f\circ g) = \kl{(\mu_{B'}\circ Th\circ f\circ g)} \\
    &= \kl{h}\circ\kl{(f\circ g)} = \kl{h}\circ\kl{(\mu_B\circ Tf\circ\eta_A\circ g)}\\
    &= \kl{h}\circ\kl{f}\circ\kl{(\eta_A\circ g)} = \kl{h}\circ\theta_{A,B}(f)\circ Fg\,.
\end{align*}
The final step in this section is to verify that the monad $(T',\eta',\mu')$ arising from this adjunction is the one we started from. The
construction of $T'$ follows the recipe given in the previous section, that is:
\begin{itemize}
  \item $T':\CC\rarr\CC:=G\circ F$\,. Thus, $T'$ maps each object $A$ to $TA$, and each arrow $f:A\rarr B$ to $\mu_B\circ T\eta_A\circ Tf=Tf$.
  \item $\eta'_A:A\rarr TA:=\theta_{A,FA}^{-1}(\id[FA]^{(\CC_T)})=\theta^{-1}(\eta_{A.T})=\eta_A$\,.
  \item $\mu'_A:\TT A\rarr TA:=G\theta_{GFA,FA}(\id[GFA]^{(\CC)})=G\theta(\id[TA])=\mu_A\circ T\id[TA]=\mu_A$\,.
\end{itemize}
Thus, we have indeed obtained the initial $(T,\eta,\mu)$.

\paragraph{The Kleisli construction on a comonad}
Dually to the Kleisli category of a monad we can construct the Kleisli category of a comonad\footnote{In some texts, this is called a \emph{coKleisli}
category.}\HY and reobtain the comonad through an adjunction between the Kleisli category and the original one. Specifically, given a category $\CC$
and a comonad $(Q,\eps,\delta)$ on $\CC$, we define the category $\CC_Q$ as follows.
\begin{align*}
Ob(\CC_Q) &:= Ob(\CC) \\
\CC_Q(A,B)&:= \{ \kl[Q]{f}\ |\ f\in\CC(QA,B) \} \\
\id[A]^{(\CC_Q)} &:= \eps_{A.Q} \\
\kl[Q]{g}\circ\kl[Q]{f} &:= \kl[Q]{(g\circ Qf\circ\delta_A)}
\end{align*}
The Kleisli category of a comonad will be of use in the next sections, where comonads will be considered for modelling \emph{bang} of Linear Logic. We
end this section by showing a result that will be of use then.
\begin{myproposition}\label{p:Qprods}
Let $\CC$ be a category and $(Q,\eps,\delta)$ be a comonad on $\CC$. If $\CC$ has binary products then so does $\CC_Q$.
\end{myproposition}
\proof%
Let $A,B$ be objects in $\CC,\CC_Q$. We claim that their product in $\CC_Q$ is given by $(A\times B,p_1,p_2)$, where
\[ p_1:=\kl[Q]{\big(Q(A\times B)\lred{\eps}A\times B\lred{\pi_1}A\big)} \]
and similarly for $p_2$. Now, for each $\kl[Q]{f}:C\rarr A$ and $\kl[Q]{g}:C\rarr B$, setting $\ang{\kl[Q]{f},\kl[Q]{g}}:=\kl[Q]{\ang{f,g}}$ we have:
\[
 p_1\circ\ang{\kl[Q]{f},\kl[Q]{g}} = \kl[Q]{(\pi_1\circ\eps\circ Q\ang{f,g}\circ\delta)} = \kl[Q]{(\pi_1\circ
 \ang{f,g}\circ\eps\circ\delta)}=\kl[Q]{f}\,,
\]
and similarly $p_2\circ\ang{\kl[Q]{f},\kl[Q]{g}} =\kl[Q]{g}$\,. Finally, for any $\kl[Q]{h}:C\rarr A\times B$,
\begin{align*}
 \ang{p_1\circ\kl[Q]{h},p_2\circ\kl[Q]{h}} &= \kl[Q]{\ang{\pi_1\circ\eps\circ Qh\circ\delta,\pi_2\circ\eps\circ Qh\circ\delta}}
    = \kl[Q]{\ang{\pi_1\circ h,\pi_2\circ h}}\\ &= \kl[Q]{h}\,.
\end{align*}
\qed[-1]

\begin{myexercise}\label{e:Qterms}
Show that the Kleisli category $\CC_Q$ of a comonad $(Q,\eps,\delta)$ has a terminal object when $\CC$ does.
\end{myexercise}
\subsection{Modelling of Linear Exponentials}
In this section we employ comonads in order to model the exponential \emph{bang} operator, $\Bang$\,, of Linear Logic. Let us start by modelling
a \emph{weak bang} operator, $\WBang$\,, which involves solely the following proof rules.
\[
\infer[\someQWE{\WBang L}]{\Gamma,\WBang A\vdash B}{\Gamma,A\vdash B}\qd[3]\infer[\someQWE{\WBang R}]{\WBang B\vdash\WBang A}{\WBang B\vdash A}
\]
Observe that, compared to $\Bang$\,, $\WBang$ is weak in its Right rule, and it also misses Contraction and Weakening.

Let us now assume as given a symmetric monoidal closed category $\CC$ along with a comonad $(Q,\eps,\delta)$ on $\CC$. As seen previously, $\CC$ is a model of
($\otimes$$\linimpl$)-Linear Logic. Moreover, $(\CC,Q)$ yields a model of ($\otimes$$\linimpl$$\WBang$)-Linear Logic by modelling each formula $\WBang
A$ by $QA$ (\ie $Q$ applied to the translation of $A$). The rules for weak bang are then interpreted as follows.
\[
\infer{f\circ\id[\Gamma]\otimes\eps_A:\Gamma\otimes QA\lrar B}{f:\Gamma\otimes A\lrar B} \qd[3] %
\infer{Qf\circ\delta_B:QB\lrar QA}{f:QB\lrar A}
\]
We know that arrow-equalities in $\CC$ correspond to proof-transformations in the proof system.
Thus, the comonadic law\, $\eps_{QA}\circ\delta_A = \id[QA] = Q\eps_A\circ\delta_A$\, corresponds to the following transformations.
\[
\xymatrix{\infer[\Cut]{\WBang A\vdash\WBang A}{\infer[\someQWE{\WBang R}]{\WBang A\vdash\WBang\WBang A}{\infer[\Ax]{\WBang A\vdash\WBang A}{}}\qd[2] %
\infer[\someQWE{\WBang L}]{\WBang\WBang A\vdash\WBang A}{\infer[\Ax]{\WBang A\vdash\WBang A}{}}}\qd\;\ar@{=>}@<-3.5mm>[r] & %
\qd\;\deduce{\infer[\Ax]{\WBang A\vdash\WBang A}{}}{\vspace{8mm}}\qd\; & %
\qd\;\infer[\someQWE{\WBang R}]{\WBang A\vdash\WBang A}{\infer[\someQWE{\WBang L}]{\WBang A\vdash A}{\infer[\Ax]{A\vdash A}{}}}\ar@{=>}@<3.5mm>[l]}
\]
\begin{myexercise}
Find a proof-transformation corresponding to the comonadic law $\delta_{QA}\circ\delta_A=Q\delta_A\circ\delta_A$\,.
\end{myexercise}
In order to extend our translation to the general $\someQWE{\Bang R}$ rule, we need arrows in $\CC$ of the form
\[ \QQ A_1\otimes\cdots\otimes\QQ A_n\lrar Q(QA_1\otimes\cdots\otimes QA_n)\,. \]
Hence, we need to impose (a coherent) distributivity of the tensor\HY either binary ($\otimes$) or nullary ($I$)\HY over the comonad $Q$. This can be
formalised by stipulating that $Q$ be a \emph{symmetric monoidal} endofunctor.

\begin{mydefinition}
Let $(\CC,\otimes,I,a,l,r,s)$ and $(\CC',\otimes',I',a',l',r',s')$ be symmetric monoidal categories. A functor $F:\CC\rarr\CC'$ is called
\boldemph{symmetric monoidal} if there exist:
\begin{itemize}
  \item a morphism $m_0:I'\rarr F(I)$\,,
  \item a natural transformation $m_2:F(\uscore)\otimes'F(\uscore)\rarr F(\uscore\otimes\uscore)$\,,
\end{itemize}
such that the following diagrams commute.
\[\xymatrix@R=8mm@C=13mm{FA\otimes'(FB\otimes'FC)\ar[d]_{a'}\ar[r]^-{\id\otimes'm_2} & FA\otimes'F(B\otimes C)\ar[r]^{m_2}
    & F(A\otimes(B\otimes C))\ar[d]^{Fa} \\
    (FA\otimes'FB)\otimes'FC\ar[r]_-{m_2\otimes'\id} & F(A\otimes B)\otimes'FC\ar[r]_{m_2} & F((A\otimes B)\otimes C)} \]
\[\xymatrix@R=8mm@C=13mm{FA\otimes'I'\ar[r]^-{\id\otimes'm_0}\ar[d]_{r'} & FA\otimes'FI\ar[d]^{m_2} \\
    FA & F(A\otimes I)\ar[l]^-{Fr}}\qd[2]
\xymatrix@R=7.5mm@C=13mm{FA\otimes'FB\ar[d]_{s'}\ar[r]^{m_2} & F(A\otimes B)\ar[d]^{Fs}\\
    FB\otimes'FA\ar[r]_{m_2} & F(B\otimes A)}     \]
We may write such an $F$ as $(F,m)$. Moreover, if $(F,m),(G,n):\CC\rarr\CC'$ are (symmetric) monoidal functors then a natural transformation
$\phi:F\rarr G$ is called \boldemph{monoidal} whenever the following diagrams commute.
\[\xymatrix@R=8mm@C=13mm{I'\ar[r]^{m_0}\ar[dr]_{n_0} & FI\ar[d]^{\phi}\\ & GI}\qd[3]
\xymatrix@R=7.4mm@C=13mm{FA\otimes'FB\ar[r]^{m_2}\ar[d]_{\phi\otimes'\phi} & F(A\otimes B)\ar[d]^{\phi}\\ GA\otimes'GB\ar[r]_{n_2} & G(A\otimes B)} \]
\deq[-1]
\end{mydefinition}
For example, the identity functor is symmetric monoidal. Moreover, if $F$ and $G$ are symmetric monoidal functors then so is $G\circ F$. Other examples
are the following.
\begin{itemize}
\item The constant endofunctor $K_I$, which maps each object to $I$ and each arrow to $\id[I]$, is symmetric monoidal with structure maps:
    \[ m_0:I\lrar I:=\id[I]\,,\qd m_2:I\otimes I\lrar I:=r_I\,.\]
\item The endofunctor $\otimes\circ\ang{\Id_\CC,\Id_\CC}$, which maps each object $A$ to $A\otimes A$ and each arrow $f$ to $f\otimes f$, is symmetric
    monoidal with:
    \[ m_0:I\lrar I\otimes I:=r_I^{-1}\,,\qd
       m_2:=(A\otimes \rnode{A}{A})\otimes(\rnode{B}{B}\otimes B)\lrar(A\otimes \rnode{BB}{B})\otimes(\rnode{AA}{A}\otimes B)\,, \]
       \cc[7pt]{A}{AA}\cc[4pt]{B}{BB}
    the latter given by use of structural transformations.
\end{itemize}
\begin{myexercise}
Verify that if $F:\CC\rarr\DD$, $G:\DD\rarr\EE$ are symmetric monoidal functors then so is $G\circ F$.
\end{myexercise}
\begin{mydefinition}
A comonad $(Q,\eps,\delta)$ on a SMCC $\CC$ is called a \boldemph{monoidal comonad} if $Q$ is a symmetric monoidal functor, say $(Q,m)$, and
$\eps,\delta$ are monoidal natural transformations. We write $Q$ as $(Q,\eps,\delta,m)$. \deq
\end{mydefinition}
Now let us assume $\CC$ is a SMCC and $(Q,\eps,\delta,m)$ is a monoidal comonad on $\CC$. The coherence of $m_2$ with $a$, expressed by the first
diagram of symmetric monoidal functors, allows us to generalise $m_0$ and $m_2$ to arbitrary arities and assume arrows:
\[ m_n: QA_1\otimes\cdots\otimes QA_n\lrar Q(A_1\otimes\cdots\otimes A_n)\,. \]
We can give the interpretation of the Right rule for bang as follows.
\[
\infer{Qf\circ m_n\circ(\delta_{B_1}\otimes\cdots\otimes\delta_{B_n}):QB_1\otimes\cdots\otimes QB_n\lrar QA}{f:QB_1\otimes\cdots\otimes QB_n\lrar A}
\]

\paragraph{Contraction and Weakening}
Our discussion on the categorical modelling of linear exponentials has only touched the issues of Right and Left rules. However, we also need adequate structure for translating Contraction and Weakening.
\[ \infer[\someQWE{Contr}]{\Gamma,\Bang A\vdash B}{\Gamma,\Bang A,\Bang A\vdash B}\qd[2]\infer[\someQWE{Weak}]{\Gamma,\Bang A\vdash B}{\Gamma\vdash B}
\]
For these rules we can use appropriate (monoidal) natural transformations.
For Contraction, we stipulate a transformation with components $d_A:QA\rarr QA\otimes QA$\,, \ie
\[ d: Q\lrar \otimes\circ\ang{Q,Q}\,. \]
For Weakening, a transformation with components $e_A:QA\rarr I$, \ie
\[ e: Q\lrar K_I. \]
Although the above allow the categorical interpretation of the proof-rules, they do not necessarily preserve the intended proof-transformations. For
that, we need to impose some further coherence conditions, which are epitomised in the following notion.
\begin{mydefinition}
Let $\CC$ be a SMCC. A monoidal comonad $(Q,\eps,\delta,m)$ on $\CC$ is called a \boldemph{linear exponential comonad} if there exist monoidal natural
transformations
\[ d: Q\lrar \otimes\circ\ang{Q,Q}\,,\qd e: Q\lrar K_I, \]
such that:
\begin{enumerate}\renewcommand{\theenumi}{(\alph{enumi})}\renewcommand\labelenumi{\theenumi}
  \item for each object $A$, the triple $(QA,d_A,e_A)$ is a \emph{commutative comonoid in $\CC$}, \ie the following diagrams commute,
  \[ \xymatrix@C=13mm@R=8.75mm{QA\ar[r]^-{d_A}\ar[dr]^{d_A} & QA\otimes QA\ar[d]^{s_{QA,QA}}\\ I\otimes QA\ar[u]^{l_{QA}}
  & QA\otimes QA\ar[l]^{e_A\otimes\id[QA]}}\qd
  \xymatrix@C=16mm{QA\ar[r]^-{d_A}\ar[d]_{d_A} & QA\otimes QA\ar[d]^{d_A\otimes\id[QA]}\\
  QA\otimes QA\ar[r]_-{a\circ(\id[QA]\otimes d_A)} & (QA\otimes QA)\otimes QA} \]
  \item for each object $A$, the following diagrams commute.
  \[\begin{array}{@{\!\!\!\!\!\!}c l}
  \xymatrix@C=12mm{QA\ar[r]^{\delta_A}\ar[d]_{e_A} & \QQ A\ar[d]^{Qe_A}\\ I\ar[r]_{m_0} & QI} &
  \xymatrix@C=11mm{QA\ar[r]^-{\delta_A}\ar[d]_{d_A} & \QQ A\ar[dr]^{Qd_A}\\
    QA\otimes QA\ar[r]_-{\delta_A\otimes\delta_A} & \QQ A\otimes\QQ A\ar[r]_{m_2} & Q(QA\otimes QA)} \\
  \xymatrix@C=12mm{QA\ar[r]^{\delta_A}\ar[dr]_{e_A} & \QQ A\ar[d]^{e_{QA}}\\ & I} &
  \xymatrix@C=11mm{QA\ar[r]^-{\delta_A}\ar[d]_{d_A} & \QQ A\ar[d]^{d_{QA}}\\
    QA\otimes QA\ar[r]_-{\delta_A\otimes\delta_A} & \QQ A\otimes\QQ A}
  \end{array}\]
\end{enumerate}
We write $Q$ as $(Q,\eps,\delta,m,d,e)$. \deq
\end{mydefinition}

\begin{myexercise}
Express what it means concretely for $d,e$ to be monoidal natural transformations.
\end{myexercise}
\begin{myexercise}
Give the categorical interpretation of Contraction and Weakening in a SMCC $\CC$ with a linear exponential comonad.
\end{myexercise}

\subsubsection{Including Products}

We now consider the fragment of Linear Logic which includes all four linear connectives we have seen thus far, \ie $\otimes\linimpl\Bang\With$, and
their respective proof rules (see definitions~\ref{d:With},\,\ref{d:Bang}). The categorical modelling of ($\otimes\linimpl\Bang\With$)-Linear Logic
requires:
\begin{itemize}
 \item a symmetric monoidal closed category $\CC$,
 \item a linear exponential comonad $(Q,\eps,\delta,m,d,e)$ on $\CC$,
 \item finite products in $\CC$.
\end{itemize}
The above structure is adequate for modelling the proof rules as we have seen previously. Moreover, it provides rich structure for the Kleisli category $\CC_Q$. The next result and its proof demonstrate categorically the `interpretation' of ordinary logic within Linear Logic given by:
\[ A\impl B \equiv \Bang A\linimpl B\,. \]

\begin{myproposition}\label{p:CCQ}
Let $\CC$ be a SMCC with finite products and let $(Q,\eps,\delta,m,d,e)$ be a linear exponential comonad on $\CC$. Then:
\begin{enumerate}\renewcommand{\theenumi}{\rm(\alph{enumi})}\renewcommand\labelenumi{\theenumi}
  \item The Kleisli category $\CC_Q$ has finite products.
  \item There exists an isomorphism $i:Q\ena\rarr I$ and a natural isomorphism $j:Q(\uscore\times\uscore)\rarr Q(\uscore)\otimes Q(\uscore)$.
  \item $\CC_Q$ is cartesian closed, with the exponential of objects $B,C$ being $QB\linimpl C$.
\end{enumerate}
\end{myproposition}
\proof Part (a) has been shown previously (proposition~\ref{p:Qprods}, exercise~\ref{e:Qterms}), and part (b) is left as exercise. For (c), we have the
following isomorphisms:
\begin{align*}
 \CC_Q(A\times B,C) &= \CC(Q(A\times B),C) &&\text{definition of $\CC_Q$}\\
    &\cong\CC(QA\otimes QB,C) &&\text{part (b)}\\
    &\cong\CC(QA,QB\linimpl C) &&\text{monoidal closure of $\CC$}\\
    &= \CC_Q(A,QB\linimpl C) &&\text{defn of $\CC_Q$}.
\end{align*}
Concretely, we obtain $\theta_{A}:\CC_Q(A\times B,C)\lred{\cong}\CC_Q(A,QB\linimpl C)$ by:
\begin{align*}
  \theta_{A}      &:= (\kl[Q]{f}:A\times B\rarr C)   \longmapsto \kl[Q]{(\Lambda(f\circ j^{-1}_{A,B}))}\\
  \theta_{A}^{-1} &:= (\kl[Q]{g}:A\rarr QB\linimpl C)\longmapsto \kl[Q]{(\Lambda^{-1}(g)\circ j_{A,B})} \,.
\end{align*}
Clearly, $\theta_A$ is a bijection. In order to establish couniversality of the exponential, we need also show naturality in $A$ (see exercise~\ref{ex:expon}). So take $\kl[Q]{f}:A\times B\rarr C$ and $\kl[Q]{h}:A'\rarr A$. Note first that the following commutes.
\begin{equation}\tag{$*$}
\xymatrix{Q(A\times B)\ar[r]^{\delta}\ar[d]_{j} & \QQ(A\times B)\ar[rr]^{Q\ang{Q\pi_1,Q\pi_2}} && Q(QA\times QB)\ar[d]^{j}\\
QA\otimes QB\ar[rrr]_{\delta\otimes\delta} &&& \QQ A\otimes\QQ B}
\end{equation}
Note also that, for any $\kl[Q]{{h_i}}:A_i'\rarr A_i$ in $\CC_Q$, $i=1,2$, we have:
\[ \kl[Q]{{h_1}}\times\kl[Q]{{h_2}}:=\kl[Q]{\big(Q(A'_1\times A'_2)\lred{\ang{Q\pi_1,Q\pi_2}}QA'_1\times QA'_2\lred{h_1\times h_2}A_1\times A_2\big)} \]
Thus, noting that $\id[B]^{(\CC_Q)}=\eps_{B.Q}$,
\begin{align*}
 \theta_{A'}(\kl[Q]{f}\circ\kl[Q]{h}\times\id[B]^{(\CC_Q)})
    &= \kl[Q]{\big(\Lambda(f\circ Q(h\times\eps\circ\ang{Q\pi_1,Q\pi_2})\circ\delta\circ j^{-1})\big)} \\
    &= \kl[Q]{\big(\Lambda(f\circ Q(h\times\eps)\circ Q\ang{Q\pi_1,Q\pi_2}\circ\delta\circ j^{-1})\big)} \\
    &\overset{(*)}{=} \kl[Q]{\big(\Lambda(f\circ Q(h\times\eps)\circ j^{-1}\circ\delta\otimes\delta)\big)} \\
    &= \kl[Q]{\big(\Lambda(f\circ j^{-1}\circ Qh\otimes Q\eps\circ\delta\otimes\delta)\big)} \\
    &= \kl[Q]{\big(\Lambda(f\circ j^{-1}\circ (Qh\circ\delta)\otimes\id)\big)} \\
    &= \kl[Q]{\big(\Lambda(f\circ j^{-1})\circ Qh\circ\delta\big)}=\theta_A(\kl[Q]{f})\circ\kl[Q]{h}
\end{align*}
as required. \qed

\begin{myexercise}
Show part (b) of proposition~\ref{p:CCQ}. For the defined $j$, show commutativity of $(*)$.
\end{myexercise}

\subsection{Exercises}
\begin{enumerate}\renewcommand{\theenumi}{\textbf{\arabic{enumi}}}
 \item We say that a category $\CC$ is \emph{well-pointed} if it contains a terminal object $\ena$ and, for any pair of arrows $f,g:A\rarr B$,
 \[ f\neq g \implies \exists h:1\rarr A.\, f\circ h\neq g\circ h\,. \]
 Let now $\CC$ be a well-pointed category with a terminal object $\ena$ and binary coproducts, and consider the functor $G:\CC\rarr\CC$ given by:
 \[
   G:= A \mapsto A+\ena \,,\,
   f \mapsto f+\id[\ena]\,.
 \]
 If $\CC(\ena,\ena+\ena)=\{\incl[1],\incl[2]\}$ with $\incl[1]\neq\incl[2]$, show that if $(G,\eta,\mu)$ is a monad on $\CC$ then, for each object $A$:
 \[ \eta_A=A\lred{\incl[1]}A+\ena\,,\qd\mu_A=(A+\ena)+\ena\lred{[\id[A+\ena],\incl[2]]}A+\ena\,. \]
 \item Let $\CC$ be a SMCC and let $(Q,\eps,\delta)$ be a comonad on $\CC$.
 \begin{enumerate}
 \item[(a)]
 Suppose that the sequents $\WBang A\vdash B$ and $\WBang B\vdash C$ are provable and let $f:QA\rarr B$ and $g:QB\rarr C$ be their interpretations
 (\ie the interpretations of their proofs) in $\CC$. Find the interpretations of the sequent $\WBang A\vdash\WBang C$ which correspond to each of the following proofs and show that the two interpretations are equal.
 \[ \infer[\Cut]{\WBang A\vdash\WBang C}{\infer[\someQWE{\WBang R}]{\WBang A\vdash\WBang B}{\infer{\WBang A\vdash B}{\vdots}}\quad
    \infer[\someQWE{\WBang R}]{\WBang B\vdash\WBang C}{\infer{\WBang B\vdash C}{\vdots}}} \qquad\quad
 \infer[\someQWE{\WBang R}]{\WBang A\vdash\WBang C}{\infer[\Cut]{\WBang A\vdash C}{\infer[\someQWE{\WBang R}]
    {\WBang A\vdash\WBang B}{\infer{\WBang A\vdash B}{\vdots}}\quad
    {\infer{\WBang B\vdash C}{\vdots}}}}
 \]
 \item[(b)]
 Find the interpretations in $\CC$ of the following proofs; are the interpretations equal?
 \[ \infer[\someQWE{\WBang R}]{\WBang\WBang A\vdash\WBang\WBang A}{\infer[\someQWE{\WBang L}]{\WBang\WBang A\vdash\WBang A}
    {\infer[\Ax]{\WBang A\vdash\WBang A}{}}} \qquad\quad
    \infer[\someQWE{\WBang L}]{\WBang\WBang A\vdash\WBang\WBang A}{\infer[\someQWE{\WBang R}]{\WBang A\vdash\WBang\WBang A}
    {\infer[\Ax]{\WBang A\vdash\WBang A}{}}} \]
 \end{enumerate}
 \item Show that a symmetric monoidal category $\CC$ has finite products (given by $\otimes,I$, etc.) iff there are monoidal natural transformations
 \[ d:\Id_\CC\lrar\otimes\circ\ang{\Id_\CC,\Id_\CC}\,,\qd e:\Id_\CC\lrar K_I\,, \]
 such that the following diagram commutes, for any $A\in Ob(\CC)$.
 \[\xymatrix@C=12mm{A\ar[dr]^{d_A} & A\otimes I\ar[l]_{r_A} \\ I\otimes A\ar[u]^{l_A} & A\otimes A\ar[l]^{\id[A]\otimes e_A}\ar[u]_{e_A\otimes\id[A]}}\]
\end{enumerate}


\appendix

\section{Review of Sets, Functions and Relations}
Our aim in this Appendix is to provide a brief review of notions we will assume in the notes. If the first paragraph is not familiar to you, you  will need to acquire more background before being ready to read the notes.

\paragraph{Cartesian products, relations and functions}
Given sets $X$ and $Y$, their cartesian product is
\[ X \times Y = \{ (x, y) \mid x \in X \; \wedge \; y \in Y \} \, . \]
A \emph{relation} $R$ from $X$ to $Y$, written $R : X \rarr Y$, is a subset $R \subseteq X \times Y$. Given such a relation, we write $(x, y) \in R$, or equivalently $R(x, y)$. We  compose relations as follows: if $R : X \rarr Y$ and $S : Y \rarr Z$, then for all $x \in X$ and $z \in Z$:
\[ R; S (x, z) \equiv \exists y \in Y. \, R(x,y) \, \wedge \, S(y,z) \, . \]
A relation $f : X \rarr Y$ is a \emph{function} if it satisfies the following two properties:
\begin{itemize}
\item (single-valuedness): if  $(x, y) \in f$ and $(x, y') \in f$, then $y = y'$.
\item (totality): for all $x \in X$, for some $y \in Y$, $(x, y) \in f$.
\end{itemize}
If $f$ is  a function, we write $f(x) = y$ or $f: x \mapsto y$ for $(x, y) \in f$. Function composition is written as follows: if $f : X \rarr Y$ and $g : Y \rarr Z$,
\[ g \circ f (x) = g(f(x)) \, . \]
It is easily checked that $g \circ f = f;g$, viewing functions as relations.

\paragraph{Equality of functions}
Two functions $f, g : X \rarr Y$ are \emph{equal} if they are equal as relations, \ie as sets of ordered pairs. Equivalently, but more conveniently, we can write:
\[ f = g \;\; \Longleftrightarrow \;\; \forall x \in X. \, f(x) = g(x) \, . \]
The right-to-left implication is the standard tool for proving equality of functions on sets. As we shall see, when we enter the world of category theory, which takes a more general view of ``arrows'' $f : X \rarr Y$, for most purposes we have to leave this familiar tool behind!

\paragraph{Making the arrow notation for functions and relations unambiguous}
Our definitions of functions and relations, as they stand, have an unfortunate ambiguity. Given a relation $R : X \rarr Y$, we cannot uniquely recover its ``domain'' $X$ and ``codomain'' $Y$. In the case of a function, we can recover the domain, because of totality, but not the codomain.

\noindent \textbf{Example} Consider the set of ordered pairs $\{ (n, n) \mid n \in \Nat \}$, where $\Nat$ is the set of natural numbers. Is this the identity function $\id[\Nat] : \Nat \lrarr \Nat$, or the inclusion function $\mathsf{inc} : \Nat \hookrightarrow \ZZ$, where $\ZZ$ is the set of integers?

We wish to have unambiguous notions of domain and codomain for functions, and more generally relations. Thus we modify our official definition of a relation from $X$ to $Y$ to be an ordered triple $(X, R, Y)$, where $R \subseteq X \times Y$. We then define composition of $(X, R, Y)$ and $(Y, S, Z)$ in the obvious fashion, as $(X, R;S, Z)$. We treat functions similarly. We shall not belabour this point in the notes, but it is implicit when we set up perhaps the most fundamental example of a category, namely the category of sets.

\paragraph{Size} We shall avoid explicit discussion of set-theoretical foundations in the text, but we include a few remarks for the interested reader. Occasionally, distinctions of set-theoretic size do matter in category theory. One example which does arise in the notes is when we consider $\Cat$, the category of ``all'' categories. Does this category belong to (is it an object of) itself, at the risk of a Russell-type paradox? The way we avoid this is to impose some set-theoretic limitation of size on the categories gathered into $\Cat$. $\Cat$  will then be too big to fit into itself. For example, we can limit $\Cat$ to those categories whose collections of objects and arrows form \emph{sets} in the sense of some standard set theory such as ZFC. $\Cat$ will then be a proper class, and will not be an object of itself.
One assumption we do make throughout the notes is that the categories we deal with are ``locally small'', \ie that all hom-sets are indeed sets. Another place where some technical caveat would be in order is when we form functor categories. In practice, these issues never (well, hardly ever) cause problems, because of the strongly-typed nature of category theory. We leave the interested reader to delve further into these issues by consulting some of the standard texts.

\section{Guide to Further Reading}
Of the many texts on category theory, we shall only mention a few, which may be particularly useful to someone who has read these notes and wishes to learn more.

The short text~\cite{Pie} is very nicely written and gently paced; it is probably a little easier going than these notes.
A text which is written with a clarity and at a level which makes it ideal as a next step after these notes is~\cite{HS}. A text particularly useful for its large number of exercises with solutions is~\cite{BW}.

Another very nicely written text, focussing on the connections between categories and logic, and especially topos theory, is~\cite{Gol}, recently reissued by Dover Books. A classic text on categorical logic is \cite{LS}. A more advanced text on topos theory is~\cite{MacM}.

The text~\cite{Mac} is a classic by one of the founders of category theory. It assumes considerable background knowledge of mathematics to fully appreciate its wide-ranging examples, but it provides invaluable coverage of the key topics.

A stimulating text on the correspondence between computation and logic is~\cite{PT}; it  is out of print, but available online. A more recent  text on this topic is~\cite{CH}.

The 3-volume handbook~\cite{Bor} provides coverage of a broad range of topics in category theory.
The book~\cite{LS2} is somewhat idiosyncratic in style, but offers insights by one of the key contributors to category theory.

\bibliographystyle{plain}

\end{document}

%% file: LNPpreamble.tex
\usepackage[arrow,matrix,curve,frame]{xy}
\usepackage{amsmath,amssymb,amsfonts,latexsym,stmaryrd}
\usepackage{ifthen}
\usepackage{pstricks,pst-node,pst-tree}
\usepackage{proof,array}
\usepackage{caption,theorem}

\newcounter{newl_i}
\newcounter{quad_i}

\newtheorem{myproposition}[theorem]{Proposition}
\newtheorem{mycorollary}[theorem]{Corollary}
\newtheorem{fact}[theorem]{Fact}
\theorembodyfont{\rmfamily}
\newtheorem{myexample}[theorem]{Example}
\newtheorem{myexercise}[theorem]{Exercise}
\newtheorem{mydefinition}[theorem]{Definition}
\newtheorem{notation}[theorem]{Notation}
\newtheorem{myremark}[theorem]{Remark}

\newcommand{\ie}{\textit{i.e.}~}
\newcommand{\tbf}[1]{\textbf{#1}}
\newcommand{\Ob}{\mathsf{Ob}}
\newcommand{\Arr}{\mathsf{Ar}}
\newcommand{\dm}{\mathsf{dom}}
\newcommand{\cod}{\mathsf{cod}}
\newcommand{\labarrow}[1]{\stackrel{#1}{\longrightarrow}}
\newcommand{\Nat}{\mathbb{N}}

\newcommand{\myintertext}[1]{\\[-5pt]\intertext{#1}\\[-20pt]}
\newcommand{\ppi}{\boldsymbol\pi}
\newcommand{\actn}[1]{\mathop{#1\,_{^\bullet}}}
\newcommand{\boldemph}[1]{\textbf{\emph{#1}}}
\newcommand{\Ga}{\Gamma}
\newcommand{\App}[1][]{\mathsf{ev}_{#1}}
\newcommand{\Curry}[1]{\Lambda(#1)}

\newcommand{\sw}[3]{\actn{\boldsymbol{(}#2\ #3\boldsymbol{)}}#1}
\newcommand{\fv}{\textsf{fv}}
\newcommand{\objk}[1]{\widetilde{#1}}
\newcommand{\CCl}{\CC_\lambda}
\newcommand{\trt}{\,\centerdot\,}
\newcommand{\ZZ}{\mathbb{Z}}
\newcommand{\RR}{\mathbb{R}}
\newcommand{\pow}{\PP}
\newcommand{\Term}{\boldsymbol1}

\newcommand{\impl}{\Rightarrow}
\newcommand{\trn}[1]{\lsem#1\rsem}

\newcommand{\qd}[1][]{\ifthenelse{\equal{#1}{}}{\quad}{\setcounter{quad_i}{0}\whiledo{\value{quad_i}<#1}{\quad\stepcounter{quad_i}}}}
\newcommand{\nada}{$\text{}$}
\newcommand{\ang}[1]{\langle#1\rangle}
\newcommand{\incl}[1][]{\mspace{1.5mu}\mathsf{in}_#1}
\renewcommand{\proof}[1][]{\noindent\textbf{Proof#1:} }
\renewcommand{\qedsymbol}{\scriptsize$\blacksquare$}
\newcommand{\deqsymbol}{$\blacktriangle$}
\renewcommand{\qed}[1][0]{\ifthenelse{#1<0}{\nada\\[-1.9\baselineskip]\nada\hfill\makebox[0mm][r]{\rm\qedsymbol}\newlinesminusone{-#1}}%
                                                  {\nada\hfill\makebox[0mm][r]{\rm\qedsymbol}\newlines{#1}}}
\newcommand{\deq}[1][0]{\ifthenelse{#1<0}{\nada\\[-1.9\baselineskip]\nada\hfill\makebox[.75ex][r]{\rm\deqsymbol}\newlinesminusone{-#1}}%
                                                  {\nada\hfill\makebox[0mm][r]{\rm\deqsymbol}\newlines{#1}}}
\newcommand{\newlines}[1]{\setcounter{newl_i}{0}\whiledo{\value{newl_i}<#1}{\nada\\\stepcounter{newl_i}}}
\newcommand{\newlinesminusone}[1]{\setcounter{newl_i}{1}\whiledo{\value{newl_i}<#1}{\nada\\\stepcounter{newl_i}}}

\DeclareMathAlphabet{\mathbbold}{U}{bbold}{m}{n}
\renewcommand{\1}{\mathbbold{1}}
\renewcommand{\2}{\mathbbold{2}}
\newcommand{\TWO}{\2_{\scriptscriptstyle\rightrightarrows}}
\newcommand{\R}{\mathbb{R}}
\newcommand{\keno}{\emptyset}
\newcommand{\lrar}{\lrarr}
\newcommand{\lred}[1]{\xrightarrow{#1}}
\newcommand{\brar}{\longrightarrow_\beta}
\newcommand{\ena}{\mathbf{1}}
\newcommand{\zero}{\mathbf{0}}
\newcommand{\Ax}{\someQWE{Id}}
\newcommand{\someQWE}[1]{\text{\small$\mathsf{#1}$}}
\newcommand{\someQWEqwe}[2]{\someQWE{#1\,#2}}
\newcommand{\AndI}{\someQWEqwe{\wedge}{I}}
\newcommand{\AndEl}{\someQWEqwe{\wedge}{E_1}}
\newcommand{\AndEr}{\someQWEqwe{\wedge}{E_2}}
\newcommand{\ImpI}{\someQWEqwe{\mathnormal\supset}{I}}
\newcommand{\ImpE}{\someQWEqwe{\mathnormal\supset}{E}}
\newcommand{\AndL}{\someQWEqwe{\wedge}{L}}
\newcommand{\AndR}{\someQWEqwe{\wedge}{R}}
\newcommand{\ImpL}{\someQWEqwe{\mathnormal\supset}{L}}
\newcommand{\ImpR}{\someQWEqwe{\mathnormal\supset}{R}}
\newcommand{\TenL}{\someQWEqwe{\otimes}{L}}
\newcommand{\TenR}{\someQWEqwe{\otimes}{R}}
\newcommand{\LimpL}{\someQWEqwe{\mathnormal\linimpl}{L}}
\newcommand{\LimpR}{\someQWEqwe{\mathnormal\linimpl}{R}}
\newcommand{\LimpE}{\someQWEqwe{\mathnormal\linimpl}{E}}
\newcommand{\WithL}{\someQWEqwe{\&}{L}}
\newcommand{\WithR}{\someQWEqwe{\&}{R}}
\newcommand{\Cut}{\someQWE{Cut}}
\newcommand{\Exchs}{\someQWE{Exch}}
\newcommand{\Cons}{\someQWE{Contr}}
\newcommand{\Weaks}{\someQWE{Weak}}
\newcommand{\Exch}{\someQWE{Exch}}
\newcommand{\Con}{\someQWE{Contr}}
\newcommand{\Weak}{\someQWE{Weak}}

\newcommand{\lrarr}{\longrightarrow}
\newcommand{\llarr}{\longleftarrow}
\newcommand{\rarr}{\rightarrow}
\newcommand{\uscore}{\_\!\_}

\newcommand{\CC}{\mathcal{C}}
\newcommand{\PP}{\mathcal{P}}
\newcommand{\DD}{\mathcal{D}}
\newcommand{\II}{\mathcal{I}}
\newcommand{\EE}{\mathcal{E}}
\newcommand{\id}[1][]{\mspace{1.5mu}\mathsf{id}_{#1}}
\newcommand{\Id}{\mathsf{Id}}
\newcommand{\ev}[2]{\mathsf{ev}_{#1,#2}}
\newcommand{\evv}{\mathsf{ev}}
\newcommand{\op}[1]{#1^{\mathsf{op}}}
\newcommand{\CCop}{\op{\CC}}
\newcommand{\Set}{\textbf{Set}}
\newcommand{\Cat}{\textbf{Cat}}
\newcommand{\Vect}{\textbf{Vect}}
\newcommand{\Mon}{\textbf{Mon}}
\newcommand{\Top}{\textbf{Top}}
\newcommand{\Grp}{\textbf{Grp}}
\newcommand{\Rel}{\textbf{Rel}}
\newcommand{\Pos}{\textbf{Pos}}
\newcommand{\iso}{\stackrel{\cong}{\longrightarrow}}
\newcommand{\List}{\mathsf{List}}
\newcommand{\MList}{\mathsf{MList}}

\newcommand{\eps}{\varepsilon}
\newcommand{\QQ}{Q^2\!}
\newcommand{\QQQ}{Q^3\!}
\newcommand{\TT}{T^2\!}
\newcommand{\TTT}{T^3\!}
\newcommand{\YY}{\mathbf{Y}}
\newcommand{\funsp}[2]{(#1 \Rightarrow #2)}

\newcommand{\lldots}{, \ldots ,}
\newcommand{\lsem}{\llbracket}
\newcommand{\rsem}{\rrbracket}

\newcommand{\linimpl}{\multimap}
\newcommand{\With}{\mathbin{\&}}
\newcommand{\olet}[4]{\mathsf{let} \; #1 \; \mathsf{be} \; #2 \otimes #3 \; \mathsf{in} \; #4}

\newcommand{\kl}[2][T]{#2_{.#1}}
\newcommand{\Bang}{\boldsymbol!\mspace{1mu}}
\newcommand{\WBang}{\boldsymbol{\hat{!}}\mspace{1mu}}
\newcommand{\HY}{\,---\,}
    {\end{list}}
\newcommand{\cc}[3][5pt]{\ncbar[armA=#1,armB=#1,angle=-90,nodesep=3pt,linewidth=0.5pt]{#2}{#3}}

%% file: LNPnotes.bbl
\begin{thebibliography}{99}


\bibitem{BW}
M.~Barr and C.~Wells, \textsl{Category Theory for Computing Science}, Third Edition, Publications CRM, 1999.

\bibitem{Bor}
F.~Borceux, \textsl{Handbook of Categorical Algebra Volumes 1--3}, Cambridge University Press, 1994.

\bibitem{PT}
J.-Y.~Girard. P.~Taylor and Y.~Lafont, \textsl{Proofs and Types}, Cambridge University Press, 1989.

\bibitem{Gol}
R.~Goldblatt, \textsl{Topoi, the Categorial Analysis of Logic}, North-Holland 1984. Reprinted by Dover Books, 2006.

\bibitem{HS}
H.~Herrlich and G.~Strecker, \textsl{Category Theory}, Third Edition, Heldermann Verlag, 2007.

\bibitem{LS}
J.~Lambek and P.~J.~Scott, \textsl{Introduction to Higher-Order Categorical Logic}, Cambridge University Press, 1986.

\bibitem{LS2}
F.~Lawvere and S.~Schanuel, \textsl{Conceptual Mathematics: A First Introduction to Categories}, Cambridge University Press, 1997.

\bibitem{Mac}
S.~Mac Lane, \textsl{Categories for the Working Mathematician}, Second Edition, Springer, 1998.

\bibitem{MacM}
S.~Mac Lane and I.~Moerdijk, \textsl{Sheaves in Geometry and Logic: A First Introduction to Topos Theory}, Springer, 1994.

\bibitem{Pie}
B.~Pierce, \textsl{Basic Category Theory for Computer Scientists}, MIT Press, 1991.

\bibitem{CH}
M.~H.~S{\o}rensen and P.~Urzyczyn, \textsl{Lectures on the Curry-Howard Isomorphism}, Elsevier, 2006.

\end{thebibliography}
